\documentclass[preprint, 11pt]{imsart}

\RequirePackage[OT1]{fontenc}
\RequirePackage{amsthm,amsmath}
\RequirePackage[authoryear,round]{natbib}
\RequirePackage[colorlinks,citecolor=blue,urlcolor=blue,pdfstartview=FitB]{hyperref}

\usepackage{filecontents}
\usepackage{bibentry}

\usepackage{titlesec} 
\titleformat{\subsection}[runin]
{\itshape\large\bfseries}{\thesubsection}{1em}{}

\usepackage{tikz}
\usetikzlibrary{matrix}
\usepackage{multicol}
\usepackage{multirow}
\usepackage{changepage}
\usepackage{float}
\usepackage{amsmath,amssymb,epsfig}
\usepackage{graphics}
\usepackage{tabularx}
\usepackage{framed}
\usepackage[mathscr]{euscript}
\usepackage{imsart}

\usepackage{epstopdf}

\usepackage{stmaryrd}

\usepackage{titlesec}
\titleformat*{\section}{\large\bfseries}

\usepackage{lipsum}

\usepackage{graphicx}
\usepackage{bbm}

\startlocaldefs
\numberwithin{equation}{section}
\theoremstyle{plain}

\long\def\comment#1{}

\newtheorem{theorem}{Theorem}

\newtheorem{corollary}{Corollary}

\newtheorem{proposition}{Proposition}
\newtheorem{assumption}{Condition \!\!}
\theoremstyle{definition}

\newtheorem{algorithmsimple}{Algorithm}

\newcommand{\eps}{\varepsilon}

\newcommand{\be}{\begin{eqnarray}}
\newcommand{\ee}{\end{eqnarray}}

\newcommand{\ba}{\begin{array}}
\newcommand{\ea}{\end{array}}
\newcommand{\bs}{\begin{align}\begin{split}\nonumber}
\newcommand{\bsnumber}{\begin{align}\begin{split}}
\newcommand{\es}{\end{split}\end{align}}

\newcommand{\n}{n}

\renewcommand{\(}{\left(}
\renewcommand{\)}{\right)}
\renewcommand{\[}{\left[}
\renewcommand{\]}{\right]}
\renewcommand{\hat}{\widehat}

\newcommand{\G}{k}

\newcommand{\Ep}{{\mathrm{E}}}

\newcommand{\CRS}{\cite{Canay2017}}
\newcommand{\IM}{\cite{Ibragimov2010}}

\newcommand{\Test}{\mathtt{Test}}
\newcommand{\Cluster}{\mathtt{Cluster}}

\newcommand{\Ev}{{\mathcal{E}}}

\renewcommand{\Pr}{{\mathrm{P}}}

\renewcommand{\hat}{\widehat}
\renewcommand{\leq}{\leqslant}
\renewcommand{\geq}{\geqslant}

\newcommand{\BCH}{\cite{BCH-inference}}

\usepackage[flushleft]{threeparttable}
\usepackage{multirow}
\usepackage{algorithm,algorithmic}
\usepackage{mathrsfs}
\usepackage{enumitem}   
\usepackage{subfig}
\usepackage{pdflscape}

\renewcommand{\diamond}{b}

\endlocaldefs

\begin{document}

\begin{frontmatter}
\title{\Large Inference for Dependent Data with Learned Clusters 
}
\runtitle{\ \ Inference with Learned Clusters}

\runauthor{J. Cao, C. Hansen, D. Kozbur and L. Villacorta \ \ }

\begin{multicols}{2}

\author{\fnms{Jianfei} \snm{Cao}\ead[label=e1]{j.cao@northeastern.edu}}



\thankstext{t1}{First version:  October 2018.  This version is of  \today.   Christian Hansen would like to thank the National Science Foundation as well as The University of Chicago Booth School of Business for financial support of this research.  Damian Kozbur would like to thank The University of Z\"urich for financial support of this research.}
\affiliation{Chicago Booth }
\address{Northeastern University \\ \printead{e1}
}

\author{\fnms{Damian} \snm{Kozbur}\ead[label=e3]{damian.kozbur@econ.uzh.ch}}
\address{University of Zurich \\ \printead{e3}
}

\author{\fnms{Christian} \snm{Hansen}\ead[label=e2]{chansen1@chicagobooth.edu}}
\affiliation{Chicago Booth }
\address{University of Chicago Booth School of Business \\ \printead{e2}
}

\author{\fnms{Lucciano} \snm{Villacorta}\ead[label=e4]{ lvillacorta@bcentral.cl}}
\affiliation{Central Bank of Chile}
\address{Central Bank of Chile \\ \printead{e4}
}

\end{multicols}

\begin{abstract}  This paper presents and analyzes an approach to cluster-based inference for dependent data.  The primary setting considered here is with spatially indexed data in which the dependence structure of observed random variables is characterized by a known, observed dissimilarity measure over spatial indices.  Observations are partitioned into clusters with the use of an unsupervised clustering algorithm applied to the dissimilarity measure.  Once the partition into clusters is learned, a cluster-based inference procedure is applied to a statistical hypothesis testing procedure. The procedure proposed in the paper allows the number of clusters to depend on the data, which gives researchers a principled method for choosing an appropriate clustering level.  The paper gives conditions under which the proposed procedure asymptotically attains correct size. A simulation study shows that the proposed procedure attains near nominal size in finite samples in a variety of statistical testing problems with dependent data.
 \\
	

\end{abstract}


\end{frontmatter}


\section{Introduction}
\label{intro}

Conducting accurate statistical inference with data featuring serial or cross-sectional dependence requires carefully accounting for the underlying dependence structure. A variety of methods are available for researchers analyzing dependent data.  An important class of inferential methods is the class of \textit{cluster-based} methods. Cluster-based methods work with a partition of observations into clusters. Inference proceeds by treating as negligible any covariance between observations that fall in different clusters, followed by performing appropriate tests for a statistical hypothesis of interest. Cluster-based methods deliver asymptotically valid inference, in the sense of controlling size of tests and coverage of interval estimates, in a variety of dependent-data settings; see, e.g., \cite{Ibragimov2010}, \cite{BCH-inference}, \cite{Canay2017}, and \cite{HansenLee:cluster}.


\textcolor{black}{In some settings, there is a natural, ``known'' way to partition data in clusters to appropriately account for dependence. A common setting in empirical economics is an $n \times T$ panel consisting of a random sample of $i=1,...,n$ individuals followed over $T$ time periods. In this case, partitioning data into clusters consisting of all time periods for each individual is justified by the assumption that observations $i$ and $j$ are independent whenever $i \neq j$ which is motivated by random sampling. This case has been analyzed extensively in the literature; see, e.g., \cite{wooldridge:text}. Natural partitions also exist in group randomized trial settings where the groups at which randomization occur provide appropriate clusters; see, e.g., \cite{abadie2017should}.}

However, there is not an immediate natural set of partitions available in many settings. Rather, researchers often use data where observations can be thought of as being located in some space and where a measure of distance or dissimilarity between observations is available. For example, researchers routinely use data where observations are indexed geographically over an irregular geographic region and spatial correlation is a concern. As another example, researchers may be working with data in which observations carry a notion of economic distance. Examples of such applications include \cite{Conley:Dupor:Sectoral}, which constructs a notion of economic distance based on input-output tables, and  \cite{conley:topa:socio-econ-distance}, which constructs a notion of economic distance based on demographic composition of zip codes. In such cases, researchers may wish to use cluster-based inference as a robust method to account for dependence in the data but desire some guidance about exactly which observations to cluster together and how many clusters to use.

In this paper, we consider settings in which natural partitions may not exist but a distance or dissimilarity measure that is informative about the underlying correlation structure is available. \textcolor{black}{Throughout, we only consider cases where observation indexes are non-random or conditioned upon. Considering only such cases will rule out some examples such as endogenous network formation models or models with location choice that may be natural to consider in some cases where dependence is based on a notion of economic distance.} We assume that dependence between observable quantities is suitably small when the distance between them is suitably large. A main goal is then to verify that cluster-based methods remain asymptotically valid when clusters are not pre-specified but are constructed in a data-dependent fashion using the underlying notion of dissimilarity and `unsupervised learning' methods -- e.g. $k$-\texttt{medoids}, $k$-\texttt{means}, and \texttt{hierarchical clustering}.

At a high level, the data-driven cluster-based inference approach we propose for testing a hypothesis $H_0$ proceeds in three main steps. In the first step, unsupervised learning methods are applied using the observed dissimilarity between observations to produce a collection of partitions of the observations into clusters, $\mathscr C = \{\mathcal{C}^{(k)}\}_{k=2}^{k_{\max}}$ for some integer $k_{\max}$. Each element $\mathcal{C}^{(k)}$ provides a candidate structure for use with a cluster-based method and may be based on a different number of clusters. Cluster-based methods for testing $H_0$ will involve decision thresholds for determining whether to reject $H_0$. For example, we may decide to reject $H_0$ when an empirical p-value is smaller than some threshold. In the second step, we select a partition $\hat{\mathcal{C}}$ from $\mathscr{C}$ and a rejection rule for $H_0$ chosen to maximize exact weighted average power against a collection of alternatives subject to maintaining exact size control within completely specified data generating processes (DGPs) for the observed data. These DGPs may be data dependent; for example, they may involve parameters whose values are estimated in the data. While $\hat{\mathcal{C}}$ and the rejection rule are selected based on size and power calculations in pre-specified DGPs, the selection is done in such a way that asymptotic size control is formally maintained regardless of whether these DGPs are correctly specified. In the final step, inference in the actual data is performed using the rejection rule and $\hat{\mathcal{C}}$ determined in the second step. As $\mathscr{C}$ will typically contain candidate partitions with different numbers of groups, the procedure outlined above encompasses a data-driven method for choosing the appropriate level at which to cluster the data as well as a data-driven way to assign observations to clusters.

The main theoretical results in this paper provide conditions under which inference based on the general procedure outlined in either \cite{Ibragimov2010} or \cite{Canay2017} remains valid when data-dependent clusters produced by applying $k$-\texttt{medoids} are used in lieu of relying on a fixed number of pre-specified clusters. The regularity conditions involve moment and mixing restrictions, weak homogeneity assumptions on second moments of regressors and unobservables across groups, and restrictions on group boundaries. These moment and mixing conditions are implied by routine assumptions necessary for use of central limit approximations and the required homogeneity is less restrictive than covariance stationarity. The assumptions allow relatively general heterogeneous dependent processes and are akin to those imposed in other heteroskedasticity and autocorrelation robust (HAR) inference approaches; e.g., \cite{andrews1991} and \cite{KelejianPrucha:SpatialHAC}. \textcolor{black}{It is important to note that the formal results do not cover the more conventional approach based on using ``clustered'' standard errors. Validity of inference based on clustered standard errors in our setting with general weak dependence accommodated by considering a fixed number of clusters would rely on homogeneity conditions, including that clusters are asymptotically identical in size, that are unlikely to be satisfied with clusters produced by unsupervised learning methods.}

While providing the main inferential results, we also develop novel insights about the finite-sample behavior of the $k$-\texttt{medoids} clustering algorithm that are important for cluster-based inference. We verify that  $k$-\texttt{medoids}  partitions satisfy suitably defined balance and small boundary conditions without requiring any notion of consistency to a ``true'' partition. Balance and small boundary properties are important if $k$-\texttt{medoids} groups are to be input into cluster-based inference, as general dependent data settings do not necessitate existence of ``true'' partitions, and previous results for cluster-based inference with pre-specified clusters rely on small cluster boundaries.
A key ingredient to establishing this result is that locations are relatively regularly distributed in space.  Our formal conditions rule out small regions containing a large fraction of observations with remaining observations spread sparsely over the remaining space. This could be problematic in data on individuals taken from a metropolitan area with relatively small geographic area but a large population and a surrounding rural area with large geographic area but small population, when the relevant measure of dependence between individuals depends on geographic proximity. 
These properties of $k$-\texttt{medoids} clustering are, to our knowledge, new to the unsupervised learning literature.

As cluster-based methods provide a general approach to performing robust inference in the presence of dependent and heteroskedastic data, they belong to the general class of HAR inference procedures; see, e.g., \cite{Lazarus2018} for a recent review centered on time series analysis and \cite{conley:gmmhac} and \cite{KelejianPrucha:SpatialHAC} for seminal references in the spatial context. There is also recent work in combined contexts with both serial and spatial dependence; for instance \cite{bai2020standard} consider weak serial dependence and sparse spatial dependence, and propose a new standard error estimator using  Newey-West HAC along the time dimension and using the thresholding approach in \cite{bickel2008covariance} to estimate the small subgroup of $\left(i,j\right)$ pairs that have nonzero spatial dependence. A similar procedure is also proposed by \cite{Cai2021}.   

Our results cover scenarios where the number of clusters used for inference is bounded as in \cite{BCH-inference}, \cite{Ibragimov2010}, and \cite{Canay2017} and thus are in the spirit of the ``fixed-b'' approach - e.g. \cite{KVB}, Kiefer and Vogelsang (\citeyear{kiefvog1}, \citeyear{kiefvog2}), \cite{bchv:fbshac} - and the projection approach with fixed number of projections - e.g. \cite{phillips:HACautomated}, \cite{muller:robustLRV}, \cite{sun:orthonormal} - both of which are well-known to lead to improved size control relative to more traditional approaches under weak dependence. 

Within the general context of HAR inference, the choice of cluster structure is analogous to the choice of smoothing parameters in other HAR procedures, and our work thus complements the broad literature on data-driven tuning parameter choice in HAR inference. Our proposal is closely related to \cite{Sun2015}, \cite{Lazarus2018}, \cite{LLS:SizePower}, and especially \cite{MeullerWatson:2020}. \cite{Sun2015} proposes a method to select the number of projections to use for HAR inference in a spatially dependent context by minimizing asymptotic mean-square error (MSE) of the variance estimator where terms that enter the expression for the asymptotic MSE are estimated using a parametric model. \cite{LLS:SizePower} provide the size-power frontier based on the ``fixed-b'' approach in a time series context, and \cite{Lazarus2018} explicitly considers tuning parameter choice by trading off size and power along this frontier. \cite{MeullerWatson:2020} considers a novel standard error estimator for use with spatially dependent data based on spatial principal components and proposes a method to choose both the number of components to use in the method and the critical value to use in hypothesis testing based on minimizing confidence interval length subject to exactly controlling size within benchmark models. Their procedure exactly controls size in the benchmark models and asymptotically controls size under general conditions. Relative to these approaches, we focus on data-dependent construction of clusters for use with cluster-based inference methods. Similar to \cite{MeullerWatson:2020}, we choose both a partition and rejection rule for a hypothesis test to maximize weighted average power subject to exact size control within a parametric model in such a way that asymptotic size control is maintained. 



\section{Methodology: Inference with Unsupervised Cluster Learning}
\label{method}

Consider data given by $\mathscr D = \{ \zeta_i \}_{i \in \mathsf X}$.  Here, $\zeta_i$ are observable random variables or vectors and $\mathsf X$ is a (spatial) indexing set of cardinality $n$.  This paper assumes that $\mathsf X$ is equipped with a known dissimilarity measure $\mathrm d$, which is an $n \times n$ array of nonnegative real dissimilarities.  When added emphasis is helpful, $\mathsf X$ is written $(\mathsf X,\mathrm d)$.  The data $\mathscr D$ is distributed according to an unknown (joint) data generating process (DGP) -- $\mathscr D \sim \Pr_0$.  The object $\mathsf X$ will be the main object used to characterize any dependence in the data $\mathscr D$ over $i$.  

Consider testing a scalar null hypothesis, $H_0: \theta_0 = \theta^{\circledast}$, at level $\alpha \in (0,1)$.  Here $\theta^{\circledast}$ is a hypothesized value of a parameter that reflects the data generating process $\Pr_0$.  Our focus will be on the problem of testing a hypothesis about a coefficient in a linear regression model (e.g. testing $H_0:\theta_0=0$ where $\theta_0$ is a parameter in a linear regression, with $\zeta_i = (y_i, x_i)$ or $\zeta_i = (y_i, x_i, z_i)$ representing observations of an outcome variable $Y$, possibly endogenous regressor $X$, and exogenous instruments $Z$.) In the case of this example, $\Pr_0$ will be allowed to be a DGP in which observations are correlated with each other.  Failure to account for dependence in $\mathscr D$ across $i \in \mathsf X$ may lead to substantial size distortion when testing $H_0$. 

Let $\mathcal C = \{ \mathsf C_1,...,\mathsf C_\G \}$ be a partition of $\mathsf X$ of cardinality $\G \geq 2$.  The elements $\mathsf C_1,...,\mathsf C_\G$ are referred to as clusters.  A \textit{cluster-based} inferential procedure for testing $H_0$ is a (possibly random) assignment 
\begin{align} \mathtt{Test}: (\mathscr D, \mathcal C) \mapsto T \in \big \{\text{Fail to Reject, Reject} \big \}. 
\end{align}
Here, the decision rule itself is called $\mathtt{Test}$ and will generally depend on the level $\alpha \in (0,1)$ of the test.  The outcome of the test $\Test$ is referred to as $T$; i.e. $T = \Test(\mathscr D, \mathcal C)$.  The set containing the pair $(\mathscr D, \mathcal C)$ remains unnamed to avoid additional notation.

We consider the following three cluster-based inferential procedures for testing a scalar hypothesis in this paper: the procedure of \cite{Ibragimov2010} (IM), the procedure of \cite{Canay2017} (CRS), and inference based on the cluster covariance estimator as described in \cite{BCH-inference} (CCE).\footnote{Extension of formal results for testing joint hypotheses using CRS is straightforward.} For clarity, consider `$t$-statistic' based testing of the hypothesis $H_0: \theta_0 = \theta^\circledast$. Let $\mathsf C\subseteq \mathsf X$ and $\hat \theta_{\mathsf C}$ be an estimator of $\theta_0$ using only data corresponding to observations in $\mathsf C$. Now define $S \in \mathbb R^{\mathcal C}$ (where $\mathbb R^{\mathcal C}$ denotes functions $\mathcal C \rightarrow \mathbb R$) and the $t$-statistic function $t:\mathbb{R}^{\mathcal C}\rightarrow \mathbb{R}$ such that 
\begin{align}
&S = (S_\mathsf C)_{\mathsf C \in \mathcal C}, \ \ S_{\mathsf C} = (n/\G)^{1/2} (\hat \theta_{\mathsf C} - \theta^\circledast), \phantom{\Bigg |}\\ &t(S)=\frac{ \G^{-1/2}\sum_{\mathsf C \in \mathcal C} S_{\mathsf C}}{\sqrt{(\G-1)^{-1}\sum_{\mathsf C \in \mathcal C} \left (S_{\mathsf C}-\G^{-1}\sum_{\mathsf C' \in \mathcal C} S_{\mathsf C'}\right )^2}}.
\end{align} 
For a specified level $a\in(0,1)$ (which may differ from $\alpha$), there are IM, CRS, and CCE tests, denoted by $\mathtt{Test}_{\text{IM}(a)}$, $\mathtt{Test}_{\text{CRS}(a)}$, $\mathtt{Test}_{\text{CCE}(a)}$.  These tests are defined by their outcomes given data $\mathscr D$ and a partition $\mathcal C$:
\begin{align} 
	&T_{\text{IM}(a)} \ \ = \ \text{Reject}  \ \ \ \text{if} \ \ \ \left | t(S) \right | >t_{1-a/2,\G-1}, \\
	&T_{\text{CRS}(a)} = \ \text{Reject} \ \ \text{ if } \ \ \left | t(S) \right | >  \mathrm{quantile}_{1-a} ( \{   | t(hS)  | \rbrace_{h\in \mathcal H_{\mathcal C}} ),\\
	&T_{\text{CCE}(a)} = \ \text{Reject} \ \ \text{ if } \ \ \left | \frac{ \hat \theta_{\mathsf X}-\theta^\circledast}{ \hat V_{\text{CCE}, \mathcal C}^{1/2}} \right | > \sqrt{\frac{\G}{\G-1}} \times t_{1-a/2,\G-1},
\end{align}
where $t_{1-a/2,\G-1}$ is the $(1-a/2)$-quantile of a $t$-distribution with $\G-1$ degrees of freedom; the set $\{ h S \}_{h \in \mathcal H_{\mathcal C}}$ is the orbit of the action of $\{ \pm 1\}^{\mathcal C}$ on $S$, so that for each $h$,  $hS \in \mathbb R^{\mathcal C}$ has $\mathsf C^{\text{th}}$ component 
$\pm (n/\G)^{1/2}(\hat \theta_{ \mathsf C} - \theta^\circledast)$ for some sign in $\{\pm 1\}$; and $\hat V_{\text{CCE},\mathcal C}$ is the standard cluster covariance matrix estimator. $\Test_{\bullet(a)}$ and $T_{\bullet(a)}$ are used when the choice of IM, CRS, or CCE is unspecified. When we wish to be clear about explicit dependence on $\mathscr D$ and $\mathcal C$ we use the more cumbersome notation $T_{\bullet(a), \mathcal C} = \Test_{\bullet(a)}(\mathscr D, \mathcal C)$.
With prespecified, non-data-dependent clusters, each of the IM, CRS, and CCE procedures controls size asymptotically under weak dependence between observations under respective regularity conditions. 

The second important definition is that of an \textit{unsupervised clustering algorithm}, which is an assignment that returns, to every $\mathsf X = (\mathsf X,\mathrm d)$, a partition of $\mathsf X$ given by the mapping 
\begin{align} \mathtt{Cluster}: \mathsf X \mapsto \mathcal C \end{align} 
generated by trying to keep the distance between observations in the same partition small and the distance between observations in different partitions large. Then, if the dissimilarity $\mathrm d$ appropriately reflects the dependence in $\zeta_i$, the resulting partition $\mathcal C$ may have the desired property that averages of observations belonging to different clusters exhibit negligible dependence. In the formal analysis in Section \ref{clt_section}, the imposed mixing conditions imply that dependence between $\zeta_i$ and $\zeta_j$ vanishes as $\mathrm{d}(i,j)$ becomes large. Then, if $\mathcal C$ tends to place distant observations (as defined by $\mathrm d$) in different clusters, favorable properties of the test $T$ may be anticipated. 

Though there are many commonly used unsupervised clustering algorithms and we expect most to be usable as methods for forming data-dependent clusters, we consider only $k$-\texttt{medoids} in this paper for technical reasons discussed in Section \ref{balanced_small_boundary}. Note that by composition of specific $\Test$ and $\Cluster$ procedures, it is already possible to define an outcome $T$ for $H_0$ given $\mathscr D,\mathsf X$ by constructing $T = \Test(\mathscr D, \Cluster(\mathsf X)).$

The $k$-\texttt{medoids} algorithm we use for establishing our results is as follows.  For finite $(\mathsf X, \mathrm d)$ and medoids $ \mathscr I\subseteq \mathsf X$ define
$\text{cost} (\mathscr I) =\sum_{j\in \mathsf X}\min_{i \in \mathscr I} \mathrm d(i,j)^2.$

\vspace{3mm}

\noindent{$k$-\texttt{medoids}:  \textit{Input.}  $(\mathsf X, \mathrm d)$, $\G$.

\vspace{-2mm}

 \textit{Initialize} a set of medoids $\mathscr I \subseteq \mathsf X$ with cardinality $|\mathscr I | = k$.  	 
\vspace{-1mm}

 \textit{While}  $\text{cost}((\mathscr I \cup \{ j\} )\setminus \{ i\} ) <\text{cost}(\mathscr I)$ for some $i \in \mathscr I$ and $j \in \mathsf X \setminus \mathscr I$,
\vspace{-1mm}

  \textit{Replace} $\mathscr I$ with $\mathscr I \cup \{ \hat j   \} \setminus \{ \hat i   \}$ where $(\hat i, \hat j) \in \arg {\min_{(i,j)\in\mathscr I \times (\mathsf X \setminus \mathscr I)}} \text{cost}((\mathscr I \cup \{ j\} )\setminus \{ i\})$.
\vspace{-1mm}

\textit{Output.} $\mathcal C = \{\mathsf C_{i}\}_{i \in \mathscr I}$ with $j \in \mathsf C_i$ if $\mathrm d(i,j) = \min( \{\mathrm d(i',j) \}_{i' \in \mathscr I})$ .

\vspace{3mm}

This implementation has run time $O(k(n-k)^2)$ per iteration.  Additional details about computational run times in our application and simulations are presented in supplemental material.

The final layer to our proposed testing procedure is a method for data-dependent choice of the cluster-based inferential procedure by considering a collection of candidate testing and clustering procedures of the form $\Test$ and $\Cluster$. We propose making this choice on the basis of simultaneously controlling Type-I and Type-II error rates. Let $\text{Err}_{\text{Type-I}} (\mathtt{Test}, \Cluster) $ denote type-I error for the testing outcome defined by $\mathtt{Test}(\mathscr D, \Cluster(\mathsf X))$.  Next, consider a set of alternatives $\mathbf \Theta_{\text{alt}}$. 
Let $\text{Err}_{\text{Type-II}}(\Test,\Cluster)$ denote a weighted average type-II error against the alternatives in $\mathbf \Theta_{\text{alt}}$. \textcolor{black}{The choice of the alternative set and weighting function will be application specific and should depend on details of the problem. In the empirical illustration in Section \ref{empirical_application}, we chose simple average power over an equally spaced grid of values that we believe encompass all remotely plausible values for the parameter of interest. We believe this practice provides a simple default.} 
Because $\text{Err}_{\text{Type-I}} (\Test, \Cluster) $ and $\text{Err}_{\text{Type-II}  }(\Test,\Cluster)$ will typically not be known, we consider a setting in which estimates $\widehat{\text{Err}}_{\text{Type-I}}(\Test,\Cluster)$ and $\widehat{\text{Err}}_{\text{Type-II}}(\Test,\Cluster) $ are available. 

To finish the final layer, let $\mathscr T$ be a collection of pairs of the form $(\Test, \Cluster)$.  Note that the components $\mathtt{Test}$ in $\mathscr T$ are assumed to control Type I error asymptotically given suitable partitions of the data. This assumption is formalized by Condition 6 in Section \ref{sec: main results}. 
We then choose $(\hat \Test, \hat  \Cluster) \in \mathscr T$ by solving 
\begin{equation}
	\label{eq: optimization problem}
	\begin{aligned}
		&(\hat \Test, \hat \Cluster) \in  \arg\min \ \widehat{\text{Err}}_{\text{Type-II}  }(\Test, \Cluster)  \\
		&\text{s.t. } {(\Test, \Cluster) \in \mathscr T \ ,  \ \widehat{\text{Err}}_{\text{Type-I}}(\Test, \Cluster)  \leq \alpha}.
	\end{aligned}
\end{equation}
The final testing outcome for $H_0$ is then denoted \begin{align} \hat T = \hat \Test(\mathscr D, \hat \Cluster( \mathsf X)).\end{align}

 In this paper, interest will be in $\mathscr T$ of the form \begin{align}\label{eq: scriptT}
	\mathscr T_{\bullet(\alpha), \G_{\max}} =\Big \{ (\Test_{\bullet(a)},  k\text{-}\texttt{medoids}) : a \in [0,\alpha], k \in \{ 2,...,{\G_{\mathrm{max}}} \}  \Big \}
\end{align}
where $\alpha$ is the nominal testing level for $H_0$, 
${\G_{\mathrm{max}}}$ is a researcher-chosen upper bound on the number of clusters, and $\Test_{\bullet(a)}$ is defined above.\footnote{One could also include pairs involving pre-specified partitions in $\mathscr{T}$. } 
With $\mathscr T = \mathscr T_{\bullet(\alpha), \G_{\max}}$, the ``parameter space''  in \eqref{eq: scriptT} depends on two independent parameters $a$ and $\G$.  It follows that \eqref{eq: optimization problem} is a two-dimensional optimization problem with a single constraint. 
The solution is then determined by two parameters $\hat \alpha$ and $\hat \G$.  Furthermore, the testing outcome can be expressed as
$\hat T = \Test_{\bullet( \hat \alpha )}(\mathscr D, \hat {\mathcal C} \hspace{.4mm} ), \ \text{with } \hat {\mathcal C} = \mathcal C^{(\hat \G)} = \hat\G \text{-}\texttt{medoids}(\mathsf X).$  When tests are based on $\bullet$ being  IM, CRS, or CCE and it is helpful to make the overall level $\alpha$ explicit, we write $\hat T = \hat T_{\bullet(\alpha)}$ for added emphasis.

$\hat \G$ provides a data-dependent answer to how many clusters to use. Optimization over $a$ allows the parameter entering the data-dependent decision rule, $\hat \alpha$, to be smaller than the nominal level of the test, $\alpha$. That is, the data-dependent decision rule may be more conservative than would be implied by conventional, fixed rules that asymptotically control size but may fail to do so in finite samples. Finally, the constraint on $a$ in the definition of $\mathscr{T}$ in \eqref{eq: scriptT} guarantees that inference based on $\hat T$ will maintain asymptotic size control under conditions that do not require the estimator $\widehat{\text{Err}}_{\text{Type-I}}(\Test, \Cluster)$ to agree with the true Type-I error rate in finite-samples or asymptotically. 

In practice, estimates $\widehat{\text{Err}}_{\text{Type-I}}(\Test, \Cluster)$ and $\widehat{\text{Err}}_{\text{Type-II}  }(\Test, \Cluster) $ of Type-I and Type-II error rates are needed. In all results reported in the following sections, we obtain these estimates from auxiliary estimation of the dependence structure in the data based on Gaussian Quasi Maximum Likelihood Estimation (QMLE) using a simple exponential covariance function. We then use the estimated dependence structure within a Gaussian model to obtain $\widehat{\text{Err}}_{\text{Type-I}}(\Test, \Cluster)$ and $\widehat{\text{Err}}_{\text{Type-II}  }(\Test, \Cluster)$. Further details are provided in Appendix \ref{Appendix A}. We note that in principle any baseline model could be used in implementing \eqref{eq: optimization problem} and one could consider uniform size control over classes of models, $\mathcal{M}$, by replacing the size constraint in \eqref{eq: optimization problem} with $\max_{m \in \mathcal{M}} \widehat{\text{Err}}_{\text{Type-I},m}(\Test, \Cluster)  \leq \alpha$. We reemphasize that our theoretical results do not require the consistency of the estimated dependence structure of $\Pr_0$ in order to control size. Misspecification in the model for dependence asymptotically leads only to potential loss of power.

\eqref{eq: optimization problem} differs from most methods in the literature for choosing data-dependent tuning parameters for use in conducting inference with dependent data. Much of the existing literature suggests choosing a single tuning parameter to optimize a weighted combination of size distortion and power; see, for instance, \cite{LLS:SizePower}, \cite{Sun2015}, and references therein. Instead, our proposal leverages the fact that most commonly used inferential procedures for dependent data depend on two parameters - nominal size and a smoothing parameter - and focuses on maximizing power within procedures that control size. Our proposal is closely related to \cite{MeullerWatson:2020} who consider an inference approach for spatially dependent data that makes use of a tuning parameter and a critical value which are chosen by minimizing confidence interval length subject to exactly controlling size. Relative to the existing literature, both \cite{MeullerWatson:2020} and our approach offer additional flexibility by explicitly considering two choice variables and make use of criteria, minimizing interval length or maximizing power subject to maintaining size control, that we believe will be appealing to many researchers.

For ease of reference, the above procedure is restated under Algorithm \ref{Algorithm1}.  For concreteness, Algorithm \ref{Algorithm1} is specialized to {IM}, CRS, CCE with $k$-\texttt{medoids}. 
To simplify notation, write $ \widehat {\mathrm{Err}}_{\mathrm{Type-I}}(\Test_{\mathrm{\bullet}(a)}, k\text{-}\texttt{medoids}) = \widehat {\mathrm{Err}}_{\mathrm{Type-I}}(\bullet(a), k)$  and $ \widehat {\mathrm{Err}}_{\mathrm{Type-II}}(\Test_{\mathrm{\bullet}(a)}, k\text{-}\texttt{medoids}) = \widehat {\mathrm{Err}}_{\mathrm{Type-II}}(\bullet(a), k)$. We use Algorithm \ref{Algorithm1} in the empirical example and simulation study in Sections \ref{empirical_application}, \ref{simulation}. The Appendix gives full implementation details. 

\vspace{3mm}{\singlespacing
\begin{algorithmsimple} \label{Algorithm1}   (\textit{Inference with Cluster Learning with} $k$-\texttt{medoids} 
\textit{and} IM, CRS \textit{or} {CCE}). 

\vspace{2mm} 
Testing $H_0$ at level $0 < \alpha < 1$. 
\vspace{2mm}

\textit{Data:} $\mathscr D, \mathsf X$.
\vspace{2mm}

\textit{Inputs:} $k_{\max}$;  $\mathbf \Theta_{\text{alt}}$; $\bullet = $ IM, CRS, or CCE; Estimates $\widehat{\text{Err}}_{\text{Type-I} }(\bullet(a), k), \widehat{\text{Err}}_{\text{Type-II}  }(\bullet(a), k)$
\vspace{2mm}

\textit{Procedure}:  Solve \eqref{eq: optimization problem} to obtain $(\hat \alpha, \hat k)$. 
\vspace{2mm}

\textit{Output:}  Set $\hat T_{\bullet(\alpha)} = {\Test}_{\bullet(\hat \alpha)}(\mathscr D, \hat k\text{-}\texttt{medoids}( {\mathsf X}))$.
\vspace{2mm}
\end{algorithmsimple}
}

Algorithm 1 can be inverted to generate confidence sets. For a family of hypotheses $\{ H_0: \theta_0 = \theta^{\circledast} : \theta^{\circledast} \in \Theta \}$, $\mathrm {C.I.} = \{\theta^{\circledast} \in \Theta : \hat T_{\bullet(\alpha)}  \text{= Fail to Reject for }  H_0:  \theta_0 = \theta^{\circledast}  \}$ is a $(1-\alpha)$ confidence set.

A formal analysis of Algorithm 1 is given in Section \ref{clt_section}.  
As an overview, we define an asymptotic frame relying on sequences of metric spaces with a suitable growth property. We then provide conditions under which the testing procedure in Algorithm 1 is consistent in Type I error -- i.e., conditions such that the superior limit of the Type I error of Algorithm 1 is bounded by $\alpha$ under the defined asymptotic frame. One of the important high-level conditions is that the tests under consideration would asymptotically control Type I error under a sequence of fixed partitions (i.e., without cluster learning). It is important to note that our theoretical results cover only tests which use CRS and IM.  We do not believe that consistency will hold for CCE at the level of generality available for CRS and IM. This point is discussed further in the Supplemental Material.

\section{Empirical Application: Insurgent Electoral Violence}
\label{empirical_application}

In this section, we illustrate the inferential method proposed in Section \ref{method} by reexamining estimates of the effect of insurgent attacks on voters from \cite{Condra2018}. We use the same data as \cite{Condra2018} which consist of district-level observations for 205 voting districts in Afghanistan in 2014. Each observation contains information on direct morning attacks, voter turnout, and district level control variables for two separate election rounds. 

\subsection{Model and Spatial Dependence.}

We focus on column (4) of Table 2 in \cite{Condra2018} which reports the estimated effect of nearby morning attacks on voter turnout obtained from instrumental variables (IV) estimation of the linear model
\begin{align}
	\label{eq: empirical structure}
	Y_{de}&=\alpha_0+\theta_0 Attacks_{de}+W_{de}'\gamma_0+U_{de} \\
	\label{eq: empirical first}
	Attacks_{de}&=\mu_0+\pi_0 Z_{de}+W_{de}'\xi_0+V_{de}
\end{align}
where \eqref{eq: empirical structure} is the structural equation with parameter of interest $\theta_0$ and \eqref{eq: empirical first} is the first stage representation. Here, $Y_{de}$ is the turnout in district $d$ and election round $e$, and $Attacks_{de}$ is the number of morning attacks in district $d$ and election round $e$. $U_{de}$ and $V_{de}$ are unobservables where $U_{de}$ and $Attacks_{de}$ are potentially correlated. Control variables $W_{de}$ include election round dummy variables, voting hour wind conditions, population, measures of precipitation and ambient temperature, and the average of predawn and morning wind conditions during the pre-election period. The excluded instrument, $Z_{de}$, is a measure of early morning wind conditions.

The analysis in \cite{Condra2018} reports inference for the parameter $\theta_0$ based on CCE with clusters defined by the cross-sectional unit of observation, the district.  For the sake of reference, from this point, the inference and corresponding method from \cite{Condra2018} will be named ``UNIT''. This approach will be asymptotically valid in the presence of dependence within district and heteroskedasticity as long as associations across districts are ignorable.  

To provide some evidence on dependence across districts, we provide \cite{moran} $\mathcal{I}$ tests; see also \cite{prucha:moran}. To construct the tests, we first estimate model \eqref{eq: empirical structure} to obtain estimates of $U_{de}$ which are used to construct the observation specific score for $\hat\theta$, the IV estimator of $\theta_0$ based on \eqref{eq: empirical structure}-\eqref{eq: empirical first}. These scores are the data input into the Moran $\mathcal{I}$ test. We consider two formulations for the weight matrix used in the Moran $\mathcal{I}$ test. In the first, we create the weight matrix by assigning a value of one for the weight corresponding to the two nearest neighbors for each observation on the basis of Euclidean distance in the latitude and longitude coordinates of the centroid of each district and set all other weights to zero (``Spatial weights''). We apply this formulation in each election round separately as well as in the full data set using both election rounds. In the second, we create the weight matrix by assigning a value of one to the weight corresponding to observations from the same district and set all other weights to zero (``Intertemporal weights''). The first setup should have power against alternatives where dependence is related to geographic distance and the second against alternatives where dependence is intertemporal. 

Panel A of Table \ref{table: moran}, displays the Moran $\mathcal{I}$ test statistics and the associated p-value for testing the null-hypothesis of zero correlation across observations based on the weight matrices described above. 
Each calculated test statistic is large, and each p-value for testing the null of zero correlation is small. The results suggest presence of both cross-sectional dependence associated with geographic distance and intertemporal dependence. As there is no immediate unique way to group districts on the basis of geographic proximity, constructing clusters in a data-dependent fashion may be helpful.



\subsection{Data-Dependent Cluster-Based Inference.}
\label{subsection_DataDepClustBasedInf}

We now apply the approach outlined in Section \ref{method} to construct data-dependent clusters of observations and produce inferential statements robust to the presence of spatial and temporal correlation and heteroskedasticity. A key input into this process is a dissimilarity measure $\mathrm d$ capturing the ``distance'' between each pair of distinct districts, $d, d'$. For this application, we use dissimilarity defined by Euclidean distance in geography 
\begin{align}\label{eq: emp distance} 
\mathrm d(d,d') = \| L_d - L_{d'} \|_2 = \sqrt{(\mathrm{lat}_d - \mathrm{lat}_{d'})^2+ (\mathrm{long}_{d} - \mathrm{long}_{d'})^2}
\end{align}
where $L_d = (\mathrm{lat}_d, \mathrm{long}_{d})$ consists of the latitude and longitude coordinates of the centroid of the district and $\| \cdot \|_2$ denotes Euclidean distance. 

Because of the nonuniform locations of district centers (i.e., districts are not on a rectangular grid), it is unclear ex ante how to partition them or how many clusters to use. As one potential resolution to this ambiguity, we apply $k$-\texttt{medoids} to the data using dissimilarity $\mathrm{d}$. Note that \eqref{eq: emp distance} assigns dissimilarity zero to observations from the same district which means that these observations will always be assigned to the same cluster. We thus allow general dependence across election rounds within district as well as capture dependence across election rounds and across districts that belong to the same cluster. Specifically, for each $\G = 2,...,\G_{\mathrm{max}}$ where $\G_{\mathrm{max}}$ is a user-specified upper bound on the number of clusters to be considered, we apply $k$-\texttt{medoids}, yielding a set of partitions $\mathscr C = \{ \mathcal C^{(2)}, ..., \mathcal C^{(\G_{\mathrm{max}})} \}$. As an illustration, Figure \ref{fig: afghanistan} displays a partition of the data resulting from applying $k$-\texttt{medoids} with $\G=6$ with districts rendered with common markers belonging to common clusters $\mathsf C \in \mathcal C^{(6)}$.\footnote{We display the result with $\G = 6$ as this is the solution obtained in \eqref{eq: optimization problem} in this example for CRS.} 


In this empirical illustration, we set $\G_{\mathrm{max}} = 8.$  In general, we recommend choosing a small value for $\G_{\mathrm{max}}$.\footnote{See Section \ref{choice:kmax} in the Supplemental Material for additional discussion.} 
This recommendation is based on two considerations. First, IM and CRS, our preferred procedures, rely on within-cluster estimates and the number of observations per cluster decreases quickly as the number of clusters considered increases, resulting in unstable performance when large numbers of clusters are considered. Second, for dependence structures with non-trivial dependence, size distortions typically become pronounced with larger numbers of groups while improvements to power tend to rapidly diminish once a modest number of groups are considered.

With the set of potential partitions, $\mathscr C$, in hand, the final inputs we need to be able to select a partition and testing procedure using \eqref{eq: optimization problem} are feasible guesses for Type I and average Type II error rates,  $\widehat{\text{Err}}_{\text{Type-I} }(\bullet(a), k)$ and $\widehat{\text{Err}}_{\text{Type-II}  }(\bullet(a), k)$. Let $U$ and $V$ respectively denote the $410 \times 1$ vectors produced by stacking $U_{de}$ and $V_{de}$. To obtain $\widehat{\text{Err}}_{\text{Type-I} }(\bullet(a), k)$ and $\widehat{\text{Err}}_{\text{Type-II}  }(\bullet(a), k)$, we use a parametric covariance model 
\begin{align} \begin{split}
\Ep
\begin{bmatrix}
	UU' & UV' \\ VU' & VV' 
\end{bmatrix} 
&= 
\begin{bmatrix}
	\Sigma_U(\tau^U) & \rho (\Sigma_U(\tau^{U})^{1/2}) (\Sigma_V(\tau^{V})^{1/2})'\\
	\cdot & \Sigma_V(\tau^V) 
\end{bmatrix} = \Sigma(\tau^U,\tau^V,\rho)
\end{split}\end{align} 
with parameters $(\tau^U,\tau^V,\rho)$ where $\Sigma_U(\tau)$ and $\Sigma_V(\tau)$ are covariance matrices with entries 
\begin{align}
	\label{eq: exp cov}
	\[\Sigma_{\varepsilon}(\tau)\]_{de,d'e'}=f_{\tau}((L_d,e),(L_d',e'))=\exp(\tau_1)\exp (-\tau_2^{-1} \Vert L_d-L_{d'}\Vert_2-\tau_3^{-1}|e-e'|)
\end{align}
for $\varepsilon \in \{U,V\}$. We obtain calibrated values for these parameters, $(\widehat\tau^U,\widehat\tau^V,\widehat\rho)$ by forming initial estimates $\hat{U}_{de}$ and $\hat{V}_{de}$ of $U_{de}$ and $V_{de}$ from IV estimation of \eqref{eq: empirical structure}-\eqref{eq: empirical first} and then using these residuals to estimate the covariance model via Gaussian QMLE. We can then approximate the exact Type I and average Type II error rates that would result from applying each inferential procedure with each potential partition within a Gaussian model with covariance matrix given by $\Sigma(\widehat\tau^U,\widehat\tau^V,\widehat\rho)$ by simulation. In this example, we use an equal weighted average of Type II error across alternatives $\theta^{\mathrm{alt}}  \in \{-10/\sqrt{n} ,-9/\sqrt{n},\dots,-1/\sqrt{n},1/\sqrt{n},2/\sqrt{n},\dots,10/\sqrt{n}\} \approx \{-.5,-.45,...,-.05,.05,...,.45,.5\}$. These alternatives correspond to an equally spaced grid of alternatives covering an interval that encompasses essentially all plausible values for the parameter of interest in this application. Note that favoring weighted average power is a common practice in econometrics; see e.g., \cite{moreira:moreira:optimalteststheory} for an overview. Given these estimates of the Type I and Type II error rates, we then solve \eqref{eq: optimization problem} to obtain the final partition and level for defining the cutoff in the decision rule used in the data. Further details are provided in the Appendix. For illustration of the problem, we present level sets of the function being minimized over and the constraint set for CCE, IM, and CRS in the Supplement.

Panel B of Table \ref{table: moran} reports inferential results for clustering at the district level as in \cite{Condra2018} (UNIT) and for the data-dependent CCE, IM, and CRS cluster-based inferential procedures. We report results for all procedures for completeness and illustration, but note that we only provide theory for IM and CRS. Our simulation results and the fact that existing theoretical results for CCE with fixed, known clusters rely on strong conditions lead us to believe that providing the formal theory for CCE is not worthwhile and to recommend that one uses IM or CRS in practice. For each method, we report the corresponding estimate of $\theta_{0}$, the standard error, $t$-statistic, 95\% confidence interval using  $\hat \G$ and $\hat\alpha$ from \eqref{eq: optimization problem}, 95\% confidence interval using $\alpha = .05$ and $\hat \G$, and $\hat\G$. 

	
In this illustration, all 95\% level confidence intervals include 0 with the exception of the interval produced by district level clustering with the asymptotic cutoff. All of the considered procedures which allow for spatial correlation would lead to failure to reject the null hypothesis that insurgent attacks have no effect on voter turnout at the 5\% level within the posited linear IV structure. There is nonetheless substantive variation in the interval estimates provided by the various procedures. 

Both procedures based on clustered standard errors (UNIT and CCE), while producing intervals that cover zero when spatial correlation is allowed for, are relatively narrow when compared to IM and CRS. For both UNIT and CCE, we can see a mild discrepancy between the intervals produced under the asymptotic approximation and those produced by allowing a data-dependent decision threshold based on the parametric model. Specifically, within the parametric model underlying the simulation, controlling size at the 5\% level would correspond to rejecting the null only when the p-value is smaller than 0.006 or 0.022 for UNIT and CCE, respectively (using the fixed number of groups approximation of \cite{BCH-inference}.) These differences suggest that the asymptotic approximation is potentially failing and suggest some caution should be exercised in relying on these results; see also \cite{Ferman:AssessInference}. We do note that the interval provided in column C.I. will exactly control size under the specified parametric model. However, we find that both CCE and individual level clustering are relatively sensitive to misspecification of the parametric model underlying construction of the estimated Type I and Type II error rates when compared to IM and CRS in the simulation results reported in Section \ref{simulation}.

Our preferred procedures are IM and CRS. These procedures are provably valid asymptotically under relatively weak conditions and are also relatively robust to misspecification of the parametric model underlying the size and power calculations in our simulation experiments. Looking to IM and CRS, we see similar performance, in the sense of relatively wide 95\% interval estimates, though the IM interval is substantially wider than CRS in this example. For both of these procedures, the asymptotic cutoffs coincide with those produced by solving \eqref{eq: optimization problem} using the simulated Type I and Type II error rates, which is reassuring though not necessary for validity. We note that there is a danger of a researcher trying both IM and CRS (and perhaps CCE) and reporting the most optimistic result. For this reason, our recommendation is to either default to, and compute only, IM or CRS 
or to report results from all procedures to provide complete information.

\section{Simulation}
\label{simulation}

This section examines the finite sample performance of inference based on data-dependent clusters in a series of simulation experiments roughly based upon the empirical illustration of Section \ref{empirical_application}. We present results for inference on a coefficient in a linear model with all exogenous variables and results for inference on the coefficient on an endogenous variable in a linear IV model. 

\subsection{Data Generating Processes and Inferential Procedures.}
\label{subsec: DGP}

We generate data $\mathscr D = \{ Y_{de}, X_{de}, W_{de} \}_{d,e=1,1}^{205,2}$ and $\mathscr D = \{ Y_{de}, X_{de}, W_{de},Z_{de} \}_{d,e=1,1}^{205,2}$ in the OLS and IV simulations.  
Each observation is associated with a spatial location $L_d = (\mathrm{lat}_d, \mathrm{long}_{d'})$, the latitude and longitude of the centroid of a district from the observed data in the empirical example (which also has 205 districts), and a temporal index $e \in \{1,2\}$.
In each case, $X_{de}$ is the variable of interest, $W_{de}$ is a $10 \times 1$ vector of control variables, and $Y_{de}$ is the outcome variable.  In the IV case, $Z_{de}$ is a scalar instrumental variable. 
We condition on a single realization of $\{X_{de},W_{de}\}_{d,e=1,1}^{205,2}$ in the OLS simulation and a single realization of $\{W_{de},Z_{de}\}_{d,e=1,1}^{205,2}$ in the IV simulation.\footnote{$X_{de},W_{del}\sim \mathrm N(0,1)$, $\mathrm{corr}(X_{de}$,$W_{d'e'l})$, $\mathrm{corr}(W_{del},W_{d'e'm})$ $=0.5\textbf{1}_{de=d'e'}$, $\mathrm{corr}(X_{de},X_{d'e'})$, $\mathrm{corr}(W_{del},W_{d'e'l})$ $=f_{\tau_0}( (L_d,e),(L_{d'},e'))$ for OLS.
For IV, $Z_{de},W_{del}$ have the same distribution as $X_{de},W_{del}$ do in OLS.
}   Conditionally, we re-draw observations of the remaining variables over simulation iterations according to the following design.  
 We set coefficients $\theta_0=0$, $\pi_0=2$, and $\gamma_0=\xi_0=0_{10}=(0,\dots,0)'$, $\tau_{0} = (0,3,1)'$.  Recall, $f_{\tau}$ is defined in \eqref{eq: exp cov}.
\begin{align}\label{eq: OLS sim}
\textbf{OLS Design.  } \quad \quad \quad \quad \quad \quad \quad& Y_{de} = \theta_0 X_{de} + W_{de}'\gamma_0 + U_{de} \quad \quad \quad \quad \quad&
\end{align}
\begin{itemize}
	\item[ A.] Homogeneous exponential covariance ({BASELINE}).  
	\begin{align}
		U_{de}\sim \mathrm N(0,1), \  \mathrm{corr}(U_{de},U_{d'e'})=f_{\tau_0}( (L_d,e),(L_{d'},e')).
	\end{align}
	\item[B.] Spatial auto-regression ({SAR}).
\begin{align}\begin{split}
	& U_{de}=0.15 \textstyle \underset{d'\ne d}\sum U_{d'e} \textbf{1}_{\Vert L_d-L_{d'}\Vert_2<0.3 }  +\varepsilon_{de}, \   \varepsilon_{de}\sim \mathrm N(0,1), \mathrm{corr}(\varepsilon_{d1},\varepsilon_{d2})=\exp(-1)\textbf{1}_{d= d'}.
\end{split}\end{align}
\end{itemize}
\begin{align}\label{eq: IV sim}
	\begin{split} 
	\textbf{IV Design.  } \quad \quad \quad  Y_{de}&=\theta_0 X_{de}+W_{de}'\gamma_0+U_{de}, \ \ 
	X_{de}=\pi_0 Z_{de}+W_{de}'\xi_0+V_{de} \quad \quad
	\end{split} 
\end{align}
\begin{itemize}
	\item[ A.] Homogeneous exponential covariance ({BASELINE}).  
	\begin{align} \begin{split}
		&U_{de},  V_{de}\sim \mathrm N(0,1), \
		\mathrm{corr}(U_{de},U_{d'e'})=\mathrm{corr}(V_{de},V_{d'e'})=f_{\tau_0}( (L_d,e),(L_{d'},e')),\\
		&\mathrm{corr}(U_{de},V_{de})=0.8 \textbf{1}_{(d',e')= (d,e)}.
	\end{split}\end{align}
	\item[B.] Spatial auto-regression ({SAR}).
	\begin{align}\begin{split}
		& U_{de}=0.15  \textstyle \underset{d' \neq d}  \sum U_{d'e} \textbf{1}_{ \Vert L_d-L_{d'}\Vert_2<0.3 } +\varepsilon_{de}\
				V_{de}=0.15 \textstyle  \underset{d' \neq d} \sum V_{d'e} \textbf{1}_{ \Vert L_d-L_{d'}\Vert_2<0.3 } +\eta_{de},\\
		& \varepsilon_{de}, \eta_{de}\sim \mathrm N(0,1),\
		\mathrm{corr}(\varepsilon_{d1},\varepsilon_{d2}), \mathrm{corr}(\eta_{d1},\eta_{d2})=\exp(-1)\textbf{1}_{ d = d'}, 
		\mathrm{corr}(\varepsilon_{de},\eta_{de})=0.8 \textbf{1}_{d'e' = de}.\\
	\end{split}\end{align}
\end{itemize}

\noindent
\textbf{Inferential Procedures:} Within each of the four simulation designs defined above, we generate 1000 simulation replications. We then report results for point estimation and for inference, focusing on size and power of hypothesis tests, about the parameter $\theta_0$ based on the following procedures:

\begin{itemize}
	\item[ 1.] SK. Inference based on the spatial HAC estimator from \cite{Sun2015} with bandwidth selection adapted from the proposal of \cite{Lazarus2018}.
	\item[ 2.] UNIT-U. Inference based on the cluster covariance estimator with clusters defined as the cross-sectional unit of observation with a critical value from a $t_{205-1}$-distribution.
	\item[ 3.] UNIT. Inference based on the cluster covariance estimator with clusters defined as the cross-sectional unit of observation and a critical value obtained by solving \eqref{eq: optimization problem} with cluster structure given by unit-level clustering. 
	\item[ 4.] CCE. Inference based on the cluster covariance estimator with clusters and rejection threshold obtained by solving \eqref{eq: optimization problem}.
	\item[ 5.] IM. Inference based on IM with clusters and rejection threshold obtained by solving \eqref{eq: optimization problem}.
	\item[ 6.] CRS. Inference based on CRS with clusters and rejection threshold obtained by solving \eqref{eq: optimization problem}.
\end{itemize}

For UNIT, CCE, IM, and CRS, we obtain preliminary estimates of unobserved components $U_{de}$ in the OLS setting and $(U_{de},V_{de})$ in the IV setting. We then apply Gaussian QMLE with the covariance structure specified in BASELINE using these preliminary estimates as data to obtain the structure for simulating Type I and Type II error rates for use in \eqref{eq: optimization problem}. Thus, results for UNIT, CCE, IM, and CRS in the BASELINE OLS and BASELINE IV settings illustrate performance when tuning parameters result from solving \eqref{eq: optimization problem} with a correctly specified model for the covariance structure with feasible estimates of covariance parameters. In contrast, results from the SAR OLS and SAR IV setting illustrate performance when tuning parameters from solving \eqref{eq: optimization problem} with a misspecified model. The Appendix provides implementation details for SK, UNIT, CCE, IM, and CRS.


\subsection{Simulation Results.}
\label{subsec: sim results}

We report size of 5\% level tests as well as results on point estimation in Table \ref{table_summary}. In the OLS simulation, we obtain point estimates by applying OLS to estimate the parameters of \eqref{eq: OLS sim}, and we obtain point estimates in the IV simulation by applying IV to estimate the parameters in \eqref{eq: IV sim}. Recall that both IM and CRS rely on first obtaining within cluster estimates and have natural point estimator defined by $\frac{1}{\hat \G} \sum_{\mathsf C \in \hat {\mathcal C}} \hat \theta_{\mathsf C}$ where $\hat\G$ and $\hat{\mathcal{C}}$ are the results from solving \eqref{eq: optimization problem} for each procedure and $\hat \theta_{\mathsf C}$ is a point estimator that uses only the observations in cluster $\mathsf C$. For point estimation, we report bias and root mean square error (RMSE) for the OLS simulation and median bias and median absolute deviation (MAD) for the IV simulation.

Table \ref{table_summary} demonstrates that both IM and CRS with data dependent groups and rejection threshold determined by \eqref{eq: optimization problem} control size across the four considered designs. CCE using \eqref{eq: optimization problem} also does a reasonable job controlling size across the designs, but with size distortions in the IV setting. The performance of IM, CRS, and CCE is similar both in the BASELINE setting, where $\hat\G$ and $\hat\alpha$ are obtained using a correctly specified covariance model, and in the SAR setting, where the problem solved to obtain tuning parameters makes use of a misspecified covariance model.

The robustness of IM, CRS, and, to a lesser extent, CCE is not exhibited by the remaining procedures. The poor behavior of UNIT-U, which ignores spatial dependence entirely, is unsurprising. More surprisingly, SK, which attempts to account for spatial dependence, also does poorly across all designs considered. Over-rejection of the spatial HAC estimator has been previously been documented in the literature. For example, \cite{Conley:Goncalves:Kim:Perron:boot} present simulation results where the null rejection rate for a 5\% level test reaches 0.600 and suggest a spatial dependence wild bootstrap approach to improve performance. We also suspect that results could substantially be improved by choosing tuning parameters and rejection rule for the spatial HAC procedure by solving a problem similar to \eqref{eq: optimization problem} adapted to that estimator. We did not pursue that direction as we wished to compare to a benchmark from the existing literature.

The UNIT procedure treats cross-sectional units as spatially uncorrelated in constructing the standard error estimator, but uses a parametric model that has spatial correlation to adjust the decision threshold for rejecting a hypothesis according to \eqref{eq: optimization problem}. We see that using a correctly specified parametric model for this adjustment restores size control, but size is not controlled under the misspecified parametric structure. This is consistent with our theoretical results which rely on the use of a small number of clusters to maintain size control while allowing for relatively general dependence structures and not requiring correct specification of the parametric model for tuning parameter choice. Robustness to misspecification is lost when large numbers of clusters are used.

The point estimation results mirror results already available in the literature; see, e.g., \cite{Ibragimov2010}, \cite{BCH-inference}, and \cite{CGH:JAR}. Both IM and CRS dominate the full sample OLS estimator in both bias and RMSE in our simulation designs. That is, there appears to be a gain, in terms of both point estimation properties and size control, in using the IM or CRS procedure. The results are muddier in the IV simulation where the full sample estimator exhibits lower median bias at the cost of larger MAD and poorer size control. 

We report power curves for 5\% level tests of the hypothesis $H_0: \theta_0 = \theta^{\mathrm{alt}} $ for alternatives $\theta^{\mathrm{alt}} $ produced by the different procedures across the simulations in Figures \ref{fig: ols power} and \ref{fig: iv power}. In these figures, the horizontal axis gives the hypothesized value, $\theta^{\mathrm{alt}} $, so size of the test is captured by the point $\theta^{\mathrm{alt}}  = 0$. Figure \ref{fig: ols power} presents the results from the OLS simulation. Here, the power curves are symmetric and the highest power among procedures that control size is obtained by IM and CRS, both of which perform similarly. Looking to the IV results in Figure \ref{fig: iv power}, we see that power curves are asymmetric and slightly shifted due to the finite sample behavior of the IV estimator. We see that there is no longer a clear picture about which of the procedures that controls size performs better in terms of power. Each of CCE, IM, and CRS exhibits higher power over different sets of alternative values. Exploring these tradeoffs further is potentially interesting but also outside our scope. 

We report properties of the data dependent number of clusters obtained from \eqref{eq: optimization problem}, $\hat\G$, across our simulation designs in Table \ref{tab: hat k}.  5\% level tests based on CRS have trivial power (power equal size) when based on fewer than six groups, so the number of groups selected for CRS is always chosen to be six or greater. 
We see that, with two exceptions, $\hat\G = 8$ is the most likely selection. 
The exceptions are for IM in the OLS and IV settings with SAR covariance structure where five groups are selected in 55.4\% and 56.8\% of the simulation replications, respectively. 
We also note that, while $\hat\G = 8$ is the most common solution, there is non-trivial weight on other numbers of groups for all procedures. Finally, we wish to reiterate that the grouping structure is not literally correct in any of our designs but is being used as a feasible way to downweight small covariances to allow for the construction of informative and robust inferential statements. Thus, the goal here is not to find the ``correct'' number of clusters.

Table \ref{tab: hat k} also reports the simulation distribution of the data dependent p-value thresholds, $\hat\alpha$, for rejecting a 5\% level test. Recall that using the asymptotic approximation underlying each procedure would correspond to rejecting with a threshold of .05. For CCE, the distribution of the threshold that would be used to provide size control at the 5\% level has substantial mass far below .05. This suggests the asymptotic approximation does not provide a reliable guide to the performance of inference based on CCE in our settings.  This is likely due to strong departures from homogeneity assumptions needed to establish the asymptotic behavior of CCE with small numbers of groups. For both IM and CRS, we see smaller shifts of $\hat\alpha$ away from .05. For both SAR cases, the distribution of $\hat\alpha$ for IM and CRS has a substantial mass point at .05. In the BASELINE cases, solving \eqref{eq: optimization problem} for both IM and CRS does lead to systematic use of p-value thresholds that are smaller than .05 to achieve 5\% level size control. These departures are significantly less extreme than when considering inference based on CCE. 
This behavior highlights how having two tuning parameters is useful for maintaining finite sample properties as size control may not be achievable in finite samples if one may only choose the number of groups for clustering and is unable to adjust the decision threshold.

The simulation results suggest that IM and CRS with data dependent clusters constructed as in Section \ref{method} control size in the simulation designs we consider while maintaining non-trivial power. 

\section{Formal Analysis of Inference with Learned Clusters}
\label{clt_section}

This section presents formal analysis of Algorithm \ref{Algorithm1}. Subsection \ref{subsec: finite ahlfors} lays groundwork by defining a notion of regularity for dissimilarity measures.  Subsections \ref{subsec: asympt ahlfors}--\ref{extension_to_ols} then define an asymptotic theory which is helpful for the analysis of Algorithm \ref{Algorithm1}.  Together, our regularity conditions are sufficient for establishing a mixing central limit theorem, as well as relevant properties of $k$-\texttt{medoids}.  

Sections \ref{subsec: finite ahlfors}--\ref{extension_to_ols}  are organized such that they build a ladder of increasingly higher-level regularity conditions.  The penultimate rung of this ladder, Theorem \ref{iv_theorem}, concludes with concrete properties for $\hat \theta_{\mathsf C}$ for a sequence of linear models in which Algorithm \ref{Algorithm1} with $\mathsf C \in \mathscr C$ obtained by $k$-\texttt{medoids} has been applied. The conclusions of Theorem \ref{iv_theorem} are then used as a final rung of high-level conditions in Sections \ref{subsec: im theorem} and \ref{subsec: crs theorem} to allow application of Theorems  \ref{im_theorem_main} and \ref{crs_theorem_main} which analyze size control for IM and CRS procedures with inputs more general than those arising from Algorithm \ref{Algorithm1}. Versions of Theorems \ref{im_theorem_main} and \ref{crs_theorem_main} were previously established in \IM \ and \CRS; and thus the primary theoretical contribution in this paper is the development of the theory in Sections \ref{subsec: finite ahlfors}--\ref{extension_to_ols}. Theorems  \ref{im_theorem_main} and \ref{crs_theorem_main} do contain refinements relative to previous results which allow simultaneous analysis using multiple partitions. 
All proofs are given in the Supplemental Appendix.

\subsection{A notion of regularity for dissimilarity measures.}
\label{subsec: finite ahlfors}

This section defines a measure of regularity for a spatial indexing set $\mathsf X$ with dissimilarity $\mathrm d$.   The most basic condition on $\mathsf X$ is that the triangle inequality holds, making $\mathsf X$ a genuine metric space.  The next definition is used to restrict the growth rate of balls centered around any point of $\mathsf X$.  Let $|\mathsf X|$ be the cardinality of $\mathsf X$.  Define $\mathsf X$ to be $(C,\delta)$-\textit{finite-Ahlfors regular}, written $\mathsf X \in \textbf{Ahlf}_{C,\delta}$,  if it is a metric space and 
$\min( |\mathsf X| , C^{-1}r^{\delta} ) \leq |\textsf{B}_{\mathsf X,r}(i)| \leq \max (Cr^{\delta} , 1)$ for any $r>0$, for any $i$ in any $\mathsf X$, where $\textsf{B}_{\mathsf X,r}(i)$ is the $r$-ball centered at $i$ in the space $\mathsf X$. \textcolor{black}{Heuristically, the above bounds rule out irregular distributions of observations where a large number of observations are densely packed in a relatively small region.}

Using $\mathsf X \in \textbf{Ahlf}_{C,\delta}$ as a notion of regularity has advantages relative to requiring $\mathsf X$ is a subset of a Euclidean space. First, $\mathsf X \in \textbf{Ahlf}_{C,\delta}$ references only intrinsic properties of $\mathsf X$.  $\mathsf X \in \textbf{Ahlf}_{C,\delta}$ is sufficient both for realizing mixing central limit theorems and for analyzing clustering.  

The condition that $\mathsf X \in \textbf{Ahlf}_{C,\delta}$ is indeed not strong enough to guarantee that $\mathsf X$ is isometric with a subset of a Euclidean space.   However, the condition does ensure that $\mathsf X$ can be effectively regularized.  The new space $(\mathsf X,\mathrm d^{1-\eps})$, in which the exponent $1-\eps$ is applied element-wise to $\mathrm d$, is a metric space for all $\eps \in (0,1)$ and is called the $\eps$-snowflake of $\mathsf X$.  This new space can be embedded into $\mathbb R^\nu$ with bounded distortion. This fact relies on Assoud's Embedding Theorem (see \cite{assoud:diss}) and a calculation controlling the doubling dimension of $\mathsf X$ in terms of $C, \delta$.  

\subsection{Asymptotics A: Frames defined by an increasing sequence of dissimilarity measures.}
\label{subsec: asympt ahlfors} 

The analysis of Algorithm \ref{Algorithm1} relies on asymptotic theory.  
We build an asymptotic frame by delineating the classes of sequences of spatial indexing sets considered.  
A sequence of finite dissimilarity measures will be denoted with ${\mathsf X}_{\shortrightarrow} = ( \mathsf X_n )_{n \in \mathbb N}$.  
Note that in this section, unlike in the earlier sections, $n$ appears as an explicit subscript in $\mathsf X_n$ precisely because $\mathsf X_n$ are elements of a newly introduced sequence.   Previously, $n$ was not needed to define quantities like $\hat T$.  
From this point, any $\mathsf X_n$ appearing with a subscript always implicitly belongs to a sequence $\mathsf X_{\shortrightarrow}$.

\begin{assumption}
	\label{ahlfors} 
	(Ahlfors Regularity) The sequence of dissimilarity measures ${\mathsf X}_{\shortrightarrow}$ 
	 satisfies $|\mathsf X_n| = n$ and is a sequence of the form $\mathbb N \rightarrow \mathbf{Ahlf}_{C,\delta}$ for some $0<C < \infty, 1\leq \delta < \infty $.     
\end{assumption}

Condition \ref{ahlfors} defines a spatial asymptotic frame.\footnote{\textcolor{black}{Note, $|\mathsf X_n| = n$ implies each observation has a unique index. Our formal results only consider $\mathsf X$ such that $\mathrm d(i,j) = 0\Rightarrow i=j$. To reconcile this convention with the data structure of the empirical example and simulation study above, it is necessary either to consider observations $(Y_{d1}, X_{d1}, W_{d1}, Z_{d1}),(Y_{d2}, X_{d2}, W_{d2}, Z_{d2})$ as a common observation arising from a single index $d$, or to consider an expanded definition of $\textbf{Ahlf}_{C,\delta}$.  Both possibilities can be carried out without essential changes to the arguments in this paper.  For brevity, these details are omitted.}}}  
This notion is simple and gives almost enough structure to allow an analysis of $k$-\texttt{medoids} techniques analytically as well as derive dependent central limit theorems and laws of large numbers.  
 Examples of metric space sequencing within $\mathbf{Ahlf}_{C,\delta}$ for fixed $C, \delta$ include $\mathbb Z^m \ \cap \ \mathsf {Sq}_n$ where $\mathsf {Sq}_n$, $m$-dimensional cubes of side length $2n$, $\mathbb Z^m \ \cap \ \mathsf A_n$ where $\mathsf A_n$, $m$-dimensional annuli of outer radius length $2n$ and inner radius length $n$.  

\textcolor{black}{We provide our results in an environment where observation indexes are non-random. Alternatively, we could carry out the analysis conditional on the observed locations under the assumption that they are exogenous. This treatment would rule out models with endogenous network formation or location choice which may be more salient when dependence is indexed by economic location.}

Further discussion on the role of Ahlfors regularity is provided in supplemental material.  There, we concretely calculate bounds on $| \mathsf B_r(i) |$ for $1 \leq r \leq diam(X)$ and $i \in \mathsf X$ for three examples. 
As a summary, (1) in the \cite{Condra2018} example in which locations are geographic, $0.002 r^2 \leq |B_r(i)| \leq 2.000 r^2$ holds, (2) in a square subset of the integer lattice with side length $\lceil \sqrt{205} \rceil=15$, $ 0.500 r^2 \leq |B_r(i)| \leq 4.000 r^2$ holds, and (3) in the distances between sectors in \cite{Conley:Dupor:Sectoral}, defined as the Euclidean distances between its share vectors computed from the input-output tables,  $0.125 r^2 \leq |B_r(i)| \leq 3.000 r^2$ holds.     

%


\subsection{Asymptotics B: Mixing conditions and central limit theory.}
\label{subsec: mixing}

This section develops a central limit theory over $\mathsf X_n \in \textbf{Ahlf}_{C,\delta}$.
Let $\mathscr D_{\shortrightarrow} = (\mathscr D_n)_{\n \in \mathbb N}$ and  $\mathscr D_n = \{\zeta_{i}\}_{i\in \mathsf X_n}$ be an array of real random vectors on a probability space $(\Omega, \mathscr F, \Pr)$.  Let $\mathscr A$ and $\mathscr B$ be two (sub-)$\sigma$-algebras of $\mathscr F$.  Define the following \textit{mixing coefficients}.   
\begin{align}\alpha^{\mathrm{mix}}(\mathscr A, \mathscr B) = \sup \{|\Pr(A \cap B) - \Pr(A) \Pr(B) | : A \in \mathscr A, B \in \mathscr B \}.\end{align}  For $\mathsf U, \mathsf V \subseteq \mathsf X_n $ and integers $g,l,r,$ define additionally 
\begin{align} \begin{split}
& \alpha_{n}^{\mathrm{mix}}(\mathsf U,\mathsf V) = \alpha^{\mathrm{mix}}(\sigma(\zeta_i: i \in \mathsf U), \sigma(\zeta_i: i \in \mathsf V)),\\
&\alpha_{g,l,n}^{\mathrm{mix}}(r) = \sup_{\mathsf U, \mathsf V \subseteq \mathsf X_n} \left \{\alpha^{\mathrm{mix}}_{n}(\mathsf U,\mathsf V), |\mathsf U| \leq g, |\mathsf V| \leq l, \mathrm d(\mathsf U,\mathsf V) \geq r\right \} ,\  \text{and }
\bar \alpha_{g,l}^{\mathrm{mix}}(r) = \sup_{n\in \mathbb N} \alpha^{\mathrm{mix}}_{g,l,n}(r).\end{split}\end{align}

The above definitions are the same as used in \cite{Jenish2009}.   The next condition relies also on the following fact.  Let $(\mathsf X_n,\mathrm d_n )\in \mathbf{Ahlf}_{C,\delta}$ for every $n$.  Then each $(\mathsf X_n,\mathrm d_n^{3/4})$ has an $L$-bi-Lipschitz map into $ \mathbb R^{\nu}$ with $\nu$ and $L$ depending only on $C, \delta$.  

\begin{assumption}[Mixing] 
	\label{mixing_assumption}
The array $\mathscr D_{\shortrightarrow}$ is an array of random vectors $\zeta_i$ on a probability space $(\Omega, \mathscr F, \Pr)$ taking values in a finite dimensional real vector space, with spatial indices given by $\mathsf X_{\shortrightarrow}$ satisfying Condition {1} with Ahlfors constants $C< \infty, \delta\geq 1$.  Let $\nu$ correspond to $C, \delta$ as defined in the preceding paragraph.  For any $\mathsf C \subseteq \mathsf X$, where $\mathsf X = \mathsf X_n$ for some $n$, let $\sigma^2( \mathsf C )=\mathrm{var}(\sum_{i \in \mathsf C }\zeta_{{i}})$.  There are constants $\underline c, \overline c$, $\{\{c_{i}\}_{i \in \mathsf X_n}\}_{n\in \mathbb N}$, $0 < \underline c < c_i< \overline c < \infty$, and $\mu>0$ such that:

\noindent \textit{(i) }$ \Ep[\zeta_{i}]=0$;
\textit{(ii)} $\lim_{k \rightarrow \infty} \sup_{n\in \mathbb N} \sup_{i \in \mathsf X_n} \Ep[ \| \zeta_{i} / c_{i}\|_2^{2+\mu} \textbf{\em 1}_{\|\zeta_{i}/c_{i}\|_2 >k } ] = 0$; 
\textit{(iii) } $\sum_{m=1}^\infty \bar \alpha^{\mathrm{mix}}_{1,1}(m)^{\frac{\mu}{\mu+2}} m^{ \nu}< \infty$; 
\textit{(iv)} $\sum_{m=1}^\infty m^{\nu-1} \bar \alpha^{\mathrm{mix}}_{g,l}(m) < \infty$, $g+l\leq 4$; 
\textit{(v)} $\bar \alpha^{\mathrm{mix}}_{1,\infty}(m)  = O(m^{-\nu-\frac{4}{3}\mu})$; 
\textit{(vi) }$\inf_{n \in \mathbb N} \inf_{\mathsf C \subseteq \mathsf X_n}\frac{\max_{i \in \mathsf C}\lambda_{\min} c_{i} ^{-2} \lambda_{\min}\sigma(\mathsf C)^2}{ |\mathsf C|} >0.$
\end{assumption}

The conditions are similar to those given for the mixing central limit theorem in \cite{Jenish2009} which requires that each $\mathsf X_n$ be a possibly uneven lattice in a finite dimensional Euclidean space with a minimum separation between all points.   
In the case that $\mathsf X_n$ embeds isometrically into $\mathbb R^\nu$ without the need to apply the snowflake regularization construction, Condition \ref{mixing_assumption} and the conditions for Corollary 1 of \cite{Jenish2009} are mildly different. First, we take $\bar \alpha^{\mathrm{mix}}_{1,\infty}(m)  = O(m^{-\nu-\frac{4}{3}\mu})$ rather than the slightly weaker $\bar \alpha^{\mathrm{mix}}_{1,\infty}(m)  = O(m^{-\nu-\mu})$.  Second, Condition \ref{mixing_assumption}(\textit{vi}) entails an infimum over $n$ and over all subsets $\mathsf C \subseteq \mathsf X_n$ while \cite{Jenish2009} only requires the condition hold for an infimum over $n$ with the particular choice $\mathsf C  = \mathsf X_n$.  Third,  \ref{mixing_assumption}(\textit{iii}) is a mixing condition similar to that assumed in \cite{bolthausen1982} (with exponent $\nu$ replacing $\nu - 1$), which is useful for summing certain covariance terms, and which is stronger than the corresponding condition in \cite{Jenish2009} that imposes only $\sum_{m=1}^\infty \bar \alpha^{\mathrm{mix}}_{1,1}(m)m^{ \nu \times \frac{\mu+2}{\mu } - 1} < \infty$.  In terms of notation, here, $\zeta_{i}$ is used rather than the more explicit $\zeta_{i,n}$, because formally the indexes $i$ belong to distinct sets, $\mathsf X_n$, for each $n$.  As a result, Condition \ref{mixing_assumption} is a condition on arrays of random variables.  The existence of $\underline c$ and $\overline c$ are not needed for the following theorem (Theorem 1) to hold.  They are convenient for exposition, as they imply for certain sequences $\mathsf C_n \subseteq \mathsf X_n$ of interest, that  $|\mathsf C_n |^{-1} \sigma(\mathsf C_n)$ stay suitably bounded from above and away from zero.

In the next theorem, let $\Phi$ be the Gaussian cumulative distribution function.  Also, as in the statement of Condition \ref{mixing_assumption}, for $\mathsf C \subseteq \mathsf X$ for any $\mathsf X = \mathsf X_n$ for some $n$, let $\sigma(\mathsf C) = \mathrm{var}(\sum_{i\in \mathsf C}\zeta_i)$.
\begin{theorem} 
	\label{clt_mean}
	Suppose that  $\mathsf X_{\shortrightarrow}$ satisfies Condition \ref{ahlfors} and  $\mathscr D_\shortrightarrow$ satisfies Condition \ref{mixing_assumption} with scalar $\zeta_i$.  
Then for every $z \in \mathbb R$, 
$\lim_{n \rightarrow \infty} \Pr( \sigma(\mathsf X_n)^{-1}\sum_{i \in \mathsf X_n}\zeta_{i} \leq z) = \Phi(z)$.
\end{theorem}

Theorem \ref{clt_mean} proves convergence in distribution to a Gaussian measure for spatially indexed arrays $\mathscr D_n = \{\zeta_i\}_{i \in \mathsf X_n}$ satisfying Conditions \ref{ahlfors} and \ref{mixing_assumption}.   
 The proof involves regularizing $\mathsf X_n$ using the snowflake construction with $\eps = 1/4$ and applying the \cite{Jenish2009} central limit theorem.  

\subsection{Asymptotics C: Balance and small common boundary conditions.}
\label{balanced_small_boundary}

A key observation in \BCH \  is that if a fixed sequence of partitions of $\mathsf X_n$ have balanced cluster sizes and small boundaries, then under mixing and additional regularity conditions, the CCE procedure as described in \BCH \ achieves asymptotically correct size. This phenomenon also holds in our asymptotic settings.  This section formalizes notions of balanced cluster sizes and small boundaries appropriate for sequences $\mathsf X_{\shortrightarrow}$ satisfying Condition \ref{ahlfors}.  This section then gives a proposition showing that under the appropriate regularity conditions, the endpoint of the $k$-\texttt{medoids} algorithm satisfies these two desired properties.  Let $\mathcal C_{\shortrightarrow} = ( \mathcal C_n )_{n \in \mathbb N}$ be a sequence such that $\mathcal C_n$ is a partition of $\mathsf X_n$ for each $n$.

\begin{assumption} [Balance and Small Boundaries]
\label{balanced_small_bound_assumption}
The sequence of partitions $\mathcal C_{\shortrightarrow}$ satisfies $| \mathcal C_n| \leq \G_{\max}$ for some $\G_{\max}$ fixed with $n$, and is asymptotically balanced with small boundaries in that 
	\begin{enumerate}[label=(\roman*)]
		\item The cluster sizes satisfy
		$\underset{\n\rightarrow \infty}{\lim\inf} \hspace{1mm} \frac{ \hspace{1mm} \min_{\mathsf C\in \mathcal{C}_n}{|\mathsf C|} \hspace{1mm} } {{\max_{\mathsf D\in \mathcal{C}_n}{|\mathsf D|}} }>0,$
\item There is a  sequence $\bar r_{\shortrightarrow} = (\bar  r_n)_{n \in \mathbb N} \nearrow \infty$ with $\lim_{n \rightarrow \infty}\frac{\max_{\mathsf{C}\in \mathcal{C}_n} | \{ i\in \mathsf C :\mathrm  d(i, \mathsf X_n \setminus \mathsf C) \leq \bar r_n\} | }{ \min_{\mathsf D \in \mathcal C_n}|\mathsf D|} = 0. $
\end{enumerate}

\end{assumption}

The definition of small boundary provided in Condition \ref{balanced_small_bound_assumption}(\textit{ii}) differs slightly from the definition given in \BCH, who leverage the fact their spatial domain is a subset of the integer lattice to define neighbor orders for pairs of locations.  Their definition of small boundaries entails a bound on the number of first order neighbors from $\mathsf C$ to $\mathsf X_n \setminus \mathsf C$.  
 In this context, there is no available definition of first order neighbor since $\mathsf X_n$ can be irregular (even non-Euclidean).  As a result, this paper works instead with an asymptotic notion of boundary which entails a sequence $\bar r_{\shortrightarrow}$ which allows boundaries to widen as $n\rightarrow \infty.$  Theorem \ref{asymptotic_independence} gives the relevant implication of Condition \ref{balanced_small_bound_assumption}.

\begin{theorem}
	\label{asymptotic_independence}
Suppose that $\mathsf X_{\shortrightarrow}$ satisfies Condition \ref{ahlfors}, $\mathscr D_{\shortrightarrow}$ satisfies Condition \ref{mixing_assumption}, and $\mathcal C_{\shortrightarrow}$ satisfies Condition \ref{balanced_small_bound_assumption}.  
Then $\lim_{n\rightarrow \infty} \sup_{\mathsf C \neq \mathsf D \in \mathcal C_n} \mathrm{cov} \( \sigma({\mathsf C})^{-1}  \sum_{i \in \mathsf C} \zeta_{i}, \ \sigma({ \mathsf D})^{-1} \sum_{i \in \mathsf D} \zeta_{i} \) = 0.$ 

\end{theorem}

\noindent The argument establishing the above fact is related to but not identical to arguments in \BCH, which were previously also given in  \cite{Jenish2009} and \cite{bolthausen1982}. Instead of counting points in ``shells'' around the boundaries of clusters, the proof of Theorem \ref{asymptotic_independence} instead relies on the doubling structure implied by the fact that $\mathsf X_n$ are Ahlfors regular.  Both arguments leverage a bound on covariances $\text{cov}(\zeta_{i},\zeta_{j})$ for sufficiently distant locations as implied by the mixing conditions stated in Condition \ref{mixing_assumption}.

Next, we verify that clusters produced by  $k$-\texttt{medoids} satisfy Condition \ref{balanced_small_bound_assumption} under an additional convexity assumption. 
$k$-\texttt{medoids} is a popular clustering technique that is related to $k$-\texttt{means}, and both produce similar clustering results in many settings. However, $k$-\texttt{means} centroids are not necessarily defined in the space $\mathsf X_n$ when dealing with a dissimilarity that does not necessarily arise from a Euclidean space. $k$-\texttt{medoids} differs by requiring that clusters are defined around elements of $\mathsf X_n$, called medoids, which must themselves also be elements of $\mathsf X_n$. 
Let
\begin{align}\mathcal C_\shortrightarrow = ( \mathcal C_n )_{n\in\mathbb N} \text{ defined by } \mathcal C_n =k\text{-}\mathtt{medoids}(\mathsf X_n).\end{align}

\begin{theorem} 
	\label{km_prop}
	Suppose that $\mathsf X_\shortrightarrow$ satisfies Condition \ref{ahlfors}.  Assume the following additional convexity condition.  There is a constant $K$ independent of $n$ such that for each $n$, (1) $\mathsf X_n$ is $K$-coarsely isometric\footnote{ $f:(\mathsf Y,\mathrm d_\mathsf Y) \rightarrow (\mathsf Z, \mathrm  d_\mathsf Z)$ is a  $K$-coarse isometry if $\mathrm d_\mathsf Z(f(i),f(j)) -K \leq \mathrm d_{\mathsf Y}(i,j) \leq \mathrm d_\mathsf Z(f(i),f(j)) +K$.} to a subset of a Euclidean space with dimension $\nu$ independent of $n$, and (2) for any two point $i,j \in \mathsf X_n$ and any $a \in [0,1]$ there is an interpolant $g \in \mathsf X_n $ such that $ | \mathrm d_n(i,g) - a \mathrm d_n(j,g) | \leq K$ and $ | \mathrm d_n(i,g) - (1-a) \mathrm d_n(j,g) | \leq  K$.   Let $\mathcal C_n =k\text{-}\mathtt{medoids}(\mathsf X_n)$ and let $\mathcal C_\shortrightarrow = ( \mathcal C_n )_{n\in\mathbb N}$.  Then $\mathcal C_\shortrightarrow$ satisfies Condition \ref{balanced_small_bound_assumption}.
\end{theorem}    

Theorem \ref{km_prop} does not require the sequence of partitions $\mathcal C_n$ of $\mathsf X_n$ resulting from applying $k$-\texttt{medoids} to converge. 
Nor does it require there be any notion of a true partition associated to the $\mathsf X_n$.  

The convexity assumption excludes cases such as the following would-be counterexample: 
$\mathsf {Rothko}_n  = \{ i \in \mathbb Z^2 : |i_1| \leq n, |i_2| \leq 2n \} \setminus  \{ i \in \mathbb Z^2 : | i_1 | \leq n-2, 1 \leq |i_2+\frac{1}{2}|  \leq 2n-2 \}. $
When plotted, $\mathsf {Rothko}_n$ resembles the background orange region in 1961 Rothko titled \textit{Orange and Yellow.}  A feasible stopping point for $2$-\texttt{medoids} for the above example is with medoids at $(0,0)$ and $(0,-1)$ and boundary (with $\bar r = 1$) at $\{ 0,1  \} \times \{-n, -n+1,...,n\}$.  These endpoints are not unique.  $(-n,2n), (n, -2n)$ are also feasible and do lead to small boundaries.

\subsection{Asymptotics D: Almost sure representation.} \label{subsec: as rep} 
The next preliminary for the analysis of Algorithm \ref{Algorithm1} is an almost sure representation theorem.  Suppose that $\mathscr C_\shortrightarrow$ is a sequence of collections of partitions of $\mathsf X_n$.  That is $\mathscr C_{\rightarrow} = (\mathscr C_n )_{n \in \mathbb N}$ and each $\mathscr C_n$ is a set containing elements $\mathcal C$ which are partitions of $\mathsf X_n$.  The following condition concerns corresponding collections $S_{\mathcal C}$ of random vectors taking values in $\mathbb R^{\mathcal C}$ where $\mathcal C \in \mathscr C_n$ for some $n$.  Let $\{S_\mathcal C\}$ be a collection of random vectors in $\mathbb R^\mathcal C$ where $\mathcal C$ range in $\cup_{n \in \mathbb N} \mathscr C_n$.

\begin{assumption}[Almost Sure Asymptotic Gaussian Representation]\label{con: a s gaussian}
\textcolor{black}{There is a probability space $(\tilde \Omega, \tilde{ \mathscr F}, \tilde \Pr)$ on which random variables $S_{\mathcal C} \in \mathbb R^{\mathcal C}$ are simultaneously defined for $\mathcal C$ ranging in $ \cup_{n\in\mathbb N}\mathscr C_n$.}
Furthermore, $\tilde{S}_{\mathcal{C}}$ and ${S}_{\mathcal{C}}$ have the same distribution, $\tilde{S}_{\mathcal{C}}^*$ is zero-mean Gaussian with independent components, $\liminf_{n \rightarrow \infty} \min_{\mathsf C \in \mathcal C \in \mathscr C_n} \mathrm{var}(\tilde S_{\mathcal C, \mathsf C}^*)> 0$, $\limsup_{n \rightarrow \infty} \max_{\mathsf C \in \mathcal C \in \mathscr C_n} \mathrm{var}(\tilde S_{\mathcal C_n, \mathsf C}^*)< \infty$ and
	$\lim_{n \rightarrow \infty} \sup_{\mathcal C \in \mathscr C_n} \| \tilde{S}_{\mathcal C} - \tilde{S}_{\mathcal C}^*\|_2=0$ $\tilde \Pr$-almost surely. 
\end{assumption}

The representors relative to $S_{\mathcal C}$ will generally be referred to by $\tilde S_{\mathcal C}^*$ and $\tilde S_{\mathcal C}$ just as in the definition.
Convergence in distribution of a sequence of random variables implies almost sure convergence for an auxiliary sequence with the same distribution as the original sequence: Theorem 2.19 in \cite{vdV} states:  ``\textit{Suppose that the sequence of random vectors $X_n$ converges in distribution to a random vector $X_0$.  Then there exists a probability space $(\tilde \Omega, \tilde {\mathcal U}, \tilde P)$ and random vectors $\tilde X_n$ defined on it such that $\tilde X_n$ is equal in distribution to $X_n$ for every $n\geq 0$ and $\tilde X_n \rightarrow \tilde X_0$ almost surely.}''


One collection of $S_\mathcal C$ of interest are constructed from $\mathscr D_\shortrightarrow$ and partitions $\mathcal C_\shortrightarrow$ of $\mathsf X_\shortrightarrow$ by 
\begin{align} S_{\mathcal C, \mathsf C} = (n/\G)^{-1/2} \sum_{i \in \mathsf C} \zeta_i, \text{ for } \mathcal C \text{ partitioning }\mathsf X_n \text{ and } |\mathcal C|=k.\end{align}


\begin{theorem}
	\label{ASREPversion}
Suppose that $\mathsf X_\rightarrow $ satisfies Condition \ref{ahlfors}, $\mathscr D_\rightarrow$ satisfies Condition \ref{mixing_assumption}, and $\mathcal C_\rightarrow$ satisfies Condition \ref{balanced_small_bound_assumption}.  For $\mathcal C$ with $|\mathcal C |=k$, let $S_{\mathcal C, \mathsf C} =(n/\G)^{-1/2} \sum_{i \in \mathsf C} \zeta_i$.  Let $\mathscr C_n $ be the singleton $\mathscr C_n = \{ \mathcal C_n\}$.  
Then $\{ S_{\mathcal C}\}_{\mathcal C \in \cup_{n\in \mathbb N} \mathscr C_n}$ satisfies Condition \ref{con: a s gaussian}. 
\end{theorem}

\subsection{Asymptotics E: Central limit theory in the cases of OLS and IV.}
\label{extension_to_ols}

This section gives a version of Theorem \ref{ASREPversion} in the context of ordinary least squares (OLS) and instrumental variables (IV) estimators.  These are the precise central limit theorems that are relevant to subsequent analysis of cluster-based inferential procedures.
Consider data $\mathscr D_n = \{ \zeta_i \}_{i \in \mathsf X_n}$ defined on $(\Omega, \mathscr F, \Pr)$ with $\zeta_i = (Y_i, X_i, W_i, Z_i)$ in which for each $ i\in \mathsf X_n$ (and for each $n$) a linear model: 
\begin{align}Y_{i}=\theta_{0}X_{i} + W_i' \gamma_0+U_{i},\end{align}
where $Y_{i}$ is a scalar outcome variable, $X_i$ is a scalar variable of interest, $U_{i}$ is a scalar idiosyncratic disturbance term, $W_i$ is a $p$-dimensional vector of control variables, and $Z_{i}$ is a $(p+1)$-dimensional vector of instruments.  The following regularity condition is helpful.

\begin{assumption}\label{con: lin reg} $\mathscr D_\shortrightarrow$ consists of $\mathscr D_n = \{(Y_{i}, X_{i}, W_{i}, Z_{i})\}_{i\in \mathsf{X}_n}$ satisfying the linear model $Y_{i}=\theta_{0}X_{i} + W_i' \gamma_0+U_{i}$ for some $\theta_0$ and $\gamma_0$.  Let $\bar X_i = (X_i, W_i)$.  $\Ep[Z_iU_i] = 0$, $\Ep[Z_i\bar X_i']$ exists and equals some invertible $\mathrm M$ which does not depend on $i$.  $(\{
Z_iU_i, Z_i \bar X_i' - \Ep[Z_i\bar X_i'] \}_{i \in \mathsf X_n})_{n\in\mathbb N}$ satisfies Condition \ref{mixing_assumption}.  
\end{assumption}

Consider the IV estimator for $\hat \theta_{\mathsf C}$ for $\theta_{0}$ using only observations $\mathsf C\subseteq \mathsf X_n$. 
The OLS estimate is a special case where $Z_i =X_i$. 
Let $\mathcal C$ partition $\mathsf X_n$ for some $n$ into $\G$ clusters.  Define $S_\mathcal C \in \mathbb R^\mathcal C$ by  
\begin{align}\label{eq: SC}
	 S_{\mathcal C, \mathsf C} =  (n/\G)^{1/2}\hspace{.5mm} (\widehat{\theta}_{\mathsf C}-\theta_0) \ \text{for} \ \mathsf C \in \mathcal C. 
\end{align}

Let $\mathscr C_n$ again be the singleton $\mathscr C_n = \{ \mathcal C_n \}$.  The next theorem is the particular asymptotic normality result about $\hat \theta_{\mathsf C}$, $\mathsf C \in \mathcal C_n$ which is useful for input into the analysis of Algorithm \ref{Algorithm1} in the context of testing hypotheses about coefficients in a linear model.

\begin{theorem}
	\label{iv_theorem}  Suppose that $\mathsf X_{\shortrightarrow}$ satisfies Condition \ref{ahlfors}.  Suppose linear regression data $\mathscr D_\shortrightarrow$ satisfies Condition \ref{con: lin reg}. Suppose that $\mathcal C_\shortrightarrow$  satisfies Condition \ref{balanced_small_bound_assumption}.  Then for $S_\mathcal C$ defined in \eqref{eq: SC} and for $\mathscr C_n = \{ \mathcal C_n\}$,  the collection $\{S_{\mathcal C} \}_{\mathcal C \in \cup_{n \in \mathbb N} \mathscr C_n}$ satisfies Condition \ref{con: a s gaussian}.
\end{theorem}

\subsection{Analysis of Cluster-Based Inference with Learned Clusters.}
\label{sec: main results}

Consider an array $\mathscr D_\shortrightarrow$ of spatially indexed datasets $\mathscr D_n = \{\zeta_i\}_{i \in \mathsf X_n}$ all distributed according to one common $\Pr_{0}$.  Consider a sequence of statistical hypotheses $H_{0\shortrightarrow} = \{ H_{0n}\}_{n \in \mathbb N}$ where each $H_{0n}$ is a restriction on only the distribution of the $n$th sample of $\mathscr D_\shortrightarrow$ (i.e., on $\mathscr D_n$).  
Sections \ref{sec: main results}-\ref{subsec: bch discussion} discuss formal properties of sequences of (random) test outcomes $\hat T_{\shortrightarrow} = (\hat T_n)_{n\in\mathbb N}$ with primary interest in $\hat T_{\shortrightarrow} = \hat T_{\bullet(\alpha)\shortrightarrow}$ arising from $\hat T_{\bullet(\alpha)}$ defined by applying Algorithm \ref{Algorithm1} to data $\mathscr D_n, \mathsf X_n$ with $\bullet =$ IM, CRS, or CCE as in Section \ref{method}. 
For any $n$, let $\mathscr C_n$ be its associated collection of partitions 
\begin{align} 
	\mathscr C_n = \{ \mathcal C : \mathcal C = \Cluster(\mathsf X_n) \text{ and } (\Test, \Cluster) \in \mathscr T_n \text{ for some } (\Test, \Cluster) \}.
\end{align}    
For each $n$, there is an element $\widehat {\texttt{Cluster}}(\mathsf X_n) \in \mathscr C_n$ which solves the optimization in equation \eqref{eq: optimization problem}.  

An important property used in the analysis of Algorithm \ref{Algorithm1} in Sections \ref{sec:IM} and \ref{sec:CRS} is that the act of selecting partitions $\hat {\mathcal C}_{n}$ has negligible effect on $\hat T_n$. This is formalized by the following condition.

\begin{assumption}\label{con: data dep partitions}
There are $\mathscr C_{1,n} \subseteq \mathscr C_n$ for all $n$ such that for $\hat T_{\bullet(\alpha) \shortrightarrow}$, \ 
(i)  $\lim_{n \rightarrow \infty } {\Pr_0}( \widehat{\mathtt{Cluster}}( \mathsf X_n)  \notin \mathscr C_{1,n}) = 0,$
 (ii) \ $\lim_{n \rightarrow \infty} \sup_{\mathcal C \in \mathscr C_{1,n}} \Big |{\Pr_0}( T_{\bullet(\alpha)} = \mathrm{ Reject } | \widehat {\mathtt{Cluster}}(\mathsf X_n) = \mathcal C) - {\Pr_0}( T_{\bullet(\alpha)} = \mathrm{ Reject } ) \Big | = 0 .$
\end{assumption}

Condition 6 is a high-level assumption that may be non-trivial in general settings.  When specializing to the context of Algorithm \ref{Algorithm1} in Section \ref{method}, Condition \ref{con: data dep partitions} is equivalently a high-level condition on the properties $\widehat{\mathrm{Err}}_{\mathrm{Type}\text{-}\mathrm{I}}$ and $\widehat{\mathrm{Err}}_{\mathrm{Type}\text{-}\mathrm{II},\mathbf{P}_{\mathrm{alt}}}$.  Note, if a single fixed $\G_0$ is used together with $k$-\texttt{medoids} clustering, i.e., if using partitions corresponding to the singleton $\{ \G_0 \}$ rather than the full set $\{ 2,..., \G_{\max} \}$, then Condition \ref{con: data dep partitions} holds.    Second, recall that the method for selecting $\hat {\mathcal C}$ in Section \ref{method} is through a combination of QMLE, simulation, and constrained optimization of weighted average power.  
Under sufficient regularity conditions, the QMLE procedure defined in the appendix is consistent for a fixed object which may depend on $n$.  If correspondingly $\hat k$ is consistent for some element of $\{2,...,\G_{\max}\}$ then Condition \ref{con: data dep partitions} holds. 
Finally, if estimates $\widehat{\mathrm{Err}}_{\mathrm{Type}\text{-}\mathrm{I}}$ and $\widehat{\mathrm{Err}}_{\mathrm{Type}\text{-}\mathrm{II},\mathbf{P}_{\mathrm{alt}}}$ in a particular example are asymptotically independent of test statistics of interest, Condition \ref{con: data dep partitions} holds.

\subsection{General analysis for IM with learned clusters.}\label{sec:IM}
\label{subsec: im theorem}

The IM procedure at level $a$ was defined in Section \ref{method} for a scalar hypothesis about a regression coefficient in a linear model of the form $H_0: \theta_0 = \theta^\circledast$, depending on a partition $\mathcal C$ of the data's spatial indexing set.  In particular, recall that in Section \ref{method}, the IM test was defined for $a \in (0,1)$ and given by 
$ T_{\text{IM}(a)} = \text{Reject} \ \ \ \ \text{if} \ \ \ |t(S)|>t_{1-a/2,\G-1},$
where $t_{1-a/2,\G-1}$ is the $(1-a/2)$-quantile of a $t$-distribution with $\G-1$ degrees of freedom and with $S = (S_\mathsf C)_{\mathsf C \in \mathcal C} \in \mathbb R^\mathcal C$, $|\mathcal C| = \G$, $S_\mathsf C = (n/\G)^{1/2} ( \hat \theta_\mathsf C - \theta^\circledast)$. 

\begin{theorem} 
	\label{im_theorem_main}
	Suppose that $\{S_\mathcal C\}_{\mathcal C \in \cup_{n\in\mathbb N} \mathscr C_n}$ satisfy Condition \ref{con: a s gaussian} under $\Pr_0$.  Suppose $\hat T_{\mathrm{IM}(\alpha)\shortrightarrow}$ satisfies Condition \ref{con: data dep partitions}.  Suppose further that $\alpha\le 2\Phi(-\sqrt{3})=0.083...$, where $\Phi$ the standard Gaussian cumulative distribution function.  Then $ \lim \sup_{n \rightarrow \infty}  {\Pr_0}(\hat T_{\mathrm{IM}(\alpha),n}  = \text{\em Reject}) \leq \alpha.$
\end{theorem}

The condition that $\alpha\le 2\Phi(-\sqrt{3})$ is also needed in \IM.  Note that the conditions above were  verified for the  OLS and IV models under regularity conditions given in the previous section with partitions generated by $k$-\texttt{medoids} with fixed $\G$.  The properties of $\widehat{\mathrm{Err}}_{\mathrm{Type}\text{-I}}$ and $\widehat{\mathrm{Err}}_{\mathrm{Type}\text{-I},\mathbf{P}_{\mathrm{alt}}}$ do not appear explicitly in the conditions of Theorem \ref{im_theorem_main}. However, their properties are reflected through Condition \ref{con: data dep partitions}. The proof of Theorem \ref{im_theorem_main} builds largely on the arguments in \IM.  Relative to \IM, Theorem \ref{im_theorem_main} requires establishing uniform limit results over tests corresponding to partitions in $\mathscr C_\shortrightarrow.$  

Theorem \ref{im_theorem_main} is helpful in understanding the statistical performance of Algorithm \ref{Algorithm1} with IM in the context of the hypothesis tests considered in the empirical application in Section \ref{empirical_application} and simulation study in Section \ref{simulation}.  In particular, consider a sequence of regression datasets $\mathscr D_\shortrightarrow$ with $\mathscr D_n = \{\zeta_i\}_{i\in\mathsf X_n}$ with $\zeta_i = (Y_i, X_i, W_i, Z_i)$ as in Section \ref{extension_to_ols} satisfying $Y_i = \theta_0 X_i + W_i'\gamma_0 + U_i$.    

\begin{corollary}($\mathrm{IM}$ and $k$-\texttt{medoids} in the linear model).  Suppose that $\mathsf X_{\shortrightarrow}$ satisfies Condition \ref{ahlfors}.  Suppose that $\mathscr D_{\shortrightarrow}$ satisfies Condition \ref{con: lin reg}.  Let $H_{0\shortrightarrow} = (H_{0,n})_{n\in \mathbb N}$ be the sequence of hypotheses $H_{0,n}: \theta_0 = \theta^{\circledast}$.  Let $\hat T_{\mathrm{IM}(\alpha)\shortrightarrow} = (\hat T_{\mathrm{IM}(\alpha),n})_{n \in \mathbb N}$ be defined by Algorithm \ref{Algorithm1} applied to $(\mathscr D_n, \mathsf X_n)$ for hypotheses $H_{0,n}$ for each $n$ with $k$-medoids and $k_{\max}$ fixed independent of $n$.  Suppose $\widehat{\mathrm{Err}}_{\mathrm{Type}\text{-}\mathrm{I}}$ and $\widehat{\mathrm{Err}}_{\mathrm{Type}\text{-}\mathrm{II},\mathbf{P}_{\mathrm{alt}}}$ are such that $\hat T_{\mathrm{IM}(\alpha)\shortrightarrow}$ satisfies Condition \ref{con: data dep partitions}. Suppose further that $\alpha\le 2\Phi(-\sqrt{3})$. Then $\lim_{n\rightarrow \infty} \Pr_0 (\hat T_{\mathrm{IM}(\alpha),n} = \ \mathrm{Reject} \ ) \leq \alpha$.
\end{corollary}

\subsection{General analysis for CRS with learned clusters.}\label{sec:CRS}
\label{subsec: crs theorem}

The procedure defining $T_{\text{CRS}(\alpha)}$ can be stated in slightly more generality.  In this more general setup, CRS relative to a partition $\mathcal C$ depends on a real-valued function $w : \mathbb R^{\mathcal C} \rightarrow \mathbb R$, and a statistic $S \in \mathbb R^{\mathcal C}$, which is chosen such that large values of $w$ provide evidence against $H_0$.  As was described in Section \ref{method}, in all implementations in this paper, $w$ is chosen to be   $w(S) = |t({S})|$ and $S$ has components $S_{\mathsf C} = (n/\G)^{1/2}(\hat \theta_\mathsf C - \theta^\circledast)$ with $\G = |\mathcal C|$.  Note, the CRS test was developed in  \CRS \ for a single fixed partition and relies on a notion of approximate symmetry. Specifically, given a partition $\mathcal C$ of $\mathsf X_n$, the CRS procedure also depends on the availability of a finite group $\mathcal H_{\mathcal C}$ of symmetries which act within component $\mathsf C$ of  $\mathcal C$.  These symmetries are formalized by a group action  $\mathscr D_n \mapsto h\mathscr D_n$ for all $h \in \mathcal H_{\mathcal C}$.    
The action on the data is defined in such a way that 
\begin{align} \mathscr D_n \mapsto h\mathscr D_n \ \text{ induces } \ S_\mathcal C \mapsto h S_\mathcal C, \end{align} 
respects the original action, i.e., $S_{\mathcal C}(h\mathscr D_n) = hS_\mathcal C (\mathscr D_n)   $.  By abuse of notation, write $w(\mathscr D_n) = w(S_{\mathcal C})$.   

In the context of cluster-based inference in this paper, the only type of group action of interest will be set of signs $\mathcal H_{\mathcal C}= \{ -1, +1 \}^{\mathcal C}$ that operates within clusters.  Elements, $h$, which are tuples of signs $h_{\mathsf C} \in \{ -1, +1\}$ indexed by $\mathsf C \in \mathcal C$ can be defined to act in a way such that $S_{\mathcal C,\mathsf C} \mapsto h_{\mathsf C} S_{\mathcal C,\mathsf C}.$
In the context of the OLS and IV models, $(\hat \theta_\mathsf C - \theta^\circledast )\mapsto h_\mathsf C( \hat \theta_\mathsf C - \theta^\circledast)$.  
To complete the description of CRS, let $M = |\mathcal H_{\mathcal C}|$ and  $w^{1}(S_\mathcal C)\leq w^{2}(S_\mathcal C)\leq$...$\leq w^{M}(S_\mathcal C)$ denote the order statistics of the orbit  $ \left\{w\left(h\mathscr D_n\right):h\in \mathcal H\right\} $.   Let $j_a=M(1-a)$,
$M^{+}(X)=|\left\{ 1\leq j\leq M:w^{j}(S_\mathcal C)>w^{{j_{a}}}S_\mathcal C\right\} |$
and $M^{0}(X)=|\left\{ 1\leq j\leq M:w^{j}(X)=w^{{j_{a}}}(S_\mathcal C)\right\} |$.  Let $\tilde a\left(S_\mathcal C\right)=\frac{Ma-M^{+}(S_\mathcal C)}{M^{0}(S_\mathcal C)}$.
Set
$ T_{\text{CRS}(a),\mathcal C} = \text{Reject}$ if $w(S_\mathcal C)>w^{j_a}(S_\mathcal C)$ or $\( w(S_\mathcal C) = w^{j_a}(S_\mathcal C) \text{ and } \tilde a(S_\mathcal C) = 1 \).$
$\hat T_{\text{CRS}(\alpha)}$ is then defined as in Section \ref{method}.

\begin{assumption}[Regularity Conditions for CRS]
	\label{regularity_crs}
(i) $\{S_\mathcal C\}_{\mathcal C \in \cup_{n\in\mathbb N} \mathscr C_n}$ satisfy Condition \ref{con: a s gaussian}.
				(ii)   $h\tilde S^*_{\mathcal C} $ has the same distribution as $\tilde S^*_\mathcal C $ for all $h \in \mathcal H_{\mathcal C}$, for all $\mathcal C$ in all $\mathscr C_n$; for distinct $h,h'$, either $w \circ h = w \circ h'$ or $\tilde \Pr(w(h\tilde S^*_{\mathcal C}) \neq w(h'\tilde S^*_{\mathcal C})) = 1$; $w$ is continuous and the action of $h$ is continuous for each $h$ in each $\mathcal H_{\mathcal C}.$ 
(iii)  $| \mathscr C_n|$ and $\max_{\mathcal C \in \mathscr C_n}|\mathcal C|$ are each bounded by a constant independent of $n$.
\end{assumption}

\begin{theorem} 
	\label{crs_theorem_main}
	Suppose that $\{S_\mathcal C\}_{\mathcal C \in \cup_{n\in\mathbb N} \mathscr C_n}$ satisfies Condition \ref{regularity_crs} under $\Pr_0$. Suppose $\hat T_{\mathrm{CRS}(\alpha)\shortrightarrow}$ satisfies Condition \ref{con: data dep partitions}.   Then $ \lim \sup_{n \rightarrow \infty}  {\Pr_0}(\hat T_{{\mathrm{CRS}(\alpha)},n}  = \text{\em Reject}) \leq \alpha.$

\end{theorem}

In \CRS, it was shown that the conclusion of Theorem \ref{crs_theorem_main} holds under conditions which are subsumed by Condition \ref{regularity_crs}, working in the case that $\mathscr C_n$ are singletons, i.e., $\mathscr C_n = \{ \mathcal C_n\}$ consist of a unique predetermined partition of $\{1,...,n\}$ without any described dependence on any spatial indexing $\mathsf X_n$.  
Condition \ref{regularity_crs}(iii) is stronger than necessary.  Discussion about weakening Condition \ref{regularity_crs}(iii) is postponed until Subsection \ref{subsec: uniformity} about uniformity.  In particular, no requirements on $|\mathscr C_n|$ are needed when $w(S) = |t(S)|$.  The supplemental material proves Theorem 7 under a strictly weaker condition, Condition \ref{regularity_crs}(iii)', which is stated formally there.  The argument accommodating the ability to generalize to this more general case is new relative to \CRS.

As was the case in the previous section with the analysis of IM, Theorem \ref{crs_theorem_main} has an important corollary which specializes to the specific case of Algorithm \ref{Algorithm1} (with $k$-\texttt{medoids}) for testing $H_0: \theta_0 = \theta^{\circledast}$, where $\theta_0$ is an unknown coefficient in a linear regression model.  Once again, consider a sequence of regression datasets $\mathscr D_\shortrightarrow$ with $\mathscr D_n = \{\zeta_i\}_{i\in\mathsf X_n}$ with $\zeta_i = (Y_i, X_i, W_i, Z_i)$ as in Subsection \ref{extension_to_ols} above satisfying $Y_i = \theta_0 X_i + W_i'\gamma_0 + U_i$.    

\begin{corollary}($\mathrm{CRS}$ and $k$-\texttt{medoids} in the linear model).  Suppose that $\mathsf X_{\shortrightarrow}$ satisfies Condition \ref{ahlfors}.  Suppose that $\mathscr D_{\shortrightarrow}$ satisfies Condition \ref{con: lin reg}.  Let $H_{0\shortrightarrow} = (H_{0,n})_{n\in \mathbb N}$ be the sequence of hypotheses $H_{0,n}: \theta_0 = \theta^{\circledast}$.  Let $\hat T_{\mathrm{CRS}(\alpha)\shortrightarrow} = (\hat T_{\mathrm{CRS}(\alpha),n})_{n \in \mathbb N}$ be defined by Algorithm \ref{Algorithm1} applied to $(\mathscr D_n, \mathsf X_n)$ for hypotheses $H_{0,n}$ for each $n$ with $k$-medoids and $k_{\max}$ fixed independent of $n$.  Suppose $\widehat{\mathrm{Err}}_{\mathrm{Type}\text{-}\mathrm{I}}$ and $\widehat{\mathrm{Err}}_{\mathrm{Type}\text{-}\mathrm{II},\mathbf{P}_{\mathrm{alt}}}$ are such that $\hat T_{\mathrm{CRS}(\alpha)\shortrightarrow}$ satisfies Condition \ref{con: data dep partitions}.  Then $\lim_{n\rightarrow \infty} \Pr_0 (\hat T_{\mathrm{IM}(\alpha),n} = \ \mathrm{Reject} \ ) \leq \alpha$.
\end{corollary}

\allowdisplaybreaks{

\appendix
	
\section{Implementation Details}
\label{Appendix A}



This section gives full implementation details for CCE, IM, CRS,  SK methods in the main text. 
Here, $n$ denotes total sample size; let $N_{\text{pan}}$ be the number of cross-sectional and $T_{\text{pan}}$ number of times time periods.  Let $Y$, $X$, $W$, $U$, $Z$, $V$ be 
$n$-row matrices obtained by stacking $Y_{de}$, $X_{de}$, $(W_{de}',1)$, $U_{de}$,$Z_{de}$ and $V_{de}$.  Let $M_A=I_n-A(A'A)^{-1}A'$ for matrices $A$. For $B\in \mathbb{R}^{n\times n}$ with rank $\ell$, let $P_B$ be a bijective linear map from the subspace orthogonal to the column span of $B$, to $\mathbb{R}^{n-\ell}$.  First, SK is the spatial-HAC estimator of \cite{Sun2015}. 
We apply the method of \cite{Sun2015} to transform our locations to a regular integer lattice. 
\cite{Sun2015} require smoothing parameters $(K_1,K_2)$. Set $K_1=K_2=\left[\sqrt{0.4\cdot N_{\text{pan}}^{2/3}} \right]$ following \cite{Lazarus2018}. The remaining methods, CCE, IM, CRS, follow a 5 step procedure.

	\noindent \textbf{\textit{Step 1. Partitions}} For $\G = 2,...,\G_{max}$, apply $k$-\texttt{medoids} with dissimilarity matrix from the data to obtain $\mathcal{C}^{(\G)}$, a collection of $\G$ partitions of the observations.
	
	\noindent \textbf{\textit{Step 2. Parametric Covariance Estimation.}} 	
		In the OLS model, let $\widehat{U}$ be the vector of full-sample least-squares estimation residuals.  Note,
		$P_{M_{[X,W]}}\widehat{U}\sim \mathrm N(0,\Sigma_{PMU} ) \ \  \text{under $U\sim \mathrm N(0,\Sigma)$}  \text{ where} \  \ \Sigma_{PMU}=P_{M_{[X,W]}}M_{[X,W]}\Sigma M_{[X,W]}P_{M_{[X,W]}}' .$
		The covariance matrix $\Sigma_{PMU}$ is made non-singular by applying the matrix $P_{M_{[X,W]}}$ to $\widehat{U}$. 
		$\Sigma$ is estimated by QMLE using an exponential covariance model with parameter $\tau=(\tau_1,\tau_2,\tau_3)'$,
$
			\Sigma_{de,d'e'}(\tau)=cov[U_{de},U_{d'e'};\tau]=\exp(\tau_1)\exp (-\tau_2^{-1} \Vert L_d-L_{d'}\Vert_2-\tau_3^{-1}|e-e'|).
$ 
		Calculate
		$\widehat{\tau}=\underset{\tau}{\arg\max}\left\lbrace \frac{1}{2}\log \det (\Sigma_{PMU}(\tau) )+\frac{1}{2}\widehat{U}'P_{M_{[X,W]}}'(\Sigma_{PMU}(\tau) )^{-1}P_{M_{[X,W]}}\widehat{U}\right\rbrace ,$
		where $\Sigma_{PMU}(\tau)=P_{M_{[X,W]}}M_{[X,W]}\Sigma(\tau) M_{[X,W]}P_{M_{[X,W]}}'$
		and $\Sigma(\tau)=(\Sigma_{de,d'e'}(\tau))_{de,d'e'}$ is the implied covariance matrix of $U$ under $\tau$.
		The covariance matrix estimator is $\Sigma(\widehat{\tau})$. 	
		In the IV model, covariance matrices for the structural and first-stage equations are estimated separately. 
		Let $\widehat{U}=M_WY-M_WX\widehat{\theta}$ and $\widehat{V}=M_WX-M_WZ\widehat{\pi}$, where $\widehat{\theta}$ is the 2SLS estimator for $\theta_0$ and $\widehat{\pi}$ is the least-square estimator for $\pi$. 
		Solve
		$\widehat{\tau}^\varepsilon=\underset{\tau}{\arg\max}\left\lbrace \frac{1}{2}\log \det (\Sigma_{PM\varepsilon}(\tau) )+\frac{1}{2}\widehat{\varepsilon}'P_{M_W}'(\Sigma_{PM\varepsilon}(\tau) )^{-1}P_{M_W}\widehat{\varepsilon}\right\rbrace ,$
		where $\Sigma_{PM\varepsilon}=P_{M_W}M_W\Sigma(\tau)M_WP_{M_W}'$, $\Sigma(\tau)=(\Sigma_{de,d'e'}(\tau))_{de,d'e'}$, $  \Sigma_{de,d'e'}(\tau)$ is as previously defined, and $\varepsilon$ is either $U$ or $V$. 
		The covariance estimators for $U$ and $V$ are $\widehat{\Sigma}_U =  \Sigma(\widehat{\tau}^U)$ and $\widehat{\Sigma}_V = \Sigma(\widehat{\tau}^V)$. Finally, estimate $\widehat{\rho}$ with the empirical correlation between $\widehat{\Sigma}_U^{-1/2}\widehat{U}$ and $\widehat{\Sigma}_V^{-1/2}\widehat{V}$.

	\noindent \textbf{\textit{Step 3. Simulate data.}} 
	Given covariance estimator(s) from \textit{Step 2}, simulate independent copies of the observable data for each $\diamond=1,...,B$ as follows for each $\theta\in \{-10/\sqrt{n},-9/\sqrt{n},...,10/\sqrt{n} \}$. Use $B = 10000,1000$ in the empirical example and simulation study.
In the OLS model, draw $U^\diamond$ from $N(0,\Sigma(\widehat{\tau})) $. 
		Reproduce data by  $Y_{de}^\diamond=\widehat{\alpha}+\theta X_{de}+W_{de}'\widehat{\gamma}+U_{de}^\diamond,$
		where $\widehat{\alpha}$ and $\widehat{\gamma}$ are full-sample least-square estimators, and $U_{de}^\diamond$ is the $de$ element on $U^\diamond$. 	
				In the IV model, draw $(U^\diamond,V^\diamond)$ such that 
		$\begin{pmatrix}
			U^\diamond \\ V^\diamond
		\end{pmatrix}\sim \mathrm N\left( 0, \begin{bmatrix}
			\widehat{\Sigma}_U & \widehat{\rho}\widehat{\Sigma}_U^{1/2}(\widehat{\Sigma}_V^{1/2})'\\
			\widehat{\rho}\widehat{\Sigma}_V^{1/2}(\widehat{\Sigma}_U^{1/2})' & \widehat{\Sigma}_V
		\end{bmatrix}\right).$
		Reproduce data by 
$
				Y_{de}^\diamond=\widehat{\alpha}+\theta X_{de}^\diamond+W_{de}'\widehat{\gamma}+U_{de}^\diamond, \ 
				X_{de}^\diamond=\widehat{\mu}+\widehat{\pi}Z_{de}+X_{de}'\widehat{\xi}+V_{de}^\diamond,
$
		where $\widehat{\alpha}$ and $\widehat{\gamma}$ are full-sample 2SLS estimators, $\widehat{\mu}$, $\widehat{\pi}$, $\widehat{\xi}$ are full-sample least-square estimators for the first-stage, and $U_{de}^\diamond$, $V_{de}^\diamond$ are the $de$ elements of $U^\diamond$, $V^\diamond$. 
		
	\noindent \textbf{\textit{Step 4. Type I and Type II estimates.}}	For $\G=2,...,\G_{\max}$ and $a \in [0,.05]$,  compute simulated Type I error rate 
	$ \widehat {\mathrm{Err}}_{\mathrm{Type-I}}(\mathrm{IM}(a), \G )$, 
	$ \widehat {\mathrm{Err}}_{\mathrm{Type-I}}(\mathrm{CRS}(a), \G)$, 
	$ \widehat {\mathrm{Err}}_{\mathrm{Type-I}}(\mathrm{CCE}(a),\G)$ by testing $H_0:\theta_0=0$ on each simulated dataset with $\theta=0$ from \textit{Step 3}.  
	Set $\hat \alpha_{\text{IM}, \G}$, $\hat \alpha_{\text{IM}, \G}$, $\hat \alpha_{\text{IM}, \G}$ to be the largest value $ a \in [0, \alpha]$ such that  
	$\widehat {\mathrm{Err}}_{\mathrm{Type-I}}(\mathrm{IM}(a),\G)\leq \alpha$, 
	$\widehat {\mathrm{Err}}_{\mathrm{Type-I}}(\mathrm{CRS}(a),\G)\leq \alpha$,
	$\widehat {\mathrm{Err}}_{\mathrm{Type-I}}(\mathrm{CCE}(a),\G)\leq \alpha$. 
	For each $\G=2,\dots, \G_{\max}$, compute simulated average Type II error rate  $ \widehat {\mathrm{Err}}_{\mathrm{Type-II}}(\mathrm{IM}(\hat a_{\text{IM}, \G}), \G)$, $ \widehat {\mathrm{Err}}_{\mathrm{Type-II}}(\mathrm{CRS}(\hat a_{\text{CRS}, \G}), \G)$, $ \widehat {\mathrm{Err}}_{\mathrm{Type-II}}(\mathrm{CCE}(\hat a_{\text{CCE}, \G}), \G)$ by testing $H_0:\theta_0=0$ for 
	$\theta \in \{-10/\sqrt{n},-9/\sqrt{n},...,10/\sqrt{n} \}$ on each dataset from \textit{Step 3} and averaging.

	\noindent \textbf{\textit{Step 5. Solution to \eqref{eq: optimization problem}.}} Set $(\widehat\alpha_{\mathrm{CCE}} , \hat \G_{\mathrm{CCE}})$, $(\widehat\alpha_{\mathrm{IM}} ,  \hat \G_{\mathrm{IM}})$, $(\widehat\alpha_{\mathrm{CRS}} , \hat \G_{\mathrm{CRS}})$ as the solution to \eqref{eq: optimization problem}.


\bibliographystyle{apalike}
\bibliography{dkbib1}
\pagebreak

\makeatletter
\setlength{\@fptop}{0pt}
\makeatother

\comment{
\begin{table}[t]
	\caption{\cite{moran} $\mathcal{I}$ Tests based on IV Scores}
	\label{table: moran}
	\renewcommand{\arraystretch}{1.65}
	\begin{threeparttable}
		\begin{tabular}{lcc}
			\hline \hline
			& Moran $\mathcal{I}$ & p-value \tabularnewline
			\hline
			Spatial weights - First election round	& 4.338 &  $<$ 0.001 \\ 
			Spatial weights - Second election round	& 2.546 &  0.011     \\
			Spatial weights - Both election rounds	& 5.387 &  $<$ 0.001  \\
			Intertemporal weights 					& 2.755 &  0.006     \\
			\hline \\
		\end{tabular} 
		\begin{tablenotes}
	\item \linespread{1.65} \normalsize {\textbf{Notes:}} This table provides \cite{moran} $\mathcal{I}$ test results for the \cite{Condra2018} example. The test statistic value and associated p-value are respectively reported in the columns labeled ``Moran $\mathcal{I}$'' and ``p-value.'' The first three rows report tests with weights for the Moran test constructed on the basis of geographic distance as described in the text. The first row uses only observations from the first election round, the second row uses only observations from the second election round, and the third row uses all observations to construct the test statistic. The final row uses a weight matrix constructed solely on the basis of intertemporal distance as described in the text. 			
		\end{tablenotes}
	\end{threeparttable}
\end{table}
}

\begin{table}[t]
	\caption{\cite{Condra2018} example}
	\label{table: moran}
	\renewcommand{\arraystretch}{1.65}
	\begin{threeparttable}
		\begin{tabular}{p{7.75cm}p{3.cm}p{3cm}}
			\hline \hline
			\multicolumn{3}{c}{A. 	\cite{moran} $\mathcal{I}$ Tests based on IV Scores} \\
			\hline
			& Moran $\mathcal{I}$ & p-value \tabularnewline
			\hline
			Spatial weights - First election round	& 4.338 &  $<$ 0.001 \\ 
			Spatial weights - Second election round	& 2.546 &  0.011     \\
			Spatial weights - Both election rounds	& 5.387 &  $<$ 0.001  \\
			Intertemporal weights 					& 2.755 &  0.006     \\
		\end{tabular}
		\begin{tabular}{lp{1.cm}p{.9cm}p{.9cm}p{1.3cm}p{1.3cm}p{1.3cm}p{1.3cm}p{.65cm}c}
			\hline \hline
			\multicolumn{10}{c}{B. Impact of Early Morning Attacks on Voter Turnout during the 2014 Election} \\
			\hline
			&\ \ \ {\footnotesize{}$ \phantom{ \Big |}\hat{\theta}_{0}$} & \ {\footnotesize{}$s.e.$}  & {\footnotesize{}$t$-stat} &\multicolumn{2}{c}{\footnotesize{}C.I.}  &\multicolumn{2}{c}{\footnotesize{}C.I. (unadjusted cr. val.)} & {\footnotesize{}${\hat \G}$} & {\footnotesize{}${\hat \alpha}$}  \tabularnewline
			\hline
			UNIT & -0.145 & 0.061  & 2.385 & ( -0.312, &\  0.022 )  & ( -0.265, &\  -0.025 )   & 205 & \\
			CCE  & -0.145  &0.090 &  1.609 &   ( -0.404, & \ 0.114 )  &   ( -0.358, & \ 0.068 )  & 8 & 0.024\\
			IM   &  -0.122 & 0.432  & 0.282 & ( -1.320, & \ 1.077 )   & ( -1.320, & \ 1.077 )  & 5  & 0.050 \\
			CRS  & -0.242 &  &  & ( -1.497, & \ 0.084 )   & ( -1.497, & \ 0.084 )  & 6 & 0.031\\
			\hline\\
		\end{tabular}  
		\begin{tablenotes}
			\vspace{-.25in}
			\item \linespread{1.65} \normalsize {\textbf{Notes:}} This table provides results for the \cite{Condra2018} example. Panel A reports \cite{moran} $\mathcal{I}$ test results. The test statistic value and associated p-value are respectively reported in the columns labeled ``Moran $\mathcal{I}$'' and ``p-value.'' The first three rows report tests with weights for the Moran test constructed on the basis of geographic distance. The final row uses a weight matrix constructed solely on the basis of intertemporal distance. Panel B provides inferential results based on selected clusters. Row labels indicate which procedure
			is used. The column labeled $\hat{\theta}_{0}$ reports the IV estimate of $\theta_{0}$ for the full sample in the rows labeled UNIT and CCE and the average of IV estimators of $\theta_{0}$ of the five or six data-generated clusters in the rows labeled IM and CRS respectively. Column $s.e.$   
			reports the estimated standard errors. Note that CRS does not rely on an explicit standard error estimate.
			Column $t$-stat reports the $t$-statistic for testing the null hypothesis that $\theta_{0} = 0$. Column C.I.
			reports confidence intervals for $\theta_{0}$ obtained using $\hat\alpha$ and $\hat{k}$ from \eqref{eq: optimization problem}. Column C.I. (unadjusted crit. val.s)
			reports confidence intervals of the IV estimate of $\theta_{0}$ using fixed asymptotic thresholds, i.e., rejecting if the p-value is less than .05.
			Note that CRS produces the same CI in both cases because of the discreteness of the test.  Column ${\hat \G}$ and $\widehat{\alpha}$ indicate the number of
			clusters and p-value threshold selected in each procedure, separately.			
		\end{tablenotes}
	\end{threeparttable}
\end{table}

\pagebreak

\begin{figure}[t]
	\centering
	\includegraphics[scale=.7,trim={2cm 8cm 2cm 8cm},clip]{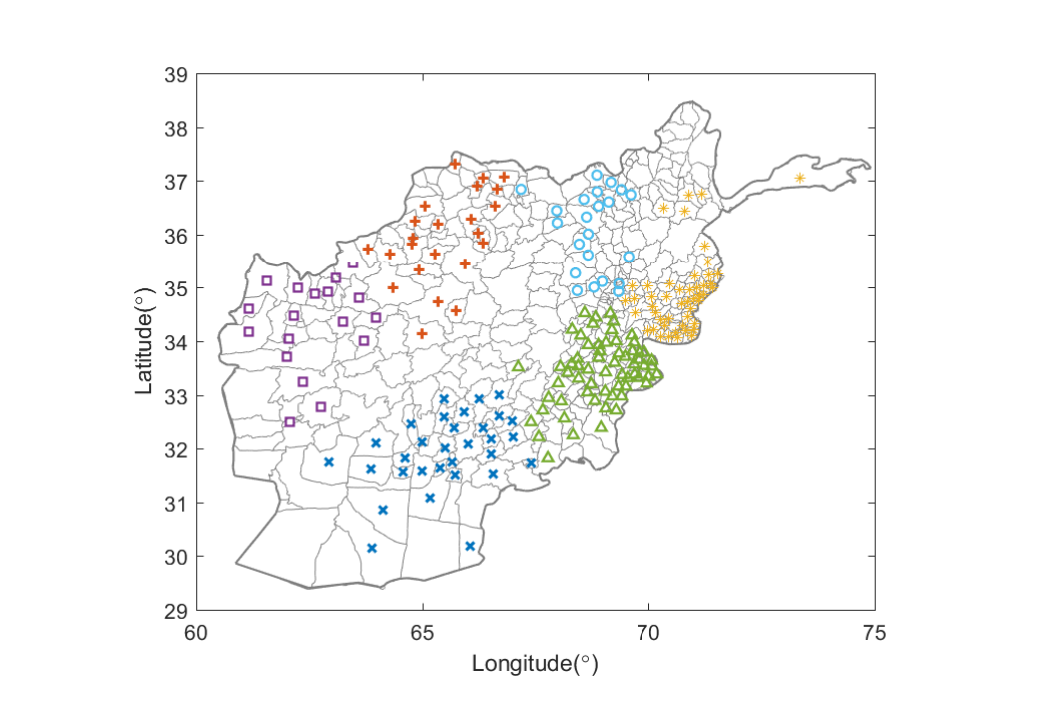}
	\caption{\linespread{1.65} \normalsize Display of partition of 205 districts in Afghanistan by $k$-\texttt{medoids} using final number of clusters given by $\G=6$. Distances are Euclidean distances based on latitude and longitude coordinates recorded at district centroids.  Different marks correspond to different clusters in the partition.  Marks are plotted at district centroids. Districts without symbols have insufficient data and were excluded in the analysis \cite{Condra2018} as well as in the current analysis.}
	\label{fig: afghanistan}
\end{figure}

\pagebreak
\clearpage

\begin{table}[t]
	\centering
	\renewcommand{\arraystretch}{1.65}
	\caption{Simulation Results}
	\label{table_summary}
	\begin{threeparttable}	
		\begin{tabular}{m{.1\textwidth} >{\centering}m{.08\textwidth} >{\centering}m{.08\textwidth} >{\centering}m{.08\textwidth} m{.1\textwidth} >{\centering}m{.08\textwidth} >{\centering}m{.08\textwidth} >{\centering\arraybackslash}m{.08\textwidth}}
			\hline
			\hline
			& \multicolumn{3}{c}{OLS} & & \multicolumn{3}{c}{IV} \\
			\cline{2-4} \cline{6-8}
			Method & Bias & RMSE & Size &  & Median Bias & MAD & Size \\
			\hline
			\multicolumn{8}{c}{A. BASELINE}\\
			\hline
			SK          & \multirow{4}{*}{0.015}  & \multirow{4}{*}{0.337} & 0.381 &  & \multirow{4}{*}{0.002}  & \multirow{4}{*}{0.114} & 0.362 \\
			UNIT-U      &        &       & 0.577 &  &        &       & 0.568 \\
			UNIT        &        &       & 0.047 &  &        &       & 0.074 \\
			CCE         &        &       & 0.046 &  &        &       & 0.062 \\
			IM          & 0.014  & 0.213 & 0.044 &  & -0.059 & 0.091 & 0.046 \\
			CRS         & 0.014  & 0.213 & 0.042 &  & -0.061 & 0.095 & 0.040 \\
			\hline
			\multicolumn{8}{c}{B. SAR}\\
			\hline
			SK          & \multirow{4}{*}{-0.013} & \multirow{4}{*}{0.866} & 0.447 &  & \multirow{4}{*}{0.002}  & \multirow{4}{*}{0.280} & 0.324 \\
			UNIT-U      &        &       & 0.639 &  &        &       & 0.502 \\
			UNIT        &        &       & 0.359 &  &        &       & 0.256 \\
			CCE         &        &       & 0.047 &  &        &       & 0.083 \\
			IM          & -0.006 & 0.385 & 0.038 &  & -0.046 & 0.158 & 0.025 \\
			CRS         & -0.003 & 0.354 & 0.049 &  & -0.095 & 0.213 & 0.041\\
			\hline 			
		\end{tabular}
		\begin{tablenotes}
			\item \linespread{1.65} \normalsize Notes: Results from the OLS Simulation (in columns labeled ``OLS'') and IV Simulation (in columns labeled ``IV'') described in Section \ref{subsec: DGP}. Row labels indicate inferential method. Panel A corresponds to the {BASELINE} design where tuning parameters for UNIT, CCE, IM, and CRS are selected based on calculating size and power from feasible estimates of a correctly specific parametric model. Panel B corresponds to the {SAR} design where tuning parameters for UNIT, CCE, IM, and CRS are selected based on calculating size and power from feasible estimates of a misspecified parametric model. For the OLS designs, we report the bias and RMSE of the point estimator associated with each procedure along with size of 5\% level tests. For the IV designs, we report the median bias and MAD of the point estimator associated with each procedure along with size of 5\% level tests.
		\end{tablenotes}
	\end{threeparttable}
\end{table}

\pagebreak

\begin{figure}[t]
	\centering
	\subfloat[BASELINE]{\includegraphics[width=.35\linewidth,trim=4.5cm 6.5cm 4.5cm 7cm]{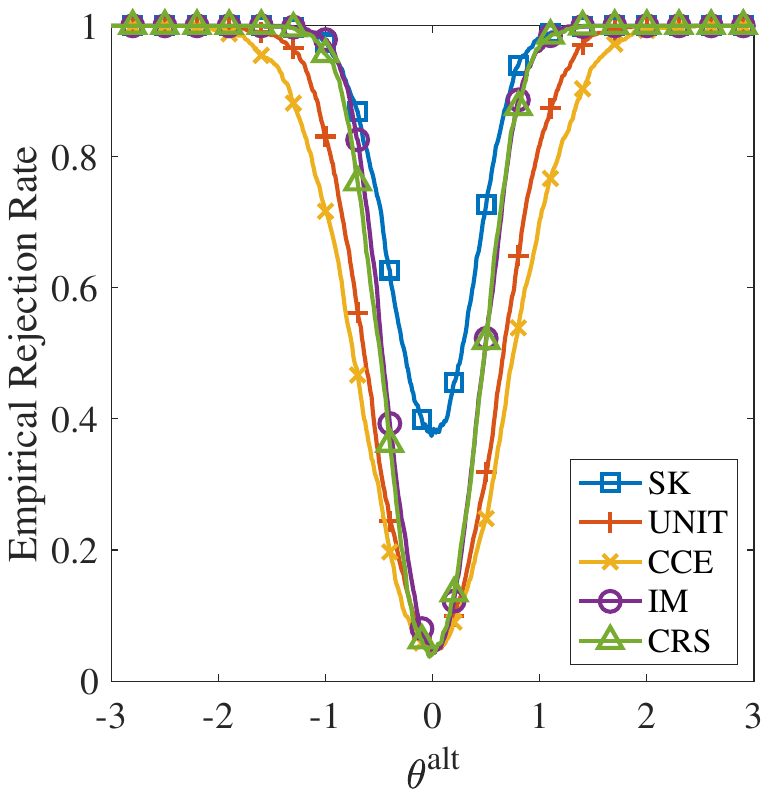}}\qquad
	\subfloat[SAR]{\includegraphics[width=.35\linewidth,trim=4.5cm 6.5cm 4.5cm 7cm]{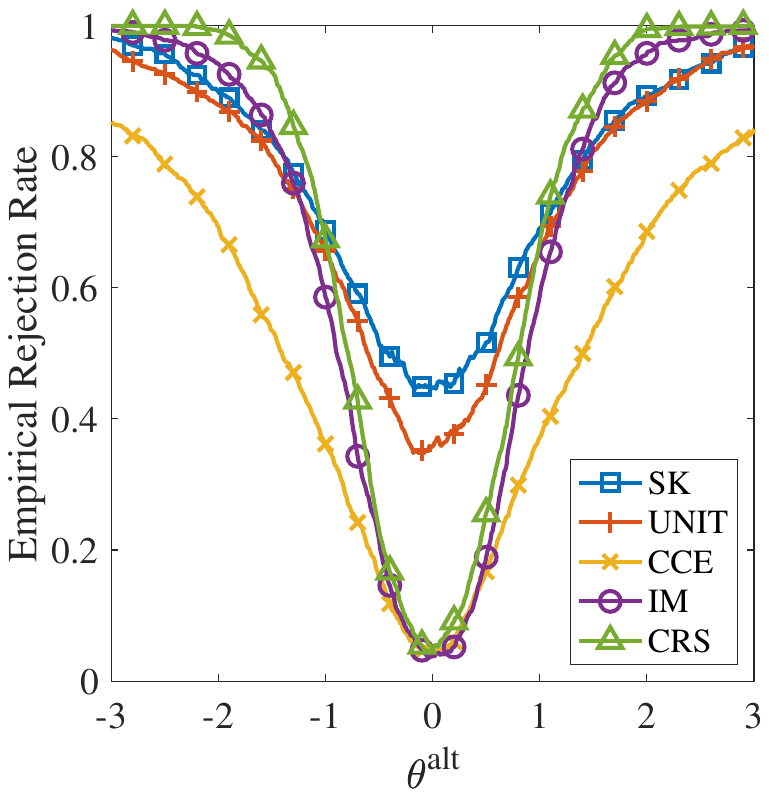}}\\
	\caption{OLS power curves.}\label{fig: ols power}
	
	\centering
	\subfloat[BASELINE]{\includegraphics[width=.35\linewidth,trim=4.5cm 6.5cm 4.5cm 7cm]{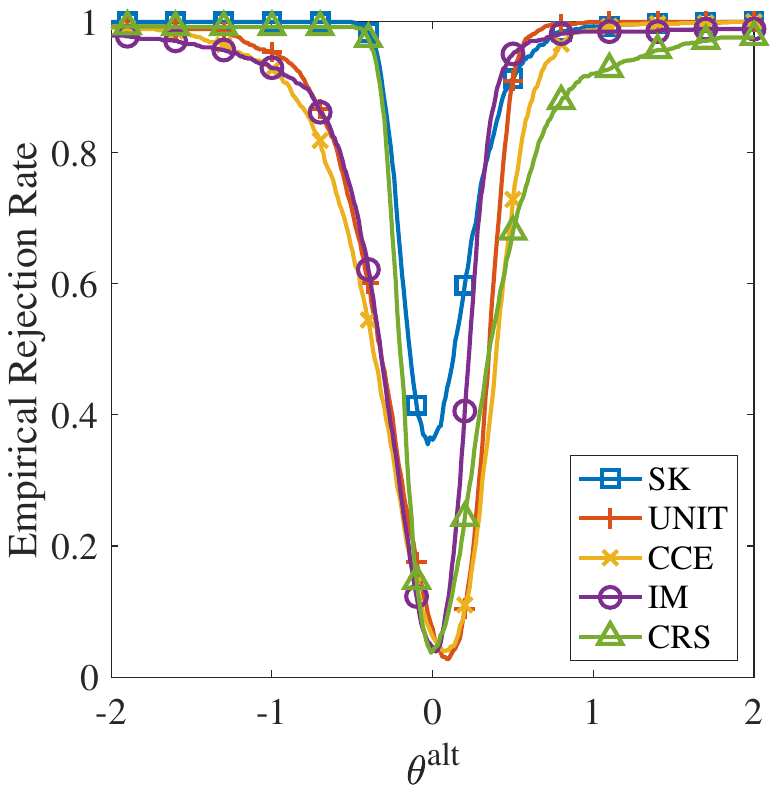}}\qquad
	\subfloat[SAR]{\includegraphics[width=.35\linewidth,trim=4.5cm 6.5cm 4.5cm 7cm]{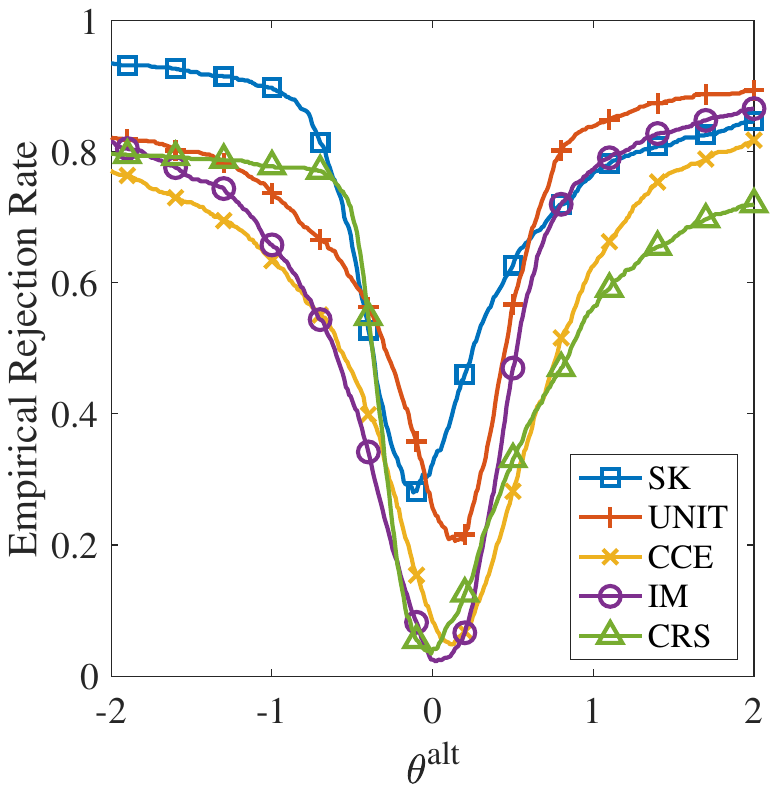}}\\
	\caption{IV power curves.}
	\label{fig: iv power}
\end{figure}

\pagebreak

\begin{table}[t]
	\caption{Distribution of $\hat\G$ and $\hat \alpha$}
	\label{tab: hat k}
	\centering
	\renewcommand{\arraystretch}{1.35}
		\begin{tabular}{  m{.05\textwidth} >{\centering}m{.08\textwidth} >{\centering}m{.08\textwidth} >{\centering}m{.08\textwidth} >{\centering}m{.08\textwidth} >{\centering}m{.08\textwidth} >{\centering}m{.08\textwidth} >{\centering}m{.08\textwidth} >{\centering\arraybackslash}m{.08\textwidth}}
			\hline
			\hline 
			\vspace{2mm}
			 \multirow{16}{*}{\rotatebox{90}{$\hat k$ Distribution}} & \multicolumn{1}{l}{$\mathrm{Pr} \ \hat \G=$} & 2  & 3 & 4 & 5 & 6 & 7 & 8\\
			\cline{2-9}
			\multicolumn{8}{c}{A. OLS - BASELINE}\\
			\cline{2-9}
			&CCE                  & 0.000          & 0.000 & 0.001 & 0.098 & 0.311 & 0.260 & 0.330 \\
			&IM                   & 0.000          & 0.000 & 0.000 & 0.007 & 0.018 & 0.113 & 0.862 \\
			&CRS                  & 0.000          & 0.000 & 0.000 & 0.000 & 0.001 & 0.099 & 0.900 \\
			\cline{2-9}
			\multicolumn{8}{c}{B. OLS - SAR}\\
			\cline{2-9}
			&CCE                  & 0.000          & 0.001 & 0.022 & 0.051 & 0.113 & 0.268 & 0.545 \\
			&IM                   & 0.000          & 0.003 & 0.069 & 0.554 & 0.040 & 0.176 & 0.158 \\
			&CRS                  & 0.000          & 0.000 & 0.000 & 0.000 & 0.001 & 0.401 & 0.598 \\
			\cline{2-9}
			\multicolumn{8}{c}{C. IV - BASELINE}\\
			\cline{2-9}
			&CCE                  & 0.000          & 0.001 & 0.001 & 0.039 & 0.242 & 0.261 & 0.456 \\
			&IM                   & 0.000          & 0.000 & 0.005 & 0.051 & 0.062 & 0.066 & 0.816 \\
			&CRS                  & 0.000          & 0.000 & 0.000 & 0.000 & 0.000 & 0.157 & 0.843 \\
			\cline{2-9}
			\multicolumn{8}{c}{D. IV - SAR}\\
			\cline{2-9}
			& CCE                  & 0.006          & 0.027 & 0.048 & 0.047 & 0.071 & 0.257 & 0.544 \\
			& IM                   & 0.005          & 0.078 & 0.240 & 0.568 & 0.030 & 0.038 & 0.041 \\
			& CRS                  & 0.000          & 0.000 & 0.000 & 0.000 & 0.032 & 0.480 & 0.488 \\
			\cline{2-9}
		\end{tabular}
			\begin{tabular}{m{.05\textwidth} >{\centering}m{.11\textwidth} >{\centering}m{.11\textwidth} >{\centering}m{.13\textwidth} >{\centering}m{.12\textwidth} >{\centering}m{.12\textwidth} >{\centering\arraybackslash}m{.1\textwidth}}
			\vspace{3mm}		
			 \multirow{16}{*}{\rotatebox{90}{$\hat \alpha \phantom{\Big |^{\Big |}} $ Distribution}} &  \multicolumn{1}{l}{$\hat \alpha$ quantile =} &  0.1  &  0.25   & 0.5    & 0.75   & 0.9    \\
			\cline{2-7}
			\multicolumn{7}{c}{A. OLS - BASELINE} \\
			\cline{2-7}
			&CCE                  & 0.004          & 0.006 & 0.008 & 0.010 & 0.012 \\
			&IM                   & 0.039          & 0.042 & 0.046 & 0.050 & 0.050 \\
			&CRS                  & 0.039          & 0.043 & 0.047 & 0.050 & 0.050 \\
			\cline{2-7}
			\multicolumn{7}{c}{B. OLS - SAR}\\
			\cline{2-7}
			&CCE                  & 0.021          & 0.024 & 0.027 & 0.031 & 0.034 \\
			&IM                   & 0.047          & 0.049 & 0.050 & 0.050 & 0.050 \\
			&CRS                  & 0.047          & 0.047 & 0.050 & 0.050 & 0.050 \\
			\cline{2-7}
			\multicolumn{7}{c}{C. IV - BASELINE}\\
			\cline{2-7}
			&CCE                  & 0.004          & 0.005 & 0.007 & 0.008 & 0.010 \\
			&IM                   & 0.045          & 0.048 & 0.050 & 0.050 & 0.050 \\
			&CRS                  & 0.039          & 0.039 & 0.047 & 0.050 & 0.050 \\
			\cline{2-7}
			\multicolumn{7}{c}{D. IV - SAR}\\
			\cline{2-7}
			&CCE                  & 0.014          & 0.020 & 0.025 & 0.029 & 0.033 \\
			&IM                   & 0.032          & 0.050 & 0.050 & 0.050 & 0.050 \\
			&CRS                  & 0.008          & 0.043 & 0.050 & 0.050 & 0.050\\
			\cline{1-7}
		\end{tabular}
		\begin{tablenotes}
			\item \linespread{1.65} \normalsize Notes: Simulation results for $(\hat k, \hat \alpha)$, in the optimization problem \eqref{eq: optimization problem}. Panels correspond to the different designs in Section \ref{subsec: DGP}.  Entries for the $\hat k$ distribution report simulation probability of $\hat \G$ being equal to the value given in the column label.  Entries for the $\hat \alpha$ distribution report .1, .25, .5, .75,  .9 quantiles of the simulation distribution of $\hat\alpha$. 
		\end{tablenotes}
\end{table}

\pagebreak

\clearpage	
\pagebreak

	\setcounter{page}{1}
	
	\section*{{SUPPLEMENT TO: ``Inference for Dependent Data with Learned Clusters''}}

\noindent
Authors:  Jianfei Cao, Christian Hansen, Damian Kozbur, Lucciano Villacorta

\noindent \today

\

\appendix

This supplement contains proofs, additional simulation results, additional discussion of optimization details including run times, and further general discussion about assumptions made in the theory and empirical calculations in the main text.

\

\textbf{Table of Contents for Supplemental Material}
{\singlespacing
\begin{align*} 
&\text{S.1.  Proofs of Theorems \ref{clt_mean}--{\ref{iv_theorem}}} & \pageref{Proof_1_5}\\
&\text{S.2.  Proof of Theorem \ref{im_theorem_main}} & \pageref{Proof_IM}\\
&\text{S.3.  Proofs of Theorem \ref{crs_theorem_main}} & \pageref{Proof_CRS}\\
&\text{S.4.  Proof of Proposition \ref{w proposition}.} & \pageref{Proof_Prop_2}\\
&\text{S.5.  Additional Discussion on the Choice of $k_{\max}$} & \pageref{choice:kmax}\\
&\text{S.6.  Bounds on Ball Cardinalities in Three Example Metric Spaces.} & \pageref{Ball_cardinalities}\\
&\text{S.7.  Discussion for BCH with learned clusters.} & \pageref{subsec: bch discussion}\\
&\text{S.8.  Discussion on Uniformity.} & \pageref{subsec: uniformity}\\
& \text{S.9.  Computation Times.} & \pageref{computation:times}\\
& \text{S.10.  Computation Times.} & \pageref{sec:opt:contours}\\
&\text{S.11.  Additional Simulation Results.} & \pageref{Additional:sim}\\
&\text{References for Supplemental Material.} & \pageref{Supp:references}\\
&\text{Tables and Figures for Supplemental Material} & \pageref{Supp:Tables} 
\end{align*}}
\setcounter{section}{19}
\setcounter{subsection}{0}
\subsection{Proofs of Theorems \ref{clt_mean}--\ref{crs_theorem_main}.}\label{Proof_1_5}

Before proving Theorems \ref{clt_mean}--\ref{crs_theorem_main}, the basic structural results about sequences of Ahlfors spaces are gathered.  The first preliminary result is due to \cite{assoud:diss} and ensures the existence of certain embedding of metric spaces of bounded doubling into Euclidean spaces.

\begin{theorem}[Assoud] Let $(\mathsf X,\mathrm d)$ be an arbitrary (not necessarily finite) metric space with finite doubling dimension $\text{dim}_{2\times}(\mathsf X) < \infty$.  Let $\eps \in (0,1)$.  Then there exists an  $L$-bi-Lipschitz map $(\mathsf X,\mathrm d^{1-\eps}) \rightarrow \mathbb R^\nu$ for some $L,\nu$ which depend only on $\eps$ and $\text{dim}_{2\times}(\mathsf X) $.
\end{theorem}

Ahlfors regularity of $\mathsf X$ implies bounded doubling as a result of standard arguments involving covering numbers of metric spaces.  Thus the next proposition follows.

\begin{proposition}  Suppose that $ \mathsf X \in \mathbf{Ahlf}_{C,\delta}$.  Then the doubling dimension of $\mathsf X$ is bounded by $\text{dim}_{2\times} (\mathsf X) \leq \delta \log_2 (3C^2)$.  Therefore, there exists an $L$-bi-Lipschitz map $(\mathsf X,\mathrm d^{3/4}) \rightarrow \mathbb R^\nu$ where  $\nu$ and the Lipschitz constant $L$ depend only on $C, \delta$. 
\end{proposition}

\subsubsection{Proof of Theorem \ref{clt_mean}.}

\begin{proof}
	
In the case that $\mathsf X_n$ do not embed isometrically, let  $L$ be the bi-Lipschitz constant from the maps $(\mathsf X_n,\mathrm d_n^{1-{1/4}}) \rightarrow \tilde {\mathsf X}_n \subseteq \mathbb R^{\nu}$.  The array $\{\{\zeta_{i}\}_{i\in \mathsf X_n}\}_{n=1}^\infty $ indexed on $\mathsf X_n$ yields an array $\{\{\tilde \zeta_{i}\}_{i \in \tilde{\mathsf X}_n}\}_{n=1}^\infty $ indexed on $\tilde {\mathsf X}_n$.  It is sufficient to check the conditions of Corollary 1 in \cite{Jenish2009} for this new process.  Apply the same set array of constants $c_{i}$ to $\tilde \zeta_{i}$.
Assumption 1 in \cite{Jenish2009} is satisfied by the fact that distances are at least $\rho_0$ in $\mathsf X_n$ for some $\rho_0>0$ by Ahlfors regularity.  Note that $L$ depends only on $C, \delta$ and in particular, does not change with $n$.  Then $\forall i,j,$ $\tilde{\mathrm  d}_n(i,j) $ is also bounded away from zero by a constant which does not depend on  $n$.   Condition \ref{mixing_assumption}(ii) is identical to Equation 3 in \cite{Jenish2009}.  

The next conditions in \cite{Jenish2009} are mixing conditions. To verify these, let $\tilde \alpha^{\mathrm{mix}}_{g , l,n}(r)$ and $\bar {\tilde  \alpha}^{\mathrm{mix}}_{g, l,n}(r)$ denote the corresponding mixing coefficients for $\tilde \zeta_{i}$ over $\tilde {\mathsf X}_n$.  Note that $\tilde {\mathrm d}(\mathsf U,\mathsf V) \geq r \Rightarrow  \mathrm d(\mathsf U,\mathsf V) \geq L^{-1}\nu^{-1/2}r^{\frac{3}{4}}$.  Let $c = L^{-1} \nu^{-1/2} $. 
Then $\bar{\tilde \alpha}^{\mathrm{mix}}_{g, l}(r) \leq\bar \alpha^{\mathrm{mix}}_{g,l}(c r^{3/4})$.
To verify Equation 4 in \cite{Jenish2009}, it is sufficient to show that $\sum_{m=1}^\infty \bar {\tilde \alpha}^{\mathrm{mix}}_{1,1} ( m) m^{\nu \times \frac{\mu+2}{\mu}-1} < \infty.$  
Next note that $\bar{\tilde \alpha}^{\mathrm{mix}}_{1,1}(m)$ is nonincreasing and defined for nonnegative real $m$ and $m^{\nu \times \frac{\mu+2}{\mu}-1}$ is a polynomial in $m$.  Thus, the above summation is bounded up to a constant by the corresponding integral  $ \int_{m=0}^\infty \bar {\tilde \alpha}^{\mathrm{mix}}_{1,1} ( m) m^{\nu \times \frac{\mu+2}{\mu}-1}dm$.
Substituting $\bar {\tilde \alpha}^{\mathrm{mix}}_{1,1}(m) \leq \bar \alpha^{\mathrm{mix}}_{1,1}(cm^{3/4})$ and a standard calculus change of variables $m' = cm^{3/4}$, $dm' = \frac{3}{4}cm^{-1/4} dm$ shows that it is sufficient to verify 
   $ \int_{m'=0}^\infty \bar { \alpha}^{\mathrm{mix}}_{1,1} (m') {m'}^{\frac{4}{3} ( \nu \times \frac{\mu+2}{\mu}-1)} {m'}^{1/3} dm'<\infty.$
This integral is in turn bounded by a constant times the summation $$\sum_{m'=1}^\infty \bar { \alpha}^{\mathrm{mix}}_{1,1} ( m') {m'}^{\frac{4}{3} \nu \times \frac{\mu+2}{\mu}-1} = \sum_{m'=1}^\infty \bar { \alpha}_{1,1}^{\mathrm{mix}} ( m') {m'}^{\nu \times \frac{\mu+2}{\mu}-1}.$$  This is seen to be finite after noting results in \cite{Jenish2009} (see discussion regarding their Equation 4 on their page 4), that $\sum_{m=1}^\infty \bar \alpha_{1,1}(m) m^{\nu \times \frac{\mu + 2}\mu - 1} < \infty$ is weaker than $\sum_{m=1}^\infty m^{\nu-1} \bar \alpha_{1,1} (m)^{\mu/(2+\mu)} < \infty$, and in turn weaker than Condition 2(iii) which assumes $\sum_{m=1}^\infty m^{\nu} \bar \alpha_{1,1} (m)^{\mu/(2+\mu)} < \infty$.  This verifies Equation 4 in \cite{Jenish2009}.  Using similar arguments, Condition {2}(iv) implies Assumption 4(2) in \cite{Jenish2009}.  Next,  Condition \ref{mixing_assumption}(v) implies that $$ \bar {\tilde \alpha}^{\mathrm{mix}}_{1,\infty}(m) \leq  \bar { \alpha}^{\mathrm{mix}}_{1,\infty}(cm^{3/4}) = O( (cm^{3/4})^{\frac{4}{3} (\nu - \mu)}) = O(  m^{\nu- \mu}  )  $$ thus verifying Assumption 4(3) in \cite{Jenish2009}.  
Finally, Condition \ref{mixing_assumption}(vi) implies Assumption 5 in \cite{Jenish2009}.  This verifies the assumptions of Corollary 1 in \cite{Jenish2009}.

\end{proof}
\subsubsection{Proof of Theorem \ref{asymptotic_independence}.}

\begin{proof}
Let $\mathsf C , \mathsf D \subseteq \mathsf X_n$ be any two partition components $\mathsf C, \mathsf D \in \mathcal C$.  Then,

$$\Big | \text{cov}\Big (\sigma({\mathsf C})^{-1} \sum_{i\in \mathsf C} \zeta_i, \sigma({\mathsf D})^{-1}\sum_{j \in \mathsf D}\zeta_j \Big ) \Big| \leq \sigma({\mathsf C})^{-1}\sigma({\mathsf D})^{-1} \sum_{(i,j) \in \mathsf C \times \mathsf D} \bar \alpha_{1,1} ( \lceil \mathrm  d_n(i,j) \rceil)^{\mu/(2+\mu)} $$
$$ =   \sigma({\mathsf C})^{-1}\sigma({\mathsf D})^{-1}    \sum_{m=1}^\infty | \{(i,j) \in \mathsf C \times \mathsf D : m-1 \leq \mathrm  d_n(i,j) <m\} | \bar \alpha_{1,1}(m) ^{\mu/(2+\mu)}. $$
The above reductions appear also in \BCH (specifically, in the proof of their Lemma 1 on page 149, their quantity $R_L$ being exactly the same as the covariance object printed immediately below), which previously also appeared in \cite{Jenish2009} and \cite{bolthausen1982}, 
  
Next, to proceed in the current setting, note that by Condition \ref{ahlfors},
$$| \{(i,j) \in \mathsf C \times \mathsf D : m-1 \leq \mathrm d_n(i,j) < m\} |  \leq | \{(i,j) \in \mathsf X \times \mathsf X : m-1 \leq \mathrm  d_n(i,j) <m\} | \leq nCm^\delta.$$  By Condition {3}, it also follows that 
$$ \max_{m \leq \bar r } | \{(i,j) \in \mathsf C \times \mathsf D : m-1 \leq \mathrm d_n(i,j) <m\} |  \leq o(n).$$ 
For ease of notation, let $\bar \alpha_{1,1} = \bar \alpha_{1,1}^{\text{mix}}$  and $\bar \alpha_{1,\infty} = \bar \alpha_{1,\infty}^{\text{mix}}$. Recall that by Condition 2(ii) and 2(v), $$\sum_{m=1}^\infty \bar \alpha_{1,1}(m) m^{\nu \times \frac{\mu + 2}\mu - 1} < \infty \ \text{and} \ \bar \alpha_{1,\infty} = O(m^{-\nu - (4/3)\mu}).$$  
Then the original covariance is bounded by 
$$ \sigma({\mathsf C})^{-1} \sigma({\mathsf D})^{-1} \[ \sum_{m=1}^{\bar r} o(n)\bar \alpha_{1,1}(m)^{\mu/(2+\mu)} +  \sum_{m=\bar r}^\infty  nCm^\delta  \bar \alpha_{1,1}(m)^{\mu/(2+\mu)} \] $$ $$ \ \ \ = n \sigma({\mathsf C})^{-1} \sigma({\mathsf D})^{-1} \[ \sum_{m=1}^{\bar r} o(1)\bar \alpha_{1,1}(m)^{\mu/(2+\mu)}  +  \sum_{m=\bar r}^\infty  Cm^\delta  \bar \alpha_{1,1}(m)^{\mu/(2+\mu)} \] .$$

\noindent
Working inside the above square brackets 
gives  $\sum_{m=1}^{\bar r} o(1) \bar \alpha_{1,1}(m)^{\mu/(2+\mu)} \leq o(1) \sum_{m=1}^\infty \bar \alpha_{1,1}(m)^{\mu/(2+\mu)} = o(1)$ as well as  $\sum_{m=1}^\infty  Cm^\delta  \bar \alpha_{1,1}(m)^{\mu/(2+\mu)} < \infty$ by Condition 2(iii), which with $\bar r \rightarrow \infty $ implies $  \sum_{m=\bar r}^\infty  Cm^\delta  \bar \alpha_{1,1}(m)^{\mu/(2+\mu)} =o(1)$.   
It is only left to bound $ n\sigma({\mathsf C})^{-1} \sigma({\mathsf D})^{-1}$.  
By Condition 2(vi), $|\mathsf C| | \mathsf D| \sigma({\mathsf C})^{-2} \sigma({\mathsf D})^{-2} = O(1)$.  Similarly, by Condition 3(i) (i.e., balance of cluster sizes), $ | \mathsf C | |\mathsf C| / n^2 = O(1)$.  Together these give $n \sigma(\mathsf C)^{-1} \sigma(\mathsf D)^{-1} = O(1)$.  The Theorem follows. 
\end{proof}

\subsubsection{Proof of Theorem \ref{km_prop}.}

\begin{proof}

Let $\bar r_n =  \log  n$.  From this point on in this proof, $n$ is excluded from notation when used as an index in a subscript and when it provides no added clarity.  
Consider two arbitrary (sequences of) points $i,j \in \mathsf X$.  
Let $M = \{ l : | \mathrm d(i,l) - \mathrm d(j,l)| \leq \bar r \}.$   
Let $M_0$ be an $\bar r^2$-net of $M$.  Note that by construction, if $|M| \neq o(n)$ then $|M_0| \neq o(n/\bar r^{2\delta})$.
Let $A \subseteq [\frac 1 {2}, \hspace{.5mm}  \frac 9 {10}]$ satisfy $|a - a'| \geq \bar r^3/\mathrm d(i,j)$ for each $a,a' \in A$.   Take $|A | \geq \frac{4}{10} \frac{\mathrm d(i,j)}{2}/ \bar r^3$.    
For $a \in A,$ let $M_a$ consist of interpolants $l_a$ such that  $ | \mathrm d(i,l_a) - a \mathrm d(i ,l) | \leq K$ and $ | \mathrm d(l_a,l) - (1-a) \mathrm d(i,l) | \leq  K$ for each $l\in M_0$.  Then by trigonometry, $\cup_{a \in A} M_a$ is $3\bar r$-separated for $n$ sufficiently large and contains $|A| \times |M_0|$ elements.  This is shown by constructing the line segments $\overline{\iota(i) \iota(l)}$ and $\overline{\iota(i) \iota(l')}$ where $\iota$ is the coarse isometry to Euclidean space $\mathsf E$, which gives points $u, u'$ belonging to the above constructed line segments with distances $d_\mathsf E( \iota (l_a), u)$, $d_\mathsf E( \iota (l_a'), u')$ bounded, implying that $u,u'$ are sufficiently separated to yield the claim.  By $3 \bar r$-separation and $| \textsf{B}_{\bar r}(l_a)|  \geq C^{-1} \bar r^{\delta}$, it follows that 
$$\left | \bigcup_{a \in A, l_a \in M_a} \textsf{B}_{\bar r}(l_a)   \right | \geq  |A| \times |M_0| C^{-1} \bar r^{\delta} \geq \frac{4}{20} \frac{\mathrm d(i,j)}{ \bar r^3}  |M_0|  C^{-1} \bar r^{\delta} .$$

As the left hand side above is  $\leq n$, it follows that the two statements (1.) $d(i,j) \geq \bar r^{4 + \delta}$ for $n$ sufficiently large, and (2.) $|M_0| \neq o(n/{\bar r}^{2\delta})$ cannot be true simultaneously.  

Next, it is shown that $\G$-\texttt{medoids} terminates with medoids $i_1,...,i_\G$ such that $\mathrm d(i_g, i_l) \geq ( \log n)^{4 + \delta}$ for $n$ sufficiently large.  This implies the small boundaries statement of Theorem 3, as a violation of this statement would require that $M$ corresponding to some pair $i_g, i_l$ satisfies $M \neq o(n)$, subsequently implying $M_0 \neq o(n/\bar r^{2\delta})$, violating the above dichotomy. 

Again for contradiction, suppose there is a sequence $\ell \rightarrow 0$ such that for infinitely many $n$, there are two clusters $\mathsf C_1, \mathsf C_2$ with medoids $i_1, i_2$ satisfying $\mathrm d(i_1, i_2) < \ell n^{1/\delta}$.   By the pigeonhole principal, there must be a cluster $\mathsf C_3$ with $n/ \G$ members.  Then $diam(\mathsf C_3)$ must be at least $C^{-1}(n/\G)^{1/\delta} $.  Let $i_3$ be the corresponding medoid.  Then there must be $i_3' \in \mathsf C_3$ such that $\mathrm d(i_3, i_3') \geq \frac{1}{4}C^{-1}(n/\G)^{1/\delta} $ and $\mathrm d(i_3', i_g) \geq \frac{1}{4}C^{-1}(n/\G)^{1/\delta} $ for any other medoid $i_g$.   Then consider the update in the partitioned medoid algorithm given by $i_2 \leftarrow i_3'$.  This update is cost reducing for $n$ sufficiently large.  To see this, note that for elements, $i \in \textsf{B}_{\frac{1}{4} C^{-1} (n/\G)^{1/\delta}}(i_3')$ the total cost reduction from being reassigned from a medoid centered around $i_2$ to a medoid centered around $i_3'$ is at least $|\textsf{B}_{\frac{1}{4} C^{-1} (n/\G)^{1/\delta}}(i_3')|{\frac{1}{4} C^{-1} (n/\G)^{1/\delta}} \geq {\frac{1}{4} C^{-1} (n/\G)^{1/\delta}} C^{-1} ({\frac{1}{4} C^{-1} (n/\G)^{1/\delta}})^{\delta}$ .  The total cost increase from reassigning elements in $\mathsf C_2$ to $\mathsf C_1$ is at most  $\ell n^{1/\delta} | \mathsf C_2| \leq \ell n^{1/\delta} n$. The difference between the above two quantities is a lower bound on the cost reduction for the update.  Comparing the above to quantities for $n$ sufficiently large, the $\G$-\texttt{medoids} algorithm could not have stopped at a step with $\mathrm d(i_1,i_2) < \ell n^{1/\delta}$ giving the desired contradiction.

Finally, note that $\mathrm d(i_g, i_l )\geq \ell  n^{1/\delta}$ for all medoids $i_g, i_l$, some $\ell$ bounded uniformly away from $0$, and for $n$ sufficiently large implies the balanced clusters condition after applying Ahlfors regularity.
\end{proof}
\subsubsection{Proof of Theorem \ref{ASREPversion}.}

\begin{proof} Fix $k \leq k_{\max}$ and consider the subsequence of $\mathcal C_{\shortrightarrow}$ for which $|\mathcal C_n|=k$.  Let $V_k$ be a $k$-dimensional real vector space with distinguished basis indexed by $\mathcal C_{\infty,k}$.  Choose bijections $\mathcal C_n \rightarrow \mathcal C_{\infty,k}$ and induce linear maps $\ell_n: \mathbb R^{\mathcal C_n} \rightarrow V $ such that the image of $S_n$ in the direction of any vector corresponding to an element of $\mathcal C_{\infty}$ has variance 1.  Applying Theorems 1 and 2 to $S_n$, the Cram\`er-Wold device to $\ell(S_n)$ and the Almost-Sure Representation theorem (Theorem 2.19 in \cite{vdV}) then gives the $\breve S^*_{\infty,k}$ and $\breve S_{\mathcal C_n}$ such that $\| \breve S^*_{\infty,k}-\breve S_{\mathcal C_n}\|_2 \rightarrow 0$, $\tilde \Pr_k$- a.s. for some $\tilde \Pr_k$ for all $n$ in the subsequence.  Set, for $n$ in the subsequence, $\tilde S_{\mathcal C_n}^* =\ell_n^{-1} \tilde S^*_{\infty,k}$ and $\tilde S_{\mathcal C_n} = \ell_n^{-1} \breve S_{\mathcal C_n}$.  Repeat this construction for all $k$ and take all $\tilde \Pr_k$ to be a common $\tilde \Pr$.  

It remains only to show that required variances are suitably bounded above and away from zero.  Let $\zeta_i$ belong to an array satisfying Condition 1.
Note that Condition 1 assumes accompanying $\underline c \geq c_i$ 
for all required $c_i$, 
then Condition 1(\textit{vi}) 
gives that $\inf_{n \geq 1} \inf_{\mathsf C \subseteq \mathsf X_n} \sigma(\mathsf C)^2 /|\mathsf C|  > 0$. 
Also $\sigma(\mathsf C)^2 /|\mathsf C|$ is bounded above as discussed in the proof of Theorem 2. Thus, with $k$ fixed, for all $\mathsf C \in \mathsf X_n$ for some $n$, using the balanced clusters condition in Condition 3, it follows that $\lim \inf_{n \rightarrow \infty} (n/k)^{1/2} / |\mathsf C|/\sigma(\mathsf C)^2>0.$  A corresponding bound on the limit superior also holds.   Specifically, as argued in the proof of Theorem 2, 
$| \{(i,j) \in \mathsf C \times \mathsf C : m-1 \leq \mathrm d_n(i,j) < m\} | \leq nCm^\delta.$  
Thus $\n^{-1} \sigma(\mathsf C)^{-2}$
$ \leq n \sum_{m=1}^\infty  Cm^{\delta} ({\bar \alpha_{1,1}^{\text{mix}}})^{\mu/(\mu+2)}$ which is $O(n)$ by Condition 2(iii).

  \end{proof}

\subsubsection{Proof of Theorem \ref{iv_theorem}.}
\begin{proof}  

Let $\bar X_i = (X_i, W_i')$.  Theorem 4 applies to all quantities of the form $A_i B_i - \Ep[A_i B_i]$ for $A,B$ ranging over $ Z_{ij}, \bar X_{ij'}, U_i$.  Recall that by Condition 5, $\Ep[A_i B_i]$ do not depend on $i$.  
Let $a = e_1  \mathrm M^{-1} $, where $e_1$ is a vector with $1$ at the component corresponding to the $X_i$ and $0$ at components corresponding to components of the $W_i$.  Let $\zeta^a_i = \sum_{j}a_j Z_{ij} U_i$.  Let also $\zeta_{ij} = Z_{ij}U_i$.  By assumed invertibility of $\mathrm M$, Theorem 4 applies to $\zeta_i^a$ and to each $\zeta_{ij}$. 
For any $\mathsf C$, let $\hat M_{\mathsf C}$ be the matrix with components  $\frac{1}{|\mathsf C| } \sum_{i \in \mathsf C}Z_{ij} \bar X_{ij'}$.   Let $\hat a_{\mathsf C} = e_1 \hat {\mathrm M_{\mathsf C}}^{-1}$, (and $\hat a_{\mathsf C} = 0$ if $\hat M_{\mathsf C}$ is not invertible).  Let $\hat Q_{\mathsf C}$ be a vector with components $\frac{1}{|\mathsf C|} \sum_{i \in \mathsf C} Z_{ij}U_i$.  
Note
$(n/k)^{-1/2} (\hat \theta_{\mathsf C} - \theta_0) = (n/k)^{-1/2} e_1'\hat M_{\mathsf C} ^{-1} \hat Q_{\mathsf C} = (n/k)^{-1/2} \sum_{i \in \mathsf C} \zeta_i^{\hat a} =(n/k)^{-1/2} \sum_{i \in \mathsf C}  \zeta_i^a  + r_{\mathsf C} $, where $r$ is a remainder term of the form $ (n/k)^{-1/2}\sum_{i \in \mathsf C} \sum_{j} (\hat a_{\mathsf C,j} - a_j) Z_{ij} U_i.$  By invertibility of $\mathrm M$, Continuous Mapping Theorem, and applying Theorem 4 to sequences furnished with terms $A_i B_i - \Ep[A_i B_i]$ it holds that $\max_{\mathsf C \in \mathcal C_n} | r_\mathsf C|$ vanishes in probability.  Note this temporarily abandons the probability spaces $\tilde \Omega$ constructed in the proof of Theorem 4.  Next, reapply the arguments in Theorem 4 to recover a new almost sure representation for sequences furnished with $(n/k)^{-1/2} \sum_{i \in \mathsf C}  \zeta_i^a  + r_{\mathsf C} $.  This completes the proof.
\end{proof}

\subsection{Proof of Theorem \ref{im_theorem_main}.}\label{Proof_IM}

\begin{proof}
Define the function $\psi$ by $\psi(S)  = \textbf{1}_{|t(S)| > t_{1-\alpha/2,\G-1}}$, where $S$ has $\G$ components.  Let $\tilde \Ep$ be expectation with respect to the measure $\tilde \Pr$.  Note that for a partition $\mathcal C$, $T_{\mathrm{IM}(\alpha),\mathcal C} = \psi(S_\mathcal C)$. By construction of $ \tilde S_{\mathcal{C}}^*$ and Theorem 1 of \cite{Ibragimov2010}, for any $n$ and any $\mathcal C \in \mathscr C_n$,
	$\tilde \Ep[\psi(\tilde S_{\mathcal C}^*)]\le \alpha.$
	
	
Next, show $ \sup_{\mathcal C\in \mathscr C_n}|\tilde \Ep[\psi(\tilde{S}_{\mathcal C})]-\tilde \Ep[ \psi(\tilde {S}_{\mathcal C}^*)]|\rightarrow 0.$  Define $\mathrm{sd}(\cdot)$ to be the sample standard deviation of a vector $(\cdot)$, where the sample is taken to correspond to the set of components of $(\cdot)$.  Note that $\mathrm{sd}$ corresponds to the denominator of $t$. Define the events $\mathcal E_{1,\mathcal C} = \{ | t(\tilde S^*_\mathcal C) - t_{1-\alpha/2, |\mathcal C|-1} | > b_\mathcal C\}, \ \ \mathcal E_{2,\mathcal C} = \{ | \mathrm{sd}(\tilde S^*_\mathcal C) | > b_\mathcal C'\}$
	for $b_\mathcal C , b_\mathcal C' \in \mathbb R$ specified below.    Note that by the upper and lower bounds on the variances of $\tilde S^*_\mathcal C$ in Condition \ref{regularity_crs},  and continuity properties of $t(\cdot)$ and $\mathrm{sd}(\cdot)$ one can select a common positive function, denoted $M$ which is independent of $\mathcal C$, such that $\lim_{b\rightarrow 0}M(b,b')=\lim_{b'\rightarrow 0}M(b,b')=0$ and such that $\tilde \Pr( \mathcal E_{1,\mathcal C} \cap \mathcal E_{2,\mathcal C}) \geq 1 - M(b_\mathcal C, b_\mathcal C')$ whenever $b_{\mathcal C}, b_{\mathcal C}'$ are sufficiently small.  Then on $\mathcal E_{1,\mathcal C}$, a sufficient condition for $\psi(\tilde S_\mathcal C)  = \psi(\tilde S^*_\mathcal C)$ is that $|t(\tilde S_\mathcal C) - t(\tilde S^*_\mathcal C)| < b_\mathcal C.$    By $\sup_{\mathcal C \in \mathscr C_n}\| \tilde S^*_\mathcal C - \tilde S_\mathcal C\|_2 \rightarrow 0$, and by continuity properties of $t(\cdot)$ outside of $\mathcal E_{2,\mathcal C}$, there are $b_\mathcal C < b_n\rightarrow 0, b_\mathcal C' < b_n' \rightarrow 0$ such that $|t(\tilde S_\mathcal C) - t(\tilde S^*_\mathcal C)| < b_\mathcal C$ is satisfied on $\mathcal E_{1,\mathcal C} \cap \mathcal E_{2,\mathcal C}$.  Therefore,  $|\tilde \Ep[\psi(\tilde{S}_{\mathcal C})]-\tilde \Ep[ \psi(\tilde {S}_{\mathcal C}^*)]| \leq  M(b_n,b_n')$, for all $\mathcal C \in \mathscr C_n$ from which the above desired uniform limit result follows.     
	
Finally,
$	\Pr_0( \hat T_{\text{IM}(\alpha),n} = \text{ Reject }	)= \Pr_0(  \texttt{Test}_{\text{IM}(\hat a)}{(\mathscr D, \hat {\mathcal C} )}= \text{ Reject}	) \leq  \Pr_0(  \texttt{Test}_{\text{IM}(\alpha)}{(\mathscr D, \hat {\mathcal C} )}   = \text{ Reject }	)$
$ = \sum_{\mathcal C \in \mathscr C_n} \Pr_0(   \texttt{Test}_{\text{IM}(\alpha)}{(\mathscr D,  {\mathcal C} )}  = \text{ Reject } | \hat {\mathcal C} = \mathcal C	) \Pr_0(\hat {\mathcal C} =  {\mathcal C}) .$
Bound the above expression by $ o(1) + \sum_{\mathcal C \in \mathscr C_{0,n}} \Pr_0(  \texttt{Test}_{\text{IM}(\alpha)}{(\mathscr D,  {\mathcal C} )} = \text{ Reject } ) \Pr_0(\hat {\mathcal C} =  {\mathcal C}) \leq o(1) + \alpha$  by using both Condition \ref{con: data dep partitions}(i) and Condition \ref{con: data dep partitions}(ii) and the previous bounds.  This concludes the proof of Theorem \ref{clt_mean}.
\end{proof}

\subsection{Proof of Theorem \ref{crs_theorem_main}.}\label{Proof_CRS}

Before giving the proof, note that the following more general condition is used in place of Condition \ref{regularity_crs}(iii) in the main text.  

\

\setcounter{assumption}{6}
\begin{assumption}\label{Condition 7 extended}
\noindent {(iii)'}  \textit{ Either
 \begin{itemize}
 \item[a.]  $| \mathscr C_n|$ and $\max_{\mathcal C \in \mathscr C_n}|\mathcal C|$ are each bounded by a constant independent of $n$; or \
 \item[b.] For any sequence $\check \G_n \rightarrow \infty$ sufficiently slowly, the following hold. There is a sequence $\check \delta_n \rightarrow 0$ and sets $\mathcal A_{\mathcal C}\subseteq \mathbb R^\mathcal C$ for each $\mathcal C \in \mathscr C_n$ with $|\mathcal C| > \check \G_n$ which are closed under the action of $\mathcal H_\mathcal C$ such that,  $\lim_{n\rightarrow \infty} \inf_{\mathcal C \in \mathscr C_n, |\mathcal C| > \check \G} \tilde \Pr(\mathcal A_\mathcal C) = 1 $, $w(\cdot)$ can be renormalized to have Lipschitz constant 1 respect to the Euclidean norm on all $\mathcal A_{\mathcal C}$,  and it holds that  $\sup_{\mathcal C \in \mathscr C_n} \sup_{r \in \mathbb R} \Pr_\mathcal C( |r - w(h\tilde S^*_\mathcal C)| < \check \delta_n )  \rightarrow 0$ where $\Pr_\mathcal C$ is the uniform probability measure over $h \in \mathcal H_\mathcal C$.  Finally,  $| \{ \mathcal C \in \mathscr C_{n} : | \mathcal C | \leq \check \G_n \}|$ depends only $\check \G_n$.
 \end{itemize}}	 
\end{assumption}

The next proposition verifies that the test statistic, $w$, defined in the text above and used for the empirical example and simulation study satisfies the anticoncentration requirement of Condition \ref{Condition 7 extended}(iii)'.

\begin{proposition}\label{w proposition}
For $S \in \mathbb R^{\mathcal C}$, let $w( S ) = |t(S)|$.  Let $S^*_\mathcal C$ be a Gaussian random variables with independent components with variances bounded from below and above (as in Condition \ref{regularity_crs}(i)).  Suppose that $\mathscr C_n = \{ \mathcal C^{(2)} ,... , \mathcal C^{(\G_{\mathrm{max},n})} \}$ where $\G_{\mathrm{max},n} $ may be an arbitrary function of $n$. Then Condition \ref{Condition 7 extended}(iii)'\textrm{\em a} is satisfied if $\G_{\mathrm{max},n}$ is bounded and Condition \ref{Condition 7 extended}(iii)'\textrm{\em b} is satisfied if $\G_{\mathrm{max},n} \rightarrow \infty$.
\end{proposition}

The proof below is given assuming Condition \ref{Condition 7 extended}(iii)'.

\begin{proof}
	Define functions $\phi(S)$ and $\tilde{\phi}(S,U)$ for $S\in\mathbb{R}^{\mathcal{C}}$ and $U\in \mathbb{R}$ according to
	$$
	\phi(S)=
	\begin{cases}
1, \text{\qquad \ if }w(S)>w^{({j_{\alpha}})}(S)\\
\tilde a(S), \text{\quad if }w(S)=w^{({j_{\alpha}})}(S)\\
0, \text{\qquad \ if }w(S)<w^{({j_{\alpha}})}(S)\\
	\end{cases}
$$
	and define $\tilde \phi(S,U) = \phi(S)$ if $W(S) \neq W^{({j_{\alpha}})}(S)$ and $\textbf 1_{U < a(S)}$ otherwise.

Without loss of generality, assume that $\tilde \Omega$ has a uniformly distributed random variable $\tilde U \in [0,1]$ independent of all previously defined random variables defined on it.  By Theorem 2.1 in \cite{Canay2017}, using the fact that $h(\tilde{S}^*_{\mathcal C}) =_d  \tilde{S}^*_{\mathcal C}$ for all $h$ by the fact that $h$ acts as component-wise multiplication by signs, it follows that $\Ep_0[\phi(\tilde{S}_{\mathcal{C}}^*)]=\alpha$. Note, 
	$$\Ep_0[\phi({S}_{\mathcal{C}})]-\alpha =\tilde \Ep[\tilde{\phi}(\tilde{S}_{\mathcal{C}},{\tilde{U}})-\tilde{\phi}(\tilde{S}_{\mathcal{C}}^*,{\tilde{U}})].$$

Therefore, it is sufficient to show  $\limsup_{n\rightarrow \infty }\sup_{\mathcal{C}\in \mathscr C_n}\tilde \Ep[\tilde{\phi}(\tilde{S}_{\mathcal{C}},{\tilde{U}})-\tilde{\phi}(\tilde{S}_{\mathcal{C}}^*,{\tilde{U}})]=0 $.	
	Consider an integer $\check \G$, and decompose $\mathscr C_{n} = \mathscr C_{n, \check \G} \cup \mathscr C_{n, \check \G}^c$ where $\mathscr C_{n,\check \G}  = \{ \mathcal C \in \mathscr C_n : | \mathcal C |\leq \check \G \}$ and $\mathscr C_{n,\check \G}^c =\{ \mathcal C\in \mathscr C_n : | \mathcal C | > \check \G \}$, which are handled separately.  
	For each $\mathcal C \in \mathscr C_{n,\check \G}$, let $\Ev_{\mathcal{C}}$ be the event in which the order statistics of $\{w(h\tilde{S}_{\mathcal{C}}):h\in \mathcal{H}_\mathcal{C} \}$ and  $\{w(h\tilde{S}_{\mathcal{C}}^*):h\in \mathcal{H}_\mathcal{C}\}$ correspond to the same transformations $h^{(1)},\dots,h^{(|\mathcal{H}_\mathcal{C}|)}$. Then the previous expression can be controlled by
	\begin{align*}
	|\tilde \Ep[\tilde{\phi}(\tilde{S}_{\mathcal{C}},{\tilde{U}})-\tilde{\phi}(\tilde{S}_{\mathcal{C}}^*,{\tilde{U}})]|
	=& |\tilde \Ep[\tilde{\phi}(\tilde{S}_{\mathcal{C}},{\tilde{U}}){\mathbf{1}}_{\Ev_{\mathcal C}}+\tilde{\phi}(\tilde{S}_{\mathcal{C}},{\tilde{U}}){\mathbf{1}}_{\Ev_{\mathcal C}^c}-\tilde{\phi}(\tilde{S}_{\mathcal{C}}^*,{\tilde{U}}){\mathbf{1}}_{\Ev_{\mathcal C}}-\tilde{\phi}(\tilde{S}_{\mathcal{C}}^*,{\tilde{U}}){\mathbf{1}}_{\Ev_{\mathcal C}^c}]|\\
	=&|\tilde \Ep[\tilde{\phi}(\tilde{S}_{\mathcal{C}},{\tilde{U}}){\mathbf{1}}_{\Ev_{\mathcal C}^c}-\tilde{\phi}(\tilde{S}_{\mathcal{C}}^*,{\tilde{U}}){\mathbf{1}}_{\Ev_{\mathcal C}^c}]|\\
	=&|\tilde \Ep[(\tilde{\phi}(\tilde{S}_{\mathcal{C}},{\tilde{U}})-\tilde{\phi}(\tilde{S}_{\mathcal{C}}^*,{\tilde{U}})){\mathbf{1}}_{\Ev_{\mathcal C}^c}]|
	\le 2\tilde \Ep[{\mathbf{1}}_{\Ev_{\mathcal C}^c}].
	\end{align*}
	
Furthermore, $\sup_{\mathcal{C}\in  \mathscr C_{n,\check \G} }\tilde \Ep[{\mathbf{1}}_{\Ev_{\mathcal C}^c}]$ can be controlled using Fatou's Lemma by
	$$\limsup_{n\rightarrow \infty }\sup_{\mathcal{C}\in \mathscr C_{n,\check \G}}\tilde \Ep[{\mathbf{1}}_{\Ev_{\mathcal C}^c}]  \le \limsup_{n\rightarrow \infty }\tilde \Ep[\sup_{\mathcal{C}\in  \mathscr C_{n,\check \G}}{\mathbf{1}}_{\Ev_{\mathcal C}^c}] \leq \tilde \Ep[ \limsup_{n\rightarrow \infty }\sup_{\mathcal{C}\in \mathscr C_{n,\check \G}}{\mathbf{1}}_{\Ev_{\mathcal C}^c}].$$

	To further bound the right-hand side of the above expression, consider any fixed $\omega \in \tilde \Omega$ which satisfies the conditions  $\sup_{\mathcal{C}\in \mathscr C_n}\Vert \tilde{S}_{\mathcal{C}}(\omega)-\tilde{S}_{\mathcal{C}}^*(\omega)\Vert_2 \rightarrow 0$ and  $\inf_{\mathcal{C}\in \mathscr C_{n,\check \G}}\inf_{h \nsim h' \in \mathcal H_{\mathcal C} }|w_{\mathcal C}(h(\tilde{S}_{\mathcal{C}}^*)(\omega))-w_{\mathcal C}(h'(\tilde{S}_{\mathcal{C}}^*)(\omega))|>\delta_{\omega}>0$
	 for some $\delta_{\omega}>0$.  Here, the equivalence $\sim$ is defined by $h \sim h'$ whenever $w \circ h = w \circ h'$.    By Condition \ref{regularity_crs}, the set of such $\omega\in \tilde \Omega$ has $\tilde \Pr$-measure 1.  Note, more explicitly, choosing such $\delta_\omega>0$ is possible due to (1) the continuity of the function $w$ and the continuity of the action of $h$ outside a set of $\tilde P$-measure 0, and (2) the fact that $| \mathscr C_{n,\check \G}|$ depends only on $\check \G$.

	Let $h^{(1)}_{\mathcal{C}}(\omega),\dots,h^{(|\mathcal{H}_\mathcal{C}|)}_{\mathcal{C}}(\omega)$ be such that 
	$w(h^{(1)}_{\mathcal{C}}(\omega)\tilde{S}_{\mathcal{C}}^*(\omega))\leq \dots\leq w(h^{(|\mathcal{H}_\mathcal{C}|)}_{\mathcal{C}}(\omega)\tilde{S}_{\mathcal{C}}^*(\omega)).$
	Then the expression
	\begin{align*} 
	w(h^{(j+1)}_{n,\mathcal{C}}(\omega)\tilde{S}_{\mathcal{C}}(\omega))-w(h^{({j)}}_{n,\mathcal{C}}(\omega)\tilde{S}_{\mathcal{C}}(\omega))
	=& \ w(h^{(j+1)}_{n,\mathcal{C}}(\omega)\tilde{S}_{\mathcal{C}}(\omega))-w(h^{(j+1)}_{n,\mathcal{C}}(\omega)\tilde{S}_{\mathcal{C}}^*(\omega))\\
	&+w(h^{(j+1)}_{n,\mathcal{C}}(\omega)\tilde{S}^*_{n,\mathcal{C}}(\omega))-w(h^{({j)}}_{n,\mathcal{C}}(\omega)\tilde{S}^*_{n,\mathcal{C}}(\omega))\\
	&+w(h^{({j)}}_{n,\mathcal{C}}(\omega)\tilde{S}^*_{n,\mathcal{C}}(\omega))-w(h^{({j)}}_{n,\mathcal{C}}(\omega)\tilde{S}_{\mathcal{C}}(\omega))
	\end{align*}
	is nonnegative and is furthermore strictly positive for at least one $j$ unless $h \sim h'$ for all $h, h' \in \mathcal H_{\mathcal C}$ in which case the theorem holds trivially.  Positivity of the above expression in the case of non-equivalent-under-$\sim$ action follows from noting that the first and the third terms are smaller than $\delta_\omega/2$ in absolute value for $n$ sufficiently large, while the second term is greater than $\delta_\omega$.

	The claim $\limsup_{n\rightarrow \infty} \sup_{\mathcal{C}\in \mathscr C_n} {\mathbf{1}}_{\Ev_{\mathcal C}^c}(\omega) = 0$ $\tilde \Pr$-a.s. follows.  
	Thus, $\limsup_{n\rightarrow \infty }\sup_{\mathcal{C}\in  \mathscr C_{n,\check \G}}\tilde \Ep[{\mathbf{1}}_{\Ev_{\mathcal C}^c}] =0 $.  Therefore, there is a sequence $\check \G_n \rightarrow \infty$ sufficiently slowly such that $\limsup_{n\rightarrow \infty }\sup_{\mathcal{C}\in  \mathscr C_{n,\check \G_n}}\tilde \Ep[{\mathbf{1}}_{\Ev_{\mathcal C}^c}] =0 $.
	
	Next consider the partitions in the complement $\mathcal C \in  \mathscr C_{\check \G_n}^c$. 
	Let $\mathcal A_s\subseteq  \mathbb R^\mathcal C$ be of the form $\mathcal A_S = \{ hS: h \in \mathcal H_\mathcal C\}$ for some $S \in \mathbb R^\mathcal C$.  Then $w^{j_{\alpha}}$ is constant over $\mathcal A_S$ and thus $w^{j_{\alpha}}(\mathcal A_S)$ is well defined. Then for any $\delta>0$,%
$$\tilde \Pr( |w^{j_{\alpha}}(\mathcal A_S)- w(\tilde S^*_\mathcal C)| < \delta | \tilde S^*_\mathcal C \in \mathcal A_S)  \leq \sup_{r \in \mathbb R} \tilde \Pr( |r - w(\tilde S^*_\mathcal C)| < \delta | \tilde S^*_\mathcal C \in \mathcal A_S).$$  Let $\Pr_{\mathcal C}$ be the uniform probability over $h \in \mathcal H_{\mathcal C}$.  Then the above expression is bounded by	$$p_{\mathcal C, S}(\delta):=\sup_{r \in \mathbb R} \Pr_\mathcal C (|r-w(h{S})|<\delta). $$
	
Next, for any $\mathcal A_\mathcal C$ which is a disjoint union of sets of the form $\mathcal A_S$.  Then $$\tilde \Pr(|w^{j_{\alpha}}(\tilde S^*_\mathcal C) - w(\tilde S^*_\mathcal C)| \geq \delta ) \geq \tilde \Pr(|w^{j_{\alpha}}(\tilde S^*_\mathcal C) - w(\tilde S^*_\mathcal C)| \geq \delta | \tilde S^*_\mathcal C \in \mathcal A_\mathcal C) \tilde \Pr(\mathcal A_\mathcal C) $$ $$ \geq (1 -  \sup p_{\mathcal C,S}(\delta)) \tilde \Pr(\mathcal A_\mathcal C) $$
where the supremum runs over $S$ indexing $\mathcal A_S\subseteq \mathcal A_\mathcal C$.

Note that a sufficient condition for $\tilde \phi(\tilde S^*_\mathcal C,\tilde U) = \tilde \phi(\tilde S_\mathcal C,\tilde U)$ is that $\mathrm{sign}(	w^{j_{\alpha}}(\tilde S^*_\mathcal C) - w(\tilde S^*_\mathcal C) ) = \mathrm{sign}(w^{j_{\alpha}}(\tilde S_\mathcal C) - w(\tilde S_\mathcal C))$.  A further sufficient condition for this is that for some $\delta>0$, $|w^{j_{\alpha}}(\tilde S^*_\mathcal C) - w(\tilde S^*_\mathcal C)| \geq \delta$, $|w^{j_{\alpha}}(\tilde S^*_\mathcal C) - w^{j_{\alpha}}(\tilde S_\mathcal C)| < \delta/2 $, and $|w(\tilde S^*_\mathcal C) - w(\tilde S_\mathcal C)| < \delta/2 $.  By Condition \ref{Condition 7 extended}(iii)', $w$ has Lipschitz constant 1 with respect to Euclidean norm on $\mathbb R^{\mathcal C}$ on a suitably chosen sequence of events of probability approaching 1.  Then $|w^{j_{\alpha}}(\tilde S^*_\mathcal C) - w^{j_{\alpha}}(\tilde S_\mathcal C)| < \delta/2 $ whenever $ \| \tilde S^*_\mathcal C - \tilde S_\mathcal C \|_2 < \delta/2$ in which case it also holds that $|w(\tilde S^*_\mathcal C) - w(\tilde S_\mathcal C)| < \delta/2 $.  A sequence $\check \delta_n$ may be chosen such that $\check \delta_n \rightarrow 0$ sufficiently slowly such that for $n$ sufficiently large, each of the above three inequalities may be achieved with common bounds on an event with probability $1-o_{\tilde \Pr}(1)$.   

Finally, using the same reasoning as at the conclusion of Theorem \ref{im_theorem_main}, and invoking Condition \ref{con: data dep partitions}, $
	\Pr_0( \hat T_{\text{CRS}(\alpha),n} = \text{ Reject }	) 
 =  o(1) + \sum_{\mathcal C \in \mathscr C_{0,n}} \Pr_0(   \texttt{Test}_{\text{CRS}( \hat \alpha)}{(\mathscr D,  {\mathcal C} )} = \text{ Reject } ) \Pr_0(\hat {\mathcal C} =  {\mathcal C}) \leq o(1) + \alpha.
$ concluding the proof of Theorem \ref{crs_theorem_main}.
	\end{proof}

\subsection{Proof of Proposition \ref{w proposition}.}\label{Proof_Prop_2}
\begin{proof}
Fix $S \in \mathbb R^{\mathcal C}$.  Let $D$ be a constant chosen later, and consider a rescaled version of $w$ defined by $w(S) = |\mathcal C|^{-1/2} | \bar S | / (D \mathrm{sd}(S) )$.  	Next bound the following anticoncentration quantity:

	$$\sup_{r \in \mathbb R}  \Pr_{\mathcal C} ( w(hS) \in [r ,r+ \delta] ).$$

	Because $\mathcal H_{\mathcal C}$ is finite, the range of the supremum above can be restricted and the above quantity is further reduced by
	$$ = \sup_{r \in [ \min_{h\in \mathcal H_{\mathcal C}} w(hS), \max_{h\in \mathcal H_{\mathcal C}} w(hS)]} \Pr_{\mathcal C} ( w(hS) \in [r,r+ \delta])$$
	$$ \leq \sup_{r \in [ \min_{h\in \mathcal H_{\mathcal C}} w(hS), \max_{h\in \mathcal H_{\mathcal C}} w(hS)]} \Pr_{\mathcal C} ( |\bar{hS} |\in [\min_{h\in \mathcal H_{\mathcal C}}aD\mathrm{sd}(hS),\max_{h \in \mathcal H_{\mathcal C}} (a+ \delta)D \mathrm{sd}(hS) ])$$
	
	Note that the length of the interval inside the $\Pr_{\mathcal C}$ never exceeds  $$\ell = \ell(\delta, D, S):= \sup_{r \in [ \min_{h\in \mathcal H_{\mathcal C}} w(hS), \max_{h\in \mathcal H_{\mathcal C}} w(hS)]} \max_{h \in \mathcal H_{\mathcal C}} (r+ \delta)D \mathrm{sd}(hS) -  \min_{h\in \mathcal H_{\mathcal C}}rD\mathrm{sd}(hS).$$
Therefore, 
$\sup_{r \in \mathbb R} \Pr_{\mathcal C} ( w(hS) \in [r,r+ \delta]) \leq \sup_{r \in \mathbb R} \Pr_{\mathcal C} ( |\mathcal C|^{-1/2} | \bar {h S} | \in [r,r+\ell]).$
 By considering both branches of the absolute value function, (which gives an additional factor of 2), and expressing the above quantity more explicitly in terms of components $h_{\mathsf c}$, it follows that
$$\sup_{r \in \mathbb R} \Pr_\mathcal C( |r - |\mathcal C|^{-1/2}|  \sum_{\mathsf C \in \mathcal C} h_{\mathsf C} S_{\mathcal C})| |<  \ell ) \leq 2 \sup_{r \in \mathbb R}  \Pr_\mathcal C( |r - |\mathcal C|^{-1/2} \sum_{\mathsf C \in \mathcal C} h_{\mathsf C} S_{\mathcal C})| <  \ell).  $$
	
	By Corollary 2.9 in \cite{RUDELSON2008600}, which states an anticoncentration bound for bernoulli sums,	
	$$2 \sup_{r \in \mathbb R}  \Pr_\mathcal C(|r-  {| \mathcal C|}^{-1/2}\sum_{\mathsf C \in \mathcal C}h_\mathsf C S_\mathsf C| <  \ell ) \leq 2 \[\sqrt{\frac 2 \pi}\frac  { \ell} { \| |\mathcal C|^{-1/2}S\|_2} + C_1B\(\frac{ \| |\mathcal C|^{-1/2}S\|_3}{ \| |\mathcal C|^{-1/2}S\|_2}\)^3  \]$$
	where $C_1$ is an absolute constant and $B=1$ is the third moment of a Bernoulli random variable. 
	
	For a sequence $\check \G \rightarrow \infty$, Gaussian concentration properties coupled with the assumed bounds on the variances of the components of $\tilde S^*_{\mathcal C}$ imply that there is are sets $\mathcal A_{\mathcal C}$ closed under the action of $\mathcal H_{\mathcal C}$ with $\tilde \Pr(\mathcal A_{\mathcal C}) \rightarrow_{\check \G \rightarrow \infty} 1$ and a fixed choice $D$ which makes $w$ Lipschitz with constant 1 on all of $\mathcal A_{\mathcal C}$, such that $\ell$ is bounded by a fixed constant times $\delta$ for all $S$ in all $\mathcal A_{\mathcal C}$. Gaussian concentration properties for $ \| |\mathcal C|^{-1/2}\tilde S^*_\mathcal C\|_2$ and $(\| |\mathcal C|^{-1/2}\tilde S^*_\mathcal C\|_3/\| |\mathcal C|^{-1/2}\tilde S^*_\mathcal C\|_2)^3$ then also allow the construction of  $\check \delta_n \rightarrow 0$ such that the right hand of the above expression side also converges to 0 provided that $\check \G$ grows sufficiently slowly.	
\end{proof}
\subsection{Additional Discussion on the choice of $\G_{\mathrm{max}}$.}\label{choice:kmax}	
In practice, $\G_{\mathrm{max}}$ needs to be specified. As mentioned in the text, our recommendation is to set $\G_{\mathrm{max}}$ to a small number. First, note that smaller $\G$ generally allows for more dependence in the data. $\G = 1$ would be equivalent to allowing arbitrary dependence between all observations in which case inference would be as if there were one observation. IM and CCE are also undefined in this case due to division by $\G - 1$ and CRS trivially has no power against any hypothesis. In the case of fixed $\G$, larger $\G$ then heuristically allows for less dependence but also provides more power. However, the returns to power for increasing $\G$ diminish rapidly as $\G$ increases for usual significance levels. This diminishing return can easily be heuristically be gauged by considering how critical values from $t_{\G-1}$ distributions behave as $\G$ increases. Setting small $\G_{\mathrm{max}}$ thus provides more robustness to dependence in the most optimistic scenario - i.e. where $\hat{\G} = \G_{\mathrm{max}}$ - at the cost of potentially limiting power. 

A second consideration is that the our preferred procedures, IM and CRS, rely on constructing within cluster estimates. At a minimum, constructing such estimates relies on having sufficient observations per cluster for the within-cluster estimates to be well-defined. Of course, one would ideally have a large number of observations per cluster. A quick benchmark is to then look at $n/\G_{\mathrm{max}}$ which is the number of observations that would be available if all groups were equal-sized and one were at the boundary and set $\G_{\mathrm{max}}$ small enough that a researcher feels comfortable with this number of observations. 

Finally, note that we wish to offer a procedure that uniformly controls size in a world where one allows for moderate dependence. Independence would suggest using the largest possible number of groups $\G = n$. However, no procedure which allows for $\G = n$ or for a large $\G$ and does not rule out dependence \textit{a priori} can deliver uniform size control without imposing the condition that dependence is either extremely large or exactly 0 as there are always local scenarios where dependence is too moderate to be reliably statistically detected but would result in procedures that assume no dependence to have large size distortions. As we are in a setting where dependence is a possibility, we need to restrict ourselves to procedures that would control size under moderate dependence, again requiring $\G_{\max}$ not be too large. 

We illustrate the results for two different choices of $\G_{\mathrm{max}}$ in the empirical example from Section \ref{empirical_application} in Table \ref{tab: empirical results 2}. For our preferred procedures (IM and CRS), the results are identical. The results are also qualitatively robust for CCE.

\subsection{Bounds on ball cardinalities in three example metric spaces.}\label{Ball_cardinalities}

This section calculates explicit bounds for $|\mathsf B_{r}(i)|$ for three examples of metric spaces.  Results are plotted in Figure \ref{fig: ball growth}.

The first panel shows results for the \cite{Condra2018} example in which locations are geographic and given by latitude and longitude.  In the second panel, we plot the same quantities for a square subset of the integer lattice with side length $\lceil \sqrt{205} \rceil=15$ (with 205 being the number of districts observed in the \cite{Condra2018} sample). In the third panel, we repeat the calculation for locations in \cite{Conley:Dupor:Sectoral}.  \cite{Conley:Dupor:Sectoral} estimate firm-specific sectoral production functions and study the sectoral correlation in productivity shocks.  They define the location of each sector according to its share vectors computed from the input-output tables, and define distance by Euclidean distance between the shares. Like in  \cite{Condra2018}, locations in \cite{Conley:Dupor:Sectoral} are treated as exogenous. 
Points illustrated in each panel correspond to pairs $(i,r)$ with $i \in \mathsf X$.  Each point has (horizontal, vertical) axis coordinates equal to $(r, \ | \mathsf B_r(i)  |)$.   Distances are normalized so that the minimum distance in each example is 1.  

In each example, we suppose $\delta = 2$.  Note, for any fixed $\delta$, in any finite dataset, there is a constant $C$ which will correspondingly satisfy Condition 1 provided unique observations never have distance zero from each other.\footnote{Note that our theorems, being asymptotic in nature, prove that for any fixed pair $(C, \delta)$, and tolerable distortion $tol$ to nominal false rejection probability $\alpha$ for some test $H_0$, there is $n_0(tol)$ (depending also on mixing rates, moments, etc) such that the probability of a false rejection of $H_0$ is bounded by $\alpha + tol$ for sample sizes $n > n_0(tol)$.  Note that assessing the value of $n_0(tol)$ for a particular application involves unearthing the implied bounds implicit in the central limit theorem in \cite{Jenish2009} or extensive simulation and assumptions about population moments (which we provide in Section 4) of the main paper. } In the first two examples, $\delta=2$ is a sensible choice, as the measure of 2-dimensional Euclidean balls grow as $r^2$.  In the third example, we find that this choice also fits the growth of ball cardinalities reasonably well.
Note that in \cite{Conley:Dupor:Sectoral} the quantity $\delta$ measures growth rate of cardinalities of balls as a function of $r$ and serves as the key notion of dimensionality, even though the ambient Euclidean space has dimension equal to 17, corresponding to the number of sectors.  Valid quadratic monomial upper and lower envelopes printed with a solid line.  These are given by: 
$0.002 r^2 \leq |B_r(i)| \leq 2.000 r^2$ in Condra et al (2018) locations, $ 0.500 r^2 \leq |B_r(i)| \leq 4.000 r^2$  in square grid locations,
$0.125 r^2 \leq |B_r(i)| \leq 3.000 r^2$ in Conley and Dupor (2003) locations.  In the first panels, the presence of small (cardinality 1) balls at small radii drives the constant implied from the lower bound upwards.  We remark that this lower bound is relevant for the properties of $k$-\texttt{medoids} and that those properties continue to hold if Condition 1 is satisfied only for $r > r_0$ for some sufficiently slowly increasing $r_0$.  Correspondingly, a polynomial bound (rather than monomial) of $0.012 r^2 - 5.000 \leq |B_r(i)|$ provides a tighter envelope.

\subsection{Discussion for BCH with learned clusters.}
\label{subsec: bch discussion}

Recall that the CCE test for $H_0: \theta_0 = \theta^{\circledast}$ using a partition $\mathcal C_n$ with $\G$ clusters is defined by 
$T_{\text{CCE}(a),n} = \ \text{Reject} \ \ \text{ if } \ \ \frac{ \hat \theta_{\mathsf X_n} - \theta^{\circledast}}{ \hat V_{\text{CCE}, \mathcal C_n}^{1/2}} > \sqrt{\frac{\G}{\G-1}} \times t_{1-a/2,\G-1}$
where $\hat V_{\text{CCE},\mathcal C_n}$ is the standard cluster covariance estimator (without a degrees of freedom correction).  BCH analyzed the CCE estimator under an asymptotic frame with a fixed, finite number of clusters $\G$ (i.e. $\G$ independent of $n$).  The resulting inference is based on calculating a $t$ statistic based on an estimated standard error.   Under regularity conditions, BCH show that the resulting $t$ statistic is asymptotically pivotal, but distributed according to $\sqrt{\G/(\G-1)} \times t_{\G-1}$ where $t_{\G-1}$ is the $t$ distribution with $\G-1$ degrees of freedom.

The regularity conditions required in \BCH \ are strong, and in particular, require that the clusters have asymptotically equal numbers of observations.  The $k$-\texttt{medoids} algorithm does not generally return such a partition.  Thus, the CCE estimator is not anticipated to have asymptotically correct size.

\subsection{Discussion of uniformity.}
\label{subsec: uniformity}

The high-level conditions used for proving Theorems \ref{im_theorem_main} and \ref{crs_theorem_main} involve the assumption of uniform central limit theorems for convergences $\sup_{\mathcal C \in \mathscr C_n} \| \tilde S_{\mathcal C} - \tilde S^*_\mathcal C \|_2 \rightarrow 0$ for Gaussian random variables $\tilde S_{\mathcal C}$.  Uniform convergence results may be explicitly derived in several cases of interest.  First, when each $\mathscr C_n$ contains only finitely many partitions $\mathcal C$ all of which have cardinality bounded by some ${\G_{\mathrm{max}}}$ independent of $n$, then uniform convergence follows immediately from a pointwise convergence result like that presented in Section \ref{extension_to_ols}.

Another case of interest involves increasing sequences ${\G_{\mathrm{max},n}}$.  In particular, if ${\G_{\mathrm{max},n}}$ increases sufficiently slowly, then uniform analogues of the results in Theorems \ref{im_theorem_main} and \ref{crs_theorem_main} may be anticipated by establishing Berry-Esseen-type bounds for dependent processes.   Note, for instance, that \cite{jirak2016} establishes a unidimensional Berry-Esseen bound for sums of weakly dependent random variables.
Note also that Theorem \ref{km_prop} can be extended to sequences ${\G_{\mathrm{max},n}}$ by the same argument as given in the supplemental material with $\check n : = \lfloor n/\G \rfloor$ replacing $n$ in the appropriate places. 


The conditions imposed for Theorem \ref{im_theorem_main} do not explicitly require any bounds on $|\mathscr C_n|$ or $|\mathcal C|, \mathcal C \in \mathscr C_n$.  Condition \ref{regularity_crs}(iii) used in Theorem \ref{crs_theorem_main} does require bounds on the $|\mathscr C_n|$ or $|\mathcal C|, \mathcal C \in \mathscr C_n$ which are independent of $n$, but this requirement can be dispensed if $w(S)$ is sufficiently regular.   This regularity condition is given in Condition \ref{regularity_crs}(iii)' in the supplemental material.  There also is provided a proposition showing that $w(S) = |t(S)|$ satisfies the required regularity condition.


\subsection{Computation times.}
\label{computation:times}
The procedure defined in Algorithm 1 in the main text is computationally intensive.  This section reports various computation times associated with running Algorithm 1 as well as the alternative techniques described in the main text.

We present two tables documenting computation time.  Table \ref{inferenceruntimes} records the computation time for all inferential methods in all simulation settings.   Table \ref{kmedoidsiterations} records the number of iterations needed to reach convergence in $k$-medoids using the metric space defined by the Afghanistan district centroids as described in the main text. Note that we apply $k$-medoids to 100 sets of randomly chosen centroids and pick the one with the lowest objective function value at convergence.

The experiment is implemented in MATLAB 9.6.0.1072779 (R2019a). 
Two models of CPUs are used: Intel Xeon E5-2690 v3 2.60GHz and Intel Xeon E5-2680 v4 2.40GHz.
We assigned 10 cores and 10 GB memory for each task.
The operating system is Linux 3.10.0-1160.25.1.el7.x86\_64.

\subsection{Optimization details for Table \ref{table: moran}.}\label{sec:opt:contours}

This section provides a plot illustrating level sets in the constrained optimization problem \eqref{eq: optimization problem} described in the main text for the empirical example.

\subsection{Additional Simulation Results.}\label{Additional:sim}

This section provides additional simulation results to complement those in the main text. 
The same settings as in Section \ref{simulation} are considered; though here we provide results with number of locations given by $N_{\mathrm{pan}}=820$ in addition to $N_{\mathrm{pan}}=205$. 
In the $N_{\mathrm{pan}}=820$ settings, four copies of the locations from the empirical example are created by reflecting the original locations over the $29^\circ$ latitude and $75^\circ$ longitude lines.
The data generating process follows Section \ref{simulation}.
The maximal number of groups to be considered in CCE, IM, and CRS is chosen to be ${\G_{\mathrm{max}}} = 12.$
In all cases, we consider 1000 simulation replications.

The reported results in this section also display more detailed information about the simulation studies than given in the main text in both the $N_{\mathrm{pan}}=820$ and $N_{\mathrm{pan}}=205$ cases.

\subsection*{References for Supplemental Material}\label{Supp:references}

\

\noindent Bester, C. A., Conley, T., and Hansen, C. (2011). Inference with dependent data using cluster

covariance estimators. \textit{Journal of Econometrics}, 165(2):137 -- 151.

\noindent Bolthausen, E. (1982). On the central limit theorem for stationary mixing random fields. \textit{Annals}  

\textit{of Probability}, 10(4):1047--1050.

\noindent Canay, I. A., Romano, J. P., and Shaikh, A. M. (2017). Randomization tests under an approximate

symmetry assumption. \textit{Econometrica}, 85(3):1013--1030.

\noindent Condra, L. N., Long, J. D., Shaver, A. C., and Wright, A. L. (2018). The logic of insurgent electoral

violence. \textit{American Economic Review}, 108(11):3199--3231.

\noindent Conley, T. G. and Dupor, B. (2003). A spatial analysis of sectoral complementarity. \textit{Journal of}

\textit{Political Economy}, 111(2):311--352.

\noindent Ibragimov, R. and M\"uller, U. K. (2010). t-statistic based correlation and heterogeneity robust

inference. \textit{Journal of Business \& Economic Statistics}, 28(4):453--468.

\noindent Jenish, N. and Prucha, I. R. (2009). Central limit theorems and uniform laws of large numbers for

arrays of random fields. \textit{Journal of Econometrics}, 150(1):86--98.

\pagebreak

\subsection*{Tables and Figures for Supplemental Material}\label{Supp:Tables}

\

\begin{table}[H]
	\center
	\caption{Impact of Early Morning Attacks on Voter Turnout during the 2014 Election}
	\label{tab: empirical results 2}
		\renewcommand{\arraystretch}{1.5}
	\begin{threeparttable}
		\begin{tabular}{lp{1.cm}p{.9cm}p{.9cm}p{1.3cm}p{1.3cm}p{1.3cm}p{1.3cm}p{.65cm}c}
			\hline \hline
			\multicolumn{10}{c}{Cluster-based Inference} \\
			\hline	
			&\ \ \ {\footnotesize{}$ \phantom{ \Big |}\hat{\theta}_{0}$} & \ {\footnotesize{}$s.e.$}  & {\footnotesize{}$t$-stat} &\multicolumn{2}{c}{\footnotesize{}C.I.}  &\multicolumn{2}{c}{\footnotesize{}C.I. (unadjusted cr. val.)} & {\footnotesize{}${\hat \G}$} & {\footnotesize{}${\hat \alpha}$}  \tabularnewline
			\hline
			\multicolumn{10}{c}{$\G_{\mathrm{max}}$ = 8} \\
			\hline			
			CCE  & -0.145  &0.090 &  1.609 &   ( -0.404, & \ 0.114 )  &   ( -0.358, & \ 0.068 )  & 8 & 0.024\\
			IM   &  -0.122 & 0.432  & 0.282 & ( -1.320, & \ 1.077 )   & ( -1.320, & \ 1.077 )  & 5  & 0.050 \\
			CRS  & -0.242 &  &  & ( -1.497, & \ 0.084 )   & ( -1.497, & \ 0.084 )  & 6 & 0.031\\
			\hline
			\multicolumn{10}{c}{$\G_{\mathrm{max}}$ = 16} \\
			\hline			
			CCE  & -0.145  &0.075 &  1.942 &   ( -0.344, & \ 0.054 )  &   ( -0.306, & \ 0.016 )  & 14 & 0.009\\
			IM   &  -0.122 & 0.432  & 0.282 & ( -1.320, & \ 1.077 )   & ( -1.320, & \ 1.077 )  & 5  & 0.050 \\
			CRS  & -0.242 &  &  & ( -1.497, & \ 0.084 )   & ( -1.497, & \ 0.084 )  & 6 & 0.031\\
			\hline			\\
		\end{tabular} 
		\begin{tablenotes}
			\item {Notes:} This table presents inferential results based on selected clusters. Row labels indicate which procedure
			is used. The column labeled $\hat{\theta}_{0}$ reports the IV estimate of $\theta_{0}$ for the full sample in the rows labeled UNIT and CCE and the average of IV estimators of $\theta_{0}$ of the five or six data-generated clusters in the rows labeled IM and CRS respectively. Column $s.e.$   
			reports the estimated standard errors obtained in each procedure. Note that CRS does not rely on an explicit standard error estimate.
			Column $t$-stat reports the $t$-statistic for testing the null hypothesis that $\theta_{0} = 0$ for each of the procedures. Column C.I.
			reports confidence intervals for $\theta_{0}$ obtained using $\hat\alpha$ and $\hat{k}$ from \eqref{eq: optimization problem} for each procedure. Column C.I. (unadjusted crit. val.s)
			reports confidence intervals of the IV estimate of $\theta_{0}$ using fixed asymptotic thresholds, i.e., rejecting if the p-value is less than .05.
			Note that CRS produces the same CI in both cases because of the discreteness of the test.  Column ${\hat \G}$ and $\widehat{\alpha}$ indicate the number of
			clusters and significance level selected in each procedure, separately. 
		\end{tablenotes}
	\end{threeparttable}
\end{table}

\

\pagebreak
\clearpage

\

\begin{figure}[t]
	\includegraphics[scale=.22]{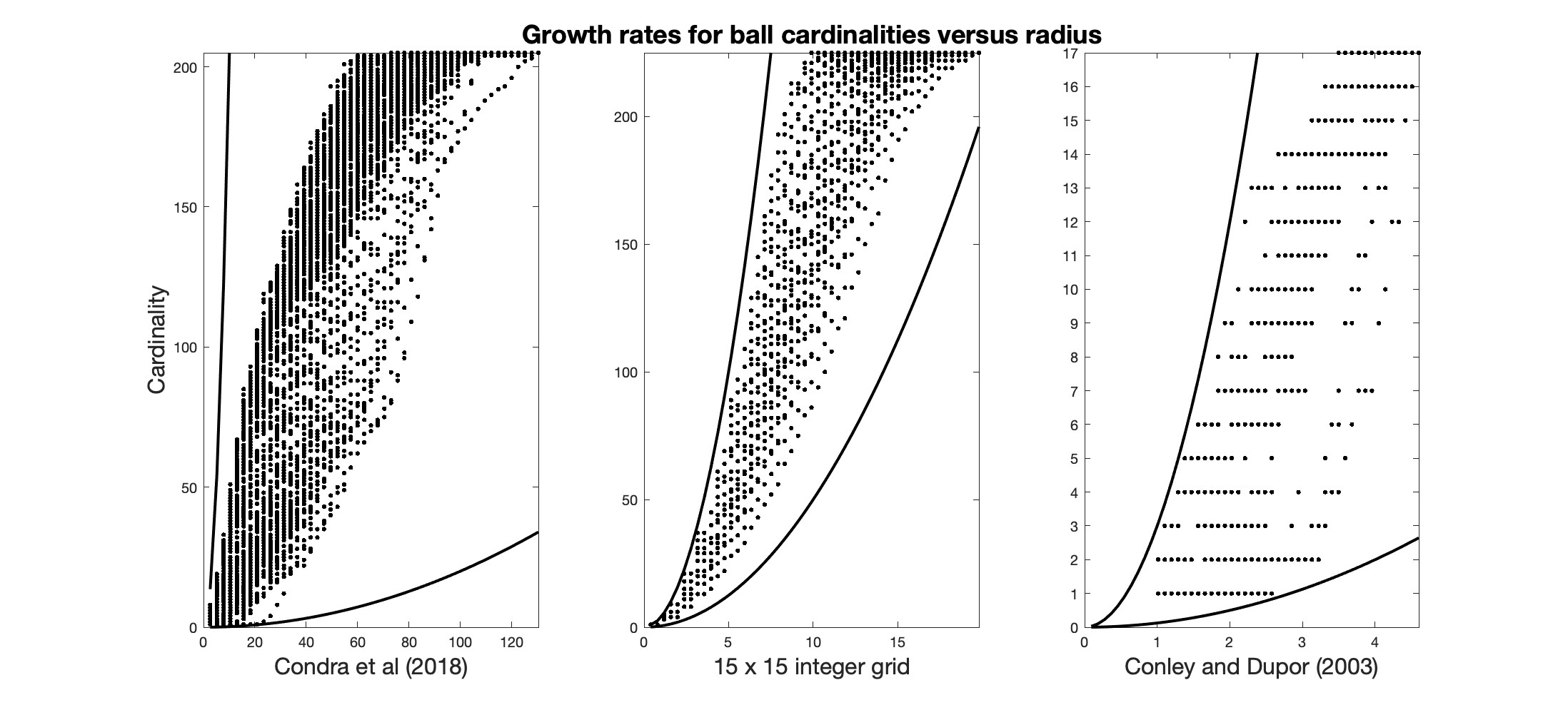}
	\caption{\linespread{1.65} \normalsize  Growth rates for ball cardinalities for three examples of metric spaces as described in the main text. Points correspond to pairs $(i,r)$ with $i \in \mathsf X$.  Each illustrated point has (horizontal, vertical) axis coordinates equal to $(r, \ | \mathsf B_r(i)  |)$.   Distances are normalized so that the minimum distance in each example is 1.  Valid quadratic monomial upper and lower envelopes printed with a solid line.  These are given by $0.002 r^2 \leq |B_r(i)| \leq 2.000 r^2$ in Condra et al (2018) locations, $ 0.500 r^2 \leq |B_r(i)| \leq 4.000 r^2$  in square grid locations, and $0.125 r^2 \leq |B_r(i)| \leq 3.000 r^2$ in Conley and Dupor (2003) locations. }
	\label{fig: ball growth}
\end{figure}


\

\pagebreak
\clearpage

\

\begin{table}[ht!]
	\centering
	\caption{Runtime in seconds}
	 \label{inferenceruntimes}
	\renewcommand{\arraystretch}{1.5}
	\begin{threeparttable}
	\begin{tabular}{cccccccccccccccc}
		\hline
		\hline
	Method	& \multicolumn{5}{c}{$N=205$} & & \multicolumn{5}{c}{$N=820$} \\
		\cline{2-6} \cline{8-12}
		& \multicolumn{2}{c}{OLS} && \multicolumn{2}{c}{IV} 		&& \multicolumn{2}{c}{OLS} && \multicolumn{2}{c}{IV}\\
		\cline{2-3} \cline{5-6} \cline{8-9} \cline{11-12}
		& BASELINE & SAR & & BASELINE & SAR &	& BASELINE & SAR & & BASELINE & SAR\\
		\hline
SK     & 52  & 38  &  & 44  & 19 &  & 19  & 38  &  & 102 & 166 \\
UNIT-U & 51  & 85  &  & 48  & 87 &  & 116 & 147 &  & 45  & 20  \\
UNIT   & 49  & 24  &  & 58  & 23 &  & 223 & 210 &  & 279 & 231 \\
CCE    & 71  & 67  &  & 62  & 36 &  & 137 & 206 &  & 427 & 293 \\
IM     & 72  & 66  &  & 61  & 60 &  & 282 & 214 &  & 308 & 358 \\
CRS    & 119 & 118 &  & 107 & 95 &  & 409 & 395 &  & 504 & 496\\
		\hline
	\end{tabular}
\begin{tablenotes}
	\item Notes: 
	This table records computation runtime (in seconds) of implementing a specific method in the corresponding setting for once. 
	\end{tablenotes}
\end{threeparttable}
\end{table}

\

\pagebreak
\clearpage

\

\begin{table}[ht!]
	\centering
	\renewcommand{\arraystretch}{1.5}
	\caption{Number of iterations in $k$-medoids}
	 \label{kmedoidsiterations}
	\begin{threeparttable}
	\begin{tabular}{cccccc}
		\hline
		\hline
$G$	& \multicolumn{2}{c}{$N=205$} & & \multicolumn{2}{c}{$N=820$}  \\
 \cline{2-3}\cline{5-6}
& best & average & & best & average \\
		\hline
2  & 3 & 2.71 &  & 3 & 2.59  \\
3  & 3 & 3.11 &  & 3 & 2.89  \\
4  & 3 & 3.31 &  & 2 & 2.58  \\
5  & 3 & 3.67 &  & 3 & 9.95  \\
6  & 4 & 3.79 &  & 3 & 22.67 \\
7  & 5 & 3.90 &  & 4 & 13.05 \\
8  & 3 & 3.93 &  & 4 & 3.50  \\
9  & -  &    -  &  & 3 & 3.50  \\
10 &  - &   -   &  & 6 & 3.70  \\
11 &  - &   -   &  & 4 & 4.14  \\
12 &  - & -     &  & 6 & 4.60 \\
\hline
	\end{tabular}
\begin{tablenotes}
\item Notes:  Iterations until convergence for $k$-medoids using Afghanistan voting districts for $N=205$ and a space derived from spatially displaced copies of Afghanistan voting districts for $N=820$ as described in the main text. 
The ``best'' column shows the number of iterations until numeric convergence for the initial centroids that achieve the lowest objective function value at convergence among the 100  sets of initial centroids. 
The ``average'' column shows the average number of iterations until numeric convergence across all 100 sets of initial centroids. 
\end{tablenotes}
\end{threeparttable}
\end{table}

\

\pagebreak
\clearpage

\

\begin{figure}[t]
	\centering
	\includegraphics[scale=.8,trim={2cm 8.5cm 0cm 9cm}]{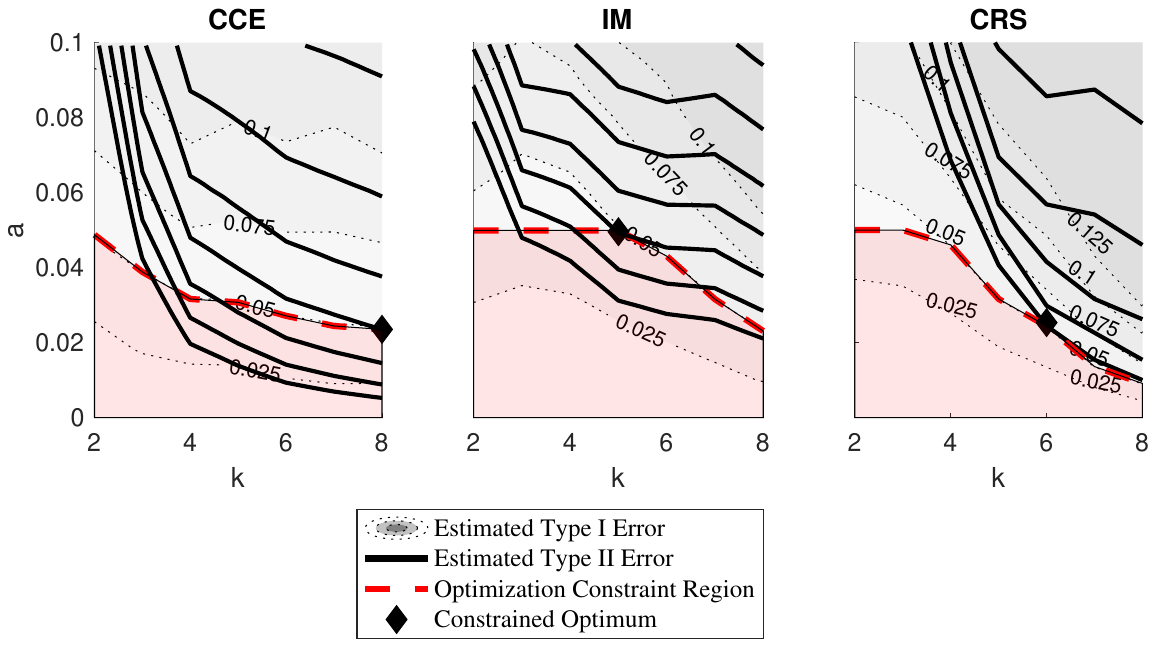} 
	\caption{\linespread{1.65} \normalsize Contours plot displaying level sets for the constrained optimization problem \eqref{eq: optimization problem} in the empirical example. The CRS panel renders using smooth interpolation because of the discreteness of CRS test.}
	\label{fig: Contours}
\end{figure}

\pagebreak
\clearpage


\begin{table}[H]
	\centering
	\renewcommand{\arraystretch}{1.2}
	\caption{Summary: OLS - BASELINE ($N_{\mathrm{pan}}=205$)}
	\label{table_summary_ols_baseline}
		\renewcommand{\arraystretch}{1.5}
	\begin{threeparttable}	
		\begin{tabular}{cccccccccccccc}
			\hline
			\hline
Method & Estim. Mean & Estim. RMSE & Size & \multicolumn{4}{c}{Power}\\
 \cline{5-8}
 & & & & -1 & -0.5 & 0.5 & 1\\
 \hline
SK     & 0.015 & 0.337 & 0.381 & 0.974 & 0.730 & 0.727 & 0.981 \\
UNIT-U & 0.015 & 0.337 & 0.577 & 0.992 & 0.844 & 0.832 & 0.987 \\
UNIT   & 0.015 & 0.337 & 0.047 & 0.831 & 0.335 & 0.319 & 0.818 \\
CCE    & 0.015 & 0.337 & 0.046 & 0.717 & 0.292 & 0.248 & 0.704 \\
IM     & 0.014 & 0.213 & 0.044 & 0.979 & 0.553 & 0.523 & 0.962 \\
CRS    & 0.014 & 0.213 & 0.042 & 0.957 & 0.514 & 0.519 & 0.970\\
 \hline 
		\end{tabular}
		\begin{tablenotes}
			\item Notes: Simulation results for estimation in the design described in Section \ref{simulation}. 
			The nominal size is 0.05. 
			Estimates are presented for the estimators, SK, UNIT-U, UNIT, CCE, IM, CRS described in the text. 
			Columns display method, estimated mean, estimated RMSE, size, and power against four alternatives (-1, -0.5, 0.5, 1). 
		\end{tablenotes}
	\end{threeparttable}
\end{table}

\

\pagebreak
\clearpage

\

\begin{table}[H]
	\centering
	\renewcommand{\arraystretch}{1.2}
	\caption{Summary: OLS- SAR ($N_{\mathrm{pan}}=205$)}
	\label{table_summary_ols_sar}
		\renewcommand{\arraystretch}{1.5}
	\begin{threeparttable}	
		\begin{tabular}{cccccccccccccc}
			\hline
			\hline
Method & Estim. Mean & Estim. RMSE & Size & \multicolumn{4}{c}{Power}\\
 \cline{5-8}
 & & & & -1 & -0.5 & 0.5 & 1\\
 \hline
SK     & -0.013 & 0.866 & 0.447 & 0.688 & 0.514 & 0.517 & 0.691 \\
UNIT-U & -0.013 & 0.866 & 0.639 & 0.761 & 0.683 & 0.681 & 0.766 \\
UNIT   & -0.013 & 0.866 & 0.359 & 0.657 & 0.461 & 0.452 & 0.656 \\
CCE    & -0.013 & 0.866 & 0.047 & 0.362 & 0.160 & 0.167 & 0.372 \\
IM     & -0.006 & 0.385 & 0.038 & 0.586 & 0.216 & 0.189 & 0.588 \\
CRS    & -0.003 & 0.354 & 0.049 & 0.674 & 0.231 & 0.256 & 0.670\\
 \hline 
		\end{tabular}
		\begin{tablenotes}
			\item Notes: Simulation results for estimation in the design described in Section \ref{simulation}. 
			The nominal size is 0.05. 
			Estimates are presented for the estimators, SK, UNIT-U, UNIT, CCE, IM, CRS described in the text. 
			Columns display method, estimated mean, estimated RMSE, size, and power against four alternatives (-1, -0.5, 0.5, 1).
		\end{tablenotes}
	\end{threeparttable}
\end{table}

\

\pagebreak
\clearpage

\

\begin{table}[H]
	\centering
	\renewcommand{\arraystretch}{1.2}
	\caption{Summary: IV - BASELINE ($N_{\mathrm{pan}}=205$)}
\label{table_summary_iv_baseline}
	\renewcommand{\arraystretch}{1.5}
	\begin{threeparttable}	
		\begin{tabular}{cccccccccccccc}
			\hline
			\hline
Method & Estim. Median & Estim. MAD & Size & \multicolumn{4}{c}{Power}\\
 \cline{5-8}
 & & & & -1 & -0.5 & 0.5 & 1\\
 \hline
SK     & 0.002  & 0.114 & 0.362 & 1.000 & 0.999 & 0.913 & 0.993 \\
UNIT-U & 0.002  & 0.114 & 0.568 & 0.996 & 0.950 & 1.000 & 1.000 \\
UNIT   & 0.002  & 0.114 & 0.074 & 0.954 & 0.711 & 0.910 & 1.000 \\
CCE    & 0.002  & 0.114 & 0.062 & 0.927 & 0.652 & 0.729 & 0.987 \\
IM     & -0.059 & 0.091 & 0.046 & 0.930 & 0.736 & 0.951 & 0.985 \\
CRS    & -0.061 & 0.095 & 0.040 & 0.992 & 0.992 & 0.680 & 0.921\\
 \hline 
		\end{tabular}
		\begin{tablenotes}
			\item Notes: Simulation results for estimation in the design described in Section \ref{simulation}. 
			The nominal size is 0.05. 
			Estimates are presented for the estimators, SK, UNIT-U, UNIT, CCE, IM, CRS described in the text. 
			Columns display method, estimated median, estimated MAD, size, and power against four alternatives (-1, -0.5, 0.5, 1). 
		\end{tablenotes}
	\end{threeparttable}
\end{table}

\

\pagebreak
\clearpage

\

\begin{table}[H]
	\centering
	\renewcommand{\arraystretch}{1.2}
	\caption{Summary: IV - SAR ($N_{\mathrm{pan}}=205$)}
	\label{table_summary_iv_sar}
		\renewcommand{\arraystretch}{1.5}
	\begin{threeparttable}	
		\begin{tabular}{cccccccccccccc}
			\hline
			\hline
Method & Estim. Median & Estim. MAD & Size & \multicolumn{4}{c}{Power}\\
 \cline{5-8}
 & & & & -1 & -0.5 & 0.5 & 1\\
 \hline
SK     & 0.002  & 0.280 & 0.324 & 0.898 & 0.655 & 0.627 & 0.768 \\
UNIT-U & 0.002  & 0.280 & 0.502 & 0.797 & 0.675 & 0.798 & 0.942 \\
UNIT   & 0.002  & 0.280 & 0.256 & 0.737 & 0.604 & 0.567 & 0.839 \\
CCE    & 0.002  & 0.280 & 0.083 & 0.634 & 0.464 & 0.282 & 0.626 \\
IM     & -0.046 & 0.158 & 0.025 & 0.658 & 0.420 & 0.470 & 0.773 \\
CRS    & -0.095 & 0.213 & 0.041 & 0.778 & 0.703 & 0.331 & 0.563\\
 \hline 
		\end{tabular}
		\begin{tablenotes}
			\item Notes: Simulation results for estimation in the design described in Section \ref{simulation}. 
			The nominal size is 0.05. 
			Estimates are presented for the estimators, SK, UNIT-U, UNIT, CCE, IM, CRS described in the text. 
			Columns display method, estimated median, estimated MAD, size, and power against four alternatives (-1, -0.5, 0.5, 1). 
		\end{tablenotes}
	\end{threeparttable}
\end{table}

\

\pagebreak
\clearpage

\

\begin{table}[H]
	\caption{Clustering: OLS - BASELINE ($N_{\mathrm{pan}}=205$)}
	\label{table_clustering_ols_baseline}
	\centering
	\renewcommand{\arraystretch}{1.2}
	\begin{threeparttable}	\begin{tabular}{cccccccccccccc}
			\hline
			\hline
			& & \multicolumn{7}{c}{$k$} & $\hat k$\\
			\cline{3-9}
			& & 2 & 3 & 4 & 5 & 6 & 7 & 8 \\
			\hline
CCE & size (usual cv)               & 0.085 & 0.126 & 0.154 & 0.176 & 0.192 & 0.201 & 0.222 &       \\
& size (simulated cv)           & 0.056 & 0.049 & 0.048 & 0.046 & 0.050 & 0.042 & 0.044 & 0.046 \\
& $\widehat{k}$ frequency                & 0.000 & 0.000 & 0.001 & 0.098 & 0.311 & 0.260 & 0.330 &       \\
& p(sim\_size\textgreater{}.05) & 1.000 & 1.000 & 1.000 & 1.000 & 1.000 & 1.000 & 1.000 & 1.000 \\
\cline{2-10}
& $\widehat{\alpha}$ quantile &       &       &       &       &       &       &       \\
& q10                           & 0.025 & 0.016 & 0.009 & 0.008 & 0.007 & 0.005 & 0.003 & 0.004 \\
& q25                           & 0.027 & 0.017 & 0.010 & 0.009 & 0.007 & 0.006 & 0.004 & 0.006 \\
& q50                           & 0.030 & 0.019 & 0.011 & 0.011 & 0.009 & 0.007 & 0.005 & 0.008 \\
& q75                           & 0.033 & 0.021 & 0.012 & 0.012 & 0.010 & 0.008 & 0.006 & 0.010 \\
& q90                           & 0.035 & 0.023 & 0.014 & 0.014 & 0.011 & 0.009 & 0.007 & 0.012 \\
\hline
IM  & size (usual cv)               & 0.041 & 0.058 & 0.063 & 0.062 & 0.054 & 0.050 & 0.054 &       \\
& size (simulated cv)           & 0.040 & 0.053 & 0.059 & 0.054 & 0.047 & 0.045 & 0.042 & 0.044 \\
& $\widehat{k}$ frequency                & 0.000 & 0.000 & 0.000 & 0.007 & 0.018 & 0.113 & 0.862 &       \\
& p(sim\_size\textgreater{}.05) & 0.383 & 0.752 & 0.832 & 0.909 & 0.888 & 0.627 & 0.791 & 0.762 \\
\cline{2-10}
& $\widehat{\alpha}$ quantile       &       &       &       &       &       &       \\
& q10                           & 0.044 & 0.038 & 0.037 & 0.035 & 0.035 & 0.040 & 0.037 & 0.039 \\
& q25                           & 0.048 & 0.041 & 0.040 & 0.038 & 0.039 & 0.044 & 0.040 & 0.042 \\
& q50                           & 0.050 & 0.046 & 0.044 & 0.042 & 0.043 & 0.048 & 0.044 & 0.046 \\
& q75                           & 0.050 & 0.050 & 0.048 & 0.046 & 0.047 & 0.050 & 0.049 & 0.050 \\
& q90                           & 0.050 & 0.050 & 0.050 & 0.050 & 0.050 & 0.050 & 0.050 & 0.050 \\
\hline
CRS & size (usual cv)               & 0.000 & 0.000 & 0.000 & 0.000 & 0.031 & 0.040 & 0.048 &       \\
& size (simulated cv)           & 0.000 & 0.000 & 0.000 & 0.000 & 0.030 & 0.033 & 0.041 & 0.042 \\
& $\widehat{k}$ frequency                & 0.000 & 0.000 & 0.000 & 0.000 & 0.001 & 0.099 & 0.900 &       \\
& p(sim\_size\textgreater{}.05) & 0.000 & 0.000 & 0.000 & 0.000 & 0.054 & 0.530 & 0.726 & 0.641 \\
\cline{2-10}
& $\widehat{\alpha}$ quantile        &       &       &       &       &       &       \\
& q10                           & 0.050 & 0.050 & 0.050 & 0.050 & 0.050 & 0.047 & 0.039 & 0.039 \\
& q25                           & 0.050 & 0.050 & 0.050 & 0.050 & 0.050 & 0.047 & 0.039 & 0.043 \\
& q50                           & 0.050 & 0.050 & 0.050 & 0.050 & 0.050 & 0.047 & 0.047 & 0.047 \\
& q75                           & 0.050 & 0.050 & 0.050 & 0.050 & 0.050 & 0.050 & 0.050 & 0.050 \\
& q90                           & 0.050 & 0.050 & 0.050 & 0.050 & 0.050 & 0.050 & 0.050 & 0.050 \\
			\hline
		\end{tabular}
		\begin{tablenotes}
			\item Notes: Simulation results for the design described in Section \ref{simulation}. 
			Inferential properties are presented for the estimators, CCE, IM, CRS, described in the text. 
			Columns under $k$ report results with the number of groups fixed at a certain $k$.
			$\hat k$ is the number of clusters chosen by the criterion on size-power tradeoff described in the text. 
						Sizes under both usual critical value (``size (usual cv)'') and adjusted critical value (``size (simulated cv)'') are reported. 
			The rows ``$\hat k$ frequency'' is the frequency of a particular $k$ achieving the highest simulated power among candidate $k$'s in the setting. 
			The rows ``p(sim\_size\textgreater.05)'' report the empirical frequency of $\widehat {\mathrm{Err}}_{\mathrm{Type-I}}(\bullet(.05), k)>.05$.
			The row ``$\widehat{\alpha}$ quantile'' report the quantiles of the selected number of groups. 
		\end{tablenotes}
	\end{threeparttable}
\end{table}

\

\pagebreak
\clearpage

\

\begin{table}[H]
	\caption{Clustering: OLS - SAR ($N_{\mathrm{pan}}=205$)}
	\label{table_clustering_ols_sar}
	\centering
	\renewcommand{\arraystretch}{1.2}
	\begin{threeparttable}	\begin{tabular}{cccccccccccccc}
			\hline
			\hline
			& & \multicolumn{7}{c}{$k$} & $\hat k$\\
			\cline{3-9}
			& & 2 & 3 & 4 & 5 & 6 & 7 & 8 \\
			\hline 
		CCE & size (usual cv)               & 0.041 & 0.053 & 0.090 & 0.067 & 0.080 & 0.082 & 0.095 &       \\
		& size (simulated cv)           & 0.034 & 0.035 & 0.037 & 0.037 & 0.049 & 0.040 & 0.043 & 0.047 \\
		& $\widehat{k}$ frequency                & 0.000 & 0.001 & 0.022 & 0.051 & 0.113 & 0.268 & 0.545 &       \\
		& p(sim\_size\textgreater{}.05) & 0.940 & 1.000 & 1.000 & 1.000 & 1.000 & 1.000 & 1.000 & 1.000 \\
		\cline{2-10}
		& $\widehat{\alpha}$ quantile &       &       &       &       &       &       &       \\
		& q10                           & 0.034 & 0.026 & 0.020 & 0.024 & 0.021 & 0.021 & 0.018 & 0.021 \\
		& q25                           & 0.037 & 0.029 & 0.022 & 0.027 & 0.023 & 0.024 & 0.021 & 0.024 \\
		& q50                           & 0.041 & 0.032 & 0.025 & 0.030 & 0.026 & 0.027 & 0.024 & 0.027 \\
		& q75                           & 0.045 & 0.035 & 0.028 & 0.033 & 0.030 & 0.031 & 0.027 & 0.031 \\
		& q90                           & 0.048 & 0.038 & 0.031 & 0.037 & 0.033 & 0.034 & 0.031 & 0.034 \\
		\hline
		IM  & size (usual cv)               & 0.022 & 0.024 & 0.021 & 0.027 & 0.038 & 0.049 & 0.051 &       \\
		& size (simulated cv)           & 0.021 & 0.024 & 0.020 & 0.026 & 0.038 & 0.044 & 0.047 & 0.038 \\
		& $\widehat{k}$ frequency                & 0.000 & 0.003 & 0.069 & 0.554 & 0.040 & 0.176 & 0.158 &       \\
		& p(sim\_size\textgreater{}.05) & 0.389 & 0.419 & 0.385 & 0.426 & 0.385 & 0.536 & 0.521 & 0.336 \\
				\cline{2-10}
		& $\widehat{\alpha}$ quantile &       &       &       &       &       &       &       \\
		& q10                           & 0.043 & 0.043 & 0.044 & 0.043 & 0.044 & 0.041 & 0.042 & 0.047 \\
		& q25                           & 0.047 & 0.047 & 0.047 & 0.047 & 0.048 & 0.045 & 0.045 & 0.049 \\
		& q50                           & 0.050 & 0.050 & 0.050 & 0.050 & 0.050 & 0.049 & 0.050 & 0.050 \\
		& q75                           & 0.050 & 0.050 & 0.050 & 0.050 & 0.050 & 0.050 & 0.050 & 0.050 \\
		& q90                           & 0.050 & 0.050 & 0.050 & 0.050 & 0.050 & 0.050 & 0.050 & 0.050 \\
		\hline
		CRS & size (usual cv)               & 0.000 & 0.000 & 0.000 & 0.000 & 0.027 & 0.050 & 0.055 &       \\
		& size (simulated cv)           & 0.000 & 0.000 & 0.000 & 0.000 & 0.027 & 0.041 & 0.046 & 0.049 \\
		& $\widehat{k}$ frequency                & 0.000 & 0.000 & 0.000 & 0.000 & 0.001 & 0.401 & 0.598 &       \\
		& p(sim\_size\textgreater{}.05) & 0.000 & 0.000 & 0.000 & 0.000 & 0.000 & 0.434 & 0.408 & 0.268 \\
				\cline{2-10}
		& $\widehat{\alpha}$ quantile &       &       &       &       &       &       &       \\
		& q10                           & 0.050 & 0.050 & 0.050 & 0.050 & 0.050 & 0.047 & 0.047 & 0.047 \\
		& q25                           & 0.050 & 0.050 & 0.050 & 0.050 & 0.050 & 0.047 & 0.047 & 0.047 \\
		& q50                           & 0.050 & 0.050 & 0.050 & 0.050 & 0.050 & 0.050 & 0.050 & 0.050 \\
		& q75                           & 0.050 & 0.050 & 0.050 & 0.050 & 0.050 & 0.050 & 0.050 & 0.050 \\
		& q90                           & 0.050 & 0.050 & 0.050 & 0.050 & 0.050 & 0.050 & 0.050 & 0.050  	 \\
			\hline
		\end{tabular}
		\begin{tablenotes}
			\item Notes: Simulation results for the design described in Section \ref{simulation}. 
Inferential properties are presented for the estimators, CCE, IM, CRS, described in the text. 
Columns under $k$ report results with the number of groups fixed at a certain $k$.
$\hat k$ is the number of clusters chosen by the criterion on size-power tradeoff described in the text. 
Sizes under both usual critical value (``size (usual cv)'') and adjusted critical value (``size (simulated cv)'') are reported. 
The rows ``$\hat k$ frequency'' is the frequency of a particular $k$ achieving the highest simulated power among candidate $k$'s in the setting. 
The rows ``p(sim\_size\textgreater.05)'' report the empirical frequency of $\widehat {\mathrm{Err}}_{\mathrm{Type-I}}(\bullet(.05), k)>.05$.
The row ``$\widehat{\alpha}$ quantile'' report the quantiles of the selected number of groups. 
		\end{tablenotes}
	\end{threeparttable}
\end{table}

\

\pagebreak
\clearpage

\

\begin{table}[H]
	\caption{Clustering: IV - BASELINE ($N_{\mathrm{pan}}=205$)}
	\label{table_clustering_iv_baseline}
	\centering
	\renewcommand{\arraystretch}{1.2}
	\begin{threeparttable}	\begin{tabular}{cccccccccccccc}
			\hline
			\hline
			& & \multicolumn{7}{c}{$k$} & $\hat k$\\
			\cline{3-9}
			& & 2 & 3 & 4 & 5 & 6 & 7 & 8 \\
			\hline
CCE & size (usual cv)               & 0.079 & 0.123 & 0.157 & 0.162 & 0.176 & 0.186 & 0.213 &       \\
& size (simulated cv)           & 0.041 & 0.053 & 0.053 & 0.052 & 0.058 & 0.059 & 0.055 & 0.062 \\
& $\widehat{k}$ frequency                & 0.000 & 0.001 & 0.001 & 0.039 & 0.242 & 0.261 & 0.456 &       \\
& p(sim\_size\textgreater{}.05) & 1.000 & 1.000 & 1.000 & 1.000 & 1.000 & 1.000 & 1.000 & 1.000 \\
		\cline{2-10}
& $\widehat{\alpha}$ quantile &       &       &       &       &       &       &       \\
& q10                           & 0.025 & 0.016 & 0.009 & 0.008 & 0.006 & 0.004 & 0.003 & 0.004 \\
& q25                           & 0.028 & 0.017 & 0.010 & 0.009 & 0.007 & 0.005 & 0.004 & 0.005 \\
& q50                           & 0.030 & 0.019 & 0.011 & 0.010 & 0.008 & 0.006 & 0.004 & 0.007 \\
& q75                           & 0.033 & 0.021 & 0.012 & 0.012 & 0.010 & 0.008 & 0.005 & 0.008 \\
& q90                           & 0.036 & 0.023 & 0.014 & 0.014 & 0.011 & 0.010 & 0.007 & 0.010 \\
\hline
IM  & size (usual cv)               & 0.045 & 0.060 & 0.062 & 0.062 & 0.045 & 0.035 & 0.044 &       \\
& size (simulated cv)           & 0.045 & 0.055 & 0.060 & 0.054 & 0.042 & 0.033 & 0.042 & 0.046 \\
& $\widehat{k}$ frequency                & 0.000 & 0.000 & 0.005 & 0.051 & 0.062 & 0.066 & 0.816 &       \\
& p(sim\_size\textgreater{}.05) & 0.300 & 0.683 & 0.702 & 0.786 & 0.640 & 0.280 & 0.414 & 0.375 \\
		\cline{2-10}
& $\widehat{\alpha}$ quantile &       &       &       &       &       &       &       \\
& q10                           & 0.045 & 0.039 & 0.038 & 0.037 & 0.040 & 0.045 & 0.043 & 0.045 \\
& q25                           & 0.049 & 0.042 & 0.042 & 0.040 & 0.043 & 0.049 & 0.047 & 0.048 \\
& q50                           & 0.050 & 0.047 & 0.046 & 0.045 & 0.048 & 0.050 & 0.050 & 0.050 \\
& q75                           & 0.050 & 0.050 & 0.050 & 0.049 & 0.050 & 0.050 & 0.050 & 0.050 \\
& q90                           & 0.050 & 0.050 & 0.050 & 0.050 & 0.050 & 0.050 & 0.050 & 0.050 \\
\hline
CRS & size (usual cv)               & 0.000 & 0.000 & 0.000 & 0.000 & 0.032 & 0.041 & 0.047 &       \\
& size (simulated cv)           & 0.000 & 0.000 & 0.000 & 0.000 & 0.029 & 0.032 & 0.036 & 0.040 \\
& $\widehat{k}$ frequency                & 0.000 & 0.000 & 0.000 & 0.000 & 0.000 & 0.157 & 0.843 &       \\
& p(sim\_size\textgreater{}.05) & 0.000 & 0.000 & 0.000 & 0.000 & 0.038 & 0.618 & 0.820 & 0.682 \\
		\cline{2-10}
& $\widehat{\alpha}$ quantile &       &       &       &       &       &       &       \\
& q10                           & 0.050 & 0.050 & 0.050 & 0.050 & 0.050 & 0.047 & 0.039 & 0.039 \\
& q25                           & 0.050 & 0.050 & 0.050 & 0.050 & 0.050 & 0.047 & 0.039 & 0.039 \\
& q50                           & 0.050 & 0.050 & 0.050 & 0.050 & 0.050 & 0.047 & 0.047 & 0.047 \\
& q75                           & 0.050 & 0.050 & 0.050 & 0.050 & 0.050 & 0.050 & 0.047 & 0.050 \\
& q90                           & 0.050 & 0.050 & 0.050 & 0.050 & 0.050 & 0.050 & 0.050 & 0.050     \\
			\hline
		\end{tabular}
		\begin{tablenotes}
			\item Notes: Simulation results for the design described in Section \ref{simulation}. 
Inferential properties are presented for the estimators, CCE, IM, CRS, described in the text. 
Columns under $k$ report results with the number of groups fixed at a certain $k$.
$\hat k$ is the number of clusters chosen by the criterion on size-power tradeoff described in the text. 
Sizes under both usual critical value (``size (usual cv)'') and adjusted critical value (``size (simulated cv)'') are reported. 
The rows ``$\hat k$ frequency'' is the frequency of a particular $k$ achieving the highest simulated power among candidate $k$'s in the setting. 
The rows ``p(sim\_size\textgreater.05)'' report the empirical frequency of $\widehat {\mathrm{Err}}_{\mathrm{Type-I}}(\bullet(.05), k)>.05$.
The row ``$\widehat{\alpha}$ quantile'' report the quantiles of the selected number of groups. 
		\end{tablenotes}
	\end{threeparttable}
\end{table}

\

\pagebreak
\clearpage

\

\begin{table}[H]
	\caption{Clustering: IV - SAR ($N_{\mathrm{pan}}=205$)}
	\label{table_clustering_iv_sar}
	\centering
	\renewcommand{\arraystretch}{1.2}
	\begin{threeparttable}	\begin{tabular}{cccccccccccccc}
			\hline
			\hline
			& & \multicolumn{7}{c}{$k$} & $\hat k$\\
			\cline{3-9}
			& & 2 & 3 & 4 & 5 & 6 & 7 & 8 \\
			\hline
CCE & size (usual cv)               & 0.039 & 0.091 & 0.148 & 0.095 & 0.113 & 0.127 & 0.138 &       \\
& size (simulated cv)           & 0.030 & 0.044 & 0.062 & 0.058 & 0.069 & 0.079 & 0.082 & 0.083 \\
& $\widehat{k}$ frequency                & 0.006 & 0.027 & 0.048 & 0.047 & 0.071 & 0.257 & 0.544 &       \\
& p(sim\_size\textgreater{}.05) & 0.929 & 0.998 & 1.000 & 0.998 & 0.998 & 0.998 & 1.000 & 0.999 \\
		\cline{2-10}
& $\widehat{\alpha}$ quantile &       &       &       &       &       &       &       \\
& q10                           & 0.033 & 0.023 & 0.014 & 0.014 & 0.009 & 0.008 & 0.006 & 0.014 \\
& q25                           & 0.037 & 0.027 & 0.020 & 0.024 & 0.019 & 0.020 & 0.016 & 0.020 \\
& q50                           & 0.041 & 0.031 & 0.024 & 0.028 & 0.025 & 0.026 & 0.022 & 0.025 \\
& q75                           & 0.045 & 0.034 & 0.027 & 0.033 & 0.029 & 0.030 & 0.027 & 0.029 \\
& q90                           & 0.049 & 0.038 & 0.031 & 0.037 & 0.033 & 0.034 & 0.031 & 0.033 \\
\hline
IM  & size (usual cv)               & 0.025 & 0.022 & 0.027 & 0.022 & 0.031 & 0.031 & 0.034 &       \\
& size (simulated cv)           & 0.025 & 0.022 & 0.027 & 0.022 & 0.031 & 0.027 & 0.031 & 0.025 \\
& $\widehat{k}$ frequency                & 0.005 & 0.078 & 0.240 & 0.568 & 0.030 & 0.038 & 0.041 &       \\
& p(sim\_size\textgreater{}.05) & 0.364 & 0.391 & 0.257 & 0.194 & 0.159 & 0.175 & 0.179 & 0.249 \\
		\cline{2-10}
& $\widehat{\alpha}$ quantile &       &       &       &       &       &       &       \\
& q10                           & 0.043 & 0.041 & 0.042 & 0.043 & 0.031 & 0.020 & 0.013 & 0.032 \\
& q25                           & 0.047 & 0.046 & 0.050 & 0.050 & 0.050 & 0.050 & 0.050 & 0.050 \\
& q50                           & 0.050 & 0.050 & 0.050 & 0.050 & 0.050 & 0.050 & 0.050 & 0.050 \\
& q75                           & 0.050 & 0.050 & 0.050 & 0.050 & 0.050 & 0.050 & 0.050 & 0.050 \\
& q90                           & 0.050 & 0.050 & 0.050 & 0.050 & 0.050 & 0.050 & 0.050 & 0.050 \\
\hline
CRS & size (usual cv)               & 0.000 & 0.000 & 0.000 & 0.000 & 0.027 & 0.045 & 0.049 &       \\
& size (simulated cv)           & 0.000 & 0.000 & 0.000 & 0.000 & 0.017 & 0.034 & 0.037 & 0.041 \\
& $\widehat{k}$ frequency                & 0.000 & 0.000 & 0.000 & 0.000 & 0.032 & 0.480 & 0.488 &       \\
& p(sim\_size\textgreater{}.05) & 0.000 & 0.000 & 0.000 & 0.000 & 0.209 & 0.620 & 0.641 & 0.497 \\
		\cline{2-10}
& $\widehat{\alpha}$ quantile &       &       &       &       &       &       &       \\
& q10                           & 0.050 & 0.050 & 0.050 & 0.050 & 0.031 & 0.016 & 0.008 & 0.008 \\
& q25                           & 0.050 & 0.050 & 0.050 & 0.050 & 0.050 & 0.047 & 0.031 & 0.043 \\
& q50                           & 0.050 & 0.050 & 0.050 & 0.050 & 0.050 & 0.047 & 0.047 & 0.050 \\
& q75                           & 0.050 & 0.050 & 0.050 & 0.050 & 0.050 & 0.050 & 0.050 & 0.050 \\
& q90                           & 0.050 & 0.050 & 0.050 & 0.050 & 0.050 & 0.050 & 0.050 & 0.050\\
			\hline
		\end{tabular}
		\begin{tablenotes}
			\item Notes: Simulation results for the design described in Section \ref{simulation}. 
Inferential properties are presented for the estimators, CCE, IM, CRS, described in the text. 
Columns under $k$ report results with the number of groups fixed at a certain $k$.
$\hat k$ is the number of clusters chosen by the criterion on size-power tradeoff described in the text. 
Sizes under both usual critical value (``size (usual cv)'') and adjusted critical value (``size (simulated cv)'') are reported. 
The rows ``$\hat k$ frequency'' is the frequency of a particular $k$ achieving the highest simulated power among candidate $k$'s in the setting. 
The rows ``p(sim\_size\textgreater.05)'' report the empirical frequency of $\widehat {\mathrm{Err}}_{\mathrm{Type-I}}(\bullet(.05), k)>.05$.
The row ``$\widehat{\alpha}$ quantile'' report the quantiles of the selected number of groups. 
		\end{tablenotes}
	\end{threeparttable}
\end{table}

\

\pagebreak
\clearpage

\

\begin{table}[H]
	\centering
	\renewcommand{\arraystretch}{1.2}
	\caption{Summary: OLS - BASELINE ($N_{\mathrm{pan}}=820$)}
	\label{table_summary_ols_baseline_large}
		\renewcommand{\arraystretch}{1.5}
	\begin{threeparttable}	
		\begin{tabular}{cccccccccccccc}
			\hline
			\hline
Method & Estim. Mean & Estim. RMSE & Size & \multicolumn{4}{c}{Power}\\
 \cline{5-8}
 & & & & -1 & -0.5 & 0.5 & 1\\
 \hline
 SK     & 0.002  & 0.259 & 0.395 & 0.998 & 0.856 & 0.847 & 0.998 \\
 UNIT-U & 0.002  & 0.259 & 0.700 & 1.000 & 0.938 & 0.944 & 1.000 \\
 UNIT   & 0.002  & 0.259 & 0.039 & 0.962 & 0.517 & 0.496 & 0.960 \\
 CCE    & 0.002  & 0.259 & 0.052 & 0.917 & 0.449 & 0.449 & 0.899 \\
 IM     & -0.001 & 0.146 & 0.058 & 1.000 & 0.877 & 0.881 & 1.000 \\
 CRS    & -0.001 & 0.146 & 0.058 & 1.000 & 0.871 & 0.868 & 1.000\\
 \hline 
		\end{tabular}
		\begin{tablenotes}
			\item Notes: Simulation results for estimation in the design described in Section \ref{simulation}. 
			The nominal size is 0.05. 
			Estimates are presented for the estimators, SK, UNIT-U, UNIT, CCE, IM, CRS described in the text. 
			Columns display method, estimated mean, estimated RMSE, size, and power against four alternatives (-1, -0.5, 0.5, 1). 
		\end{tablenotes}
	\end{threeparttable}
\end{table}

\

\pagebreak
\clearpage

\

\begin{table}[H]
	\centering
	\renewcommand{\arraystretch}{1.2}
	\caption{Summary: OLS- SAR ($N_{\mathrm{pan}}=820$)}
	\label{table_summary_ols_sar_large}
	\begin{threeparttable}	
		\begin{tabular}{cccccccccccccc}
			\hline
			\hline
Method & Estim. Mean & Estim. RMSE & Size & \multicolumn{4}{c}{Power}\\
 \cline{5-8}
 & & & & -1 & -0.5 & 0.5 & 1\\
 \hline
 SK     & 0.003  & 0.229 & 0.154 & 0.996 & 0.780 & 0.787 & 0.990 \\
 UNIT-U & 0.003  & 0.229 & 0.405 & 0.999 & 0.909 & 0.903 & 1.000 \\
 UNIT   & 0.003  & 0.229 & 0.123 & 0.992 & 0.756 & 0.743 & 0.995 \\
 CCE    & 0.003  & 0.229 & 0.020 & 0.860 & 0.363 & 0.364 & 0.848 \\
 IM     & -0.013 & 0.249 & 0.035 & 0.908 & 0.465 & 0.524 & 0.911 \\
 CRS    & -0.013 & 0.229 & 0.062 & 0.959 & 0.612 & 0.553 & 0.938\\
 \hline 
		\end{tabular}
		\begin{tablenotes}
			\item Notes: Simulation results for estimation in the design described in Section \ref{simulation}. 
			The nominal size is 0.05. 
			Estimates are presented for the estimators, SK, UNIT-U, UNIT, CCE, IM, CRS described in the text. 
			Columns display method, estimated mean, estimated RMSE, size, and power against four alternatives (-1, -0.5, 0.5, 1).
		\end{tablenotes}
	\end{threeparttable}
\end{table}

\

\pagebreak
\clearpage

\

\begin{table}[H]
	\centering
	\renewcommand{\arraystretch}{1.2}
	\caption{Summary: IV - BASELINE ($N_{\mathrm{pan}}=820$)}
\label{table_summary_iv_baseline_large}
	\renewcommand{\arraystretch}{1.5}
	\begin{threeparttable}	
		\begin{tabular}{cccccccccccccc}
			\hline
			\hline
			Method & Estim. Median & Estim. MAD & Size & \multicolumn{4}{c}{Power}\\
			\cline{5-8}
			& & & & -1 & -0.5 & 0.5 & 1\\
			\hline
			SK     & 0.004  & 0.089 & 0.390 & 1.000 & 1.000 & 0.973 & 1.000 \\
			UNIT-U & 0.004  & 0.089 & 0.700 & 1.000 & 0.995 & 1.000 & 1.000 \\
			UNIT   & 0.004  & 0.089 & 0.048 & 0.998 & 0.870 & 1.000 & 1.000 \\
			CCE    & 0.004  & 0.089 & 0.051 & 0.989 & 0.837 & 0.963 & 1.000 \\
			IM     & -0.041 & 0.061 & 0.055 & 0.994 & 0.940 & 0.998 & 0.999 \\
			CRS    & -0.042 & 0.061 & 0.056 & 0.999 & 0.999 & 0.930 & 0.993\\
 \hline 
		\end{tabular}
		\begin{tablenotes}
			\item Notes: Simulation results for estimation in the design described in Section \ref{simulation}. 
			The nominal size is 0.05. 
			Estimates are presented for the estimators, SK, UNIT-U, UNIT, CCE, IM, CRS described in the text. 
			Columns display method, estimated median, estimated MAD, size, and power against four alternatives (-1, -0.5, 0.5, 1). 
		\end{tablenotes}
	\end{threeparttable}
\end{table}

\

\pagebreak
\clearpage

\

\begin{table}[H]
	\centering
	\renewcommand{\arraystretch}{1.2}
	\caption{Summary: IV - SAR ($N_{\mathrm{pan}}=820$)}
	\label{table_summary_iv_sar_large}
		\renewcommand{\arraystretch}{1.5}
	\begin{threeparttable}	
		\begin{tabular}{cccccccccccccc}
			\hline
			\hline
Method & Estim. Median & Estim. MAD & Size & \multicolumn{4}{c}{Power}\\
 \cline{5-8}
 & & & & -1 & -0.5 & 0.5 & 1\\
 \hline
SK     & 0.003  & 0.076 & 0.144 & 1.000 & 1.000 & 0.953 & 0.997 \\
UNIT-U & 0.003  & 0.076 & 0.389 & 1.000 & 0.983 & 1.000 & 1.000 \\
UNIT   & 0.003  & 0.076 & 0.118 & 0.998 & 0.948 & 1.000 & 1.000 \\
CCE    & 0.003  & 0.076 & 0.030 & 0.979 & 0.804 & 0.898 & 0.999 \\
IM     & -0.052 & 0.130 & 0.022 & 0.774 & 0.561 & 0.756 & 0.864 \\
CRS    & -0.052 & 0.123 & 0.037 & 0.921 & 0.891 & 0.652 & 0.847\\
 \hline 
		\end{tabular}
		\begin{tablenotes}
			\item Notes: Simulation results for estimation in the design described in Section \ref{simulation}. 
			The nominal size is 0.05. 
			Estimates are presented for the estimators, SK, UNIT-U, UNIT, CCE, IM, CRS described in the text. 
			Columns display method, estimated median, estimated MAD, size, and power against four alternatives (-1, -0.5, 0.5, 1). 
		\end{tablenotes}
	\end{threeparttable}
\end{table}

\

\pagebreak
\clearpage

\

\begin{figure}[t]
	\centering
	\subfloat[BASELINE]{\includegraphics[width=.45\linewidth,trim=4.5cm 6.5cm 4.5cm 7cm]{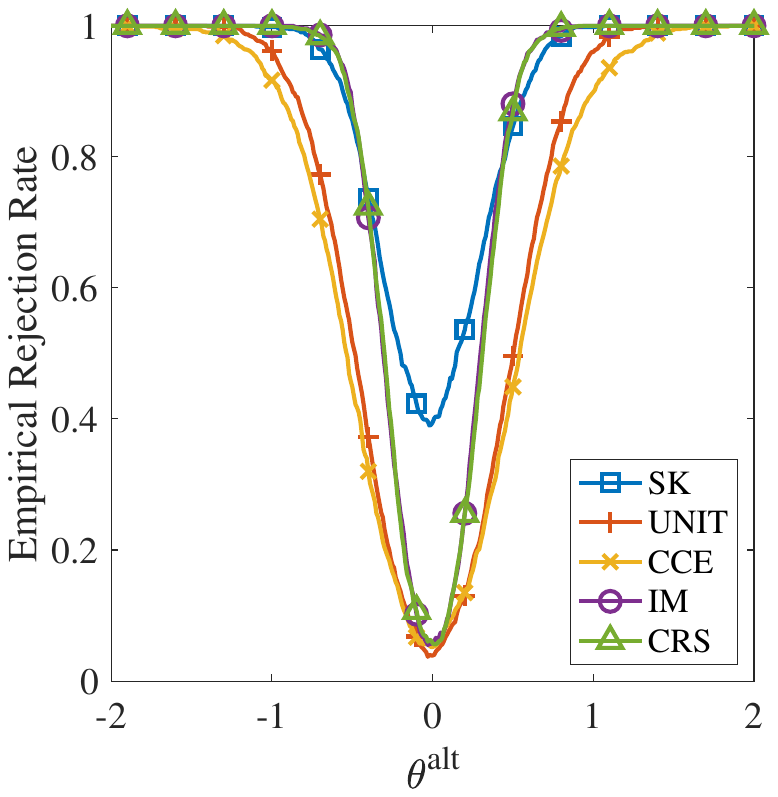}}\qquad
	\subfloat[SAR]{\includegraphics[width=.45\linewidth,trim=4.5cm 6.5cm 4.5cm 7cm]{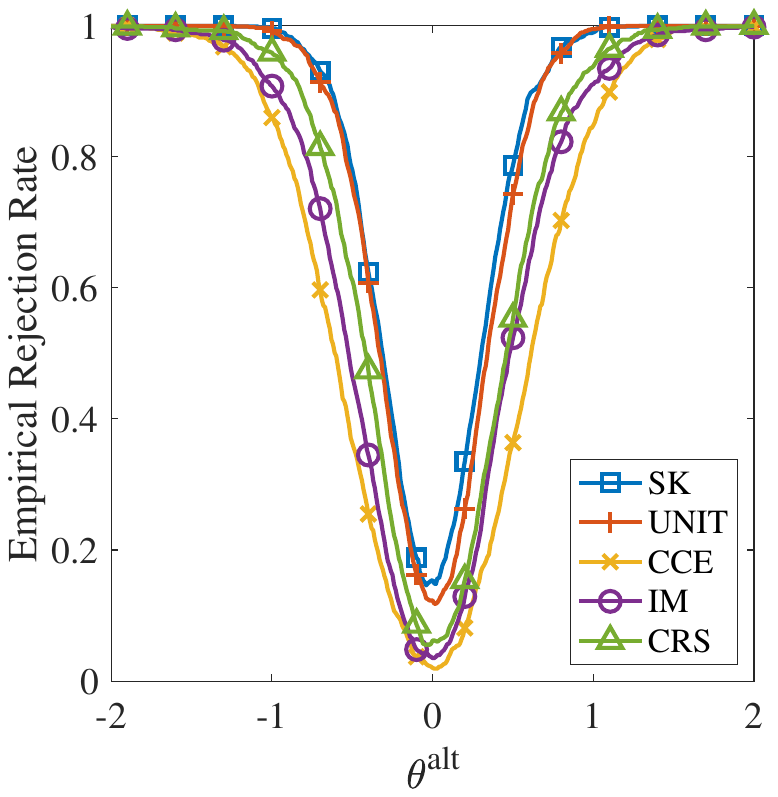}}\\
	\caption{OLS power curves ($N_{\text{pan}}=820$).}
\end{figure}

\begin{figure}[t]
	\centering
	\subfloat[BASELINE]{\includegraphics[width=.4\linewidth,trim=4.5cm 6.5cm 4.5cm 7cm]{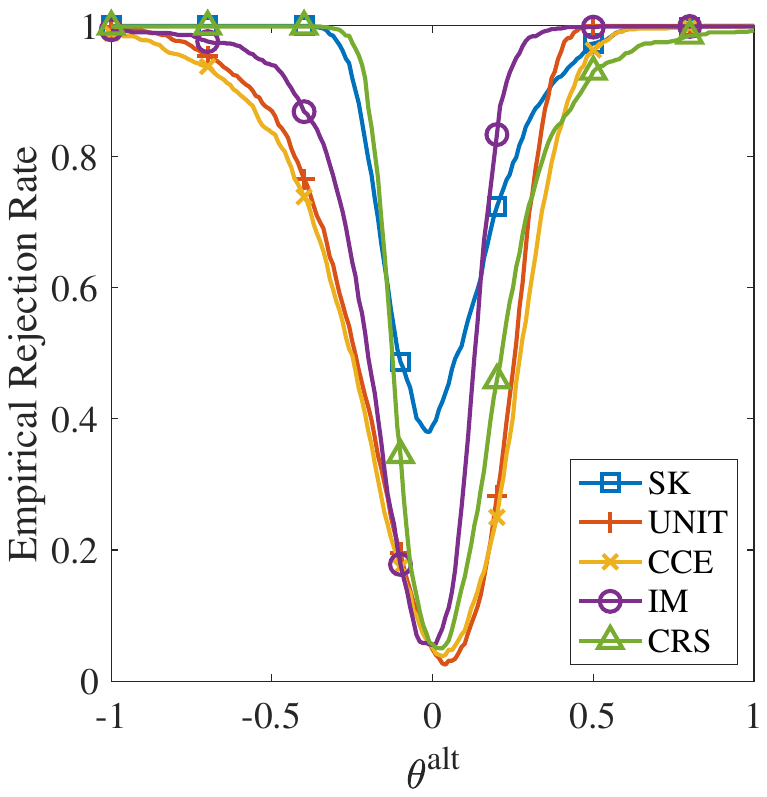}}\qquad
	\subfloat[SAR]{\includegraphics[width=.4\linewidth,trim=4.5cm 6.5cm 4.5cm 7cm]{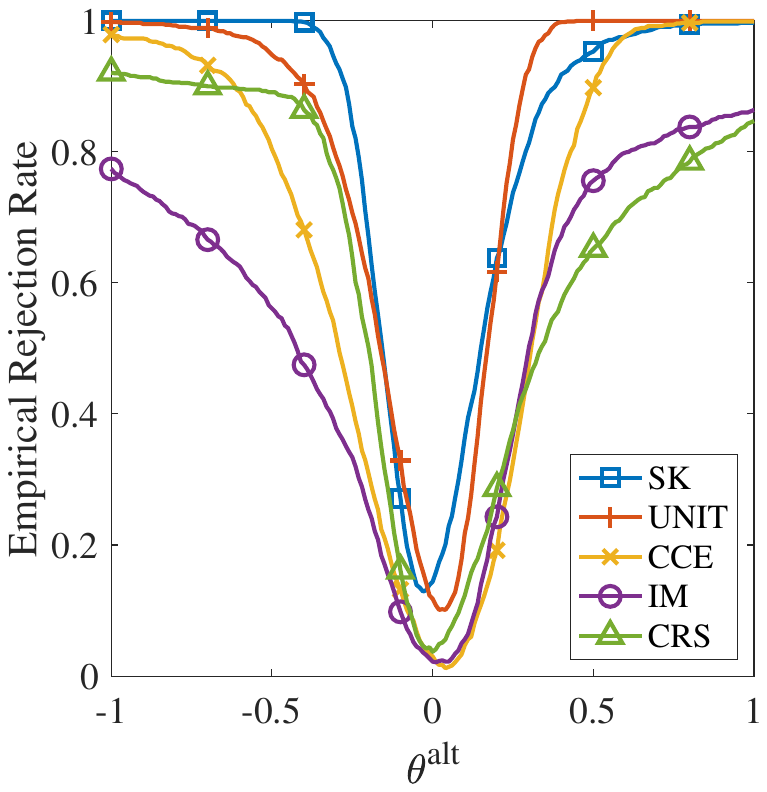}}\\
	\caption{IV power curves ($N_{\text{pan}}=820$).}
\end{figure}

\begin{landscape}

\begin{table}[H]
	\scriptsize
	\caption{Clustering: OLS - BASELINE ($N_{\mathrm{pan}}=820$)}
	\label{table_clustering_ols_baseline_large}
	\centering
	\renewcommand{\arraystretch}{1}
	\begin{threeparttable}	\begin{tabular}{cccccccccccccc}
			\hline
			\hline
			& & \multicolumn{11}{c}{$k$} & $\hat k$\\
			\cline{3-13}
			& & 2 & 3 & 4 & 5 & 6 & 7 & 8 & 9 & 10 & 11 & 12 \\
			\hline
CCE & size (usual cv)               & 0.056 & 0.093 & 0.117 & 0.146 & 0.180 & 0.178 & 0.177 & 0.188 & 0.209 & 0.209 & 0.218 &       \\
& size (simulated cv)           & 0.041 & 0.045 & 0.046 & 0.040 & 0.044 & 0.050 & 0.051 & 0.050 & 0.047 & 0.047 & 0.049 & 0.052 \\
& $\widehat{k}$ frequency                & 0.000 & 0.003 & 0.006 & 0.045 & 0.086 & 0.068 & 0.119 & 0.098 & 0.189 & 0.150 & 0.236 &       \\
& p(sim\_size\textgreater{}.05) & 0.994 & 1.000 & 1.000 & 1.000 & 1.000 & 1.000 & 1.000 & 1.000 & 1.000 & 1.000 & 1.000 & 1.000 \\
\cline{2-14}
& $\widehat{\alpha}$ quantile &       &       &       &       &       &       &       &       &       &       &       \\
& q10                           & 0.030 & 0.021 & 0.015 & 0.011 & 0.008 & 0.008 & 0.007 & 0.006 & 0.004 & 0.004 & 0.004 & 0.005 \\
& q25                           & 0.033 & 0.023 & 0.017 & 0.011 & 0.009 & 0.009 & 0.008 & 0.007 & 0.005 & 0.005 & 0.004 & 0.006 \\
& q50                           & 0.036 & 0.025 & 0.019 & 0.013 & 0.010 & 0.011 & 0.010 & 0.008 & 0.006 & 0.006 & 0.005 & 0.008 \\
& q75                           & 0.039 & 0.028 & 0.021 & 0.014 & 0.011 & 0.012 & 0.011 & 0.009 & 0.007 & 0.007 & 0.006 & 0.010 \\
& q90                           & 0.042 & 0.030 & 0.023 & 0.016 & 0.012 & 0.013 & 0.012 & 0.010 & 0.009 & 0.008 & 0.007 & 0.013 \\
\hline
IM  & size (usual cv)               & 0.046 & 0.060 & 0.059 & 0.064 & 0.064 & 0.075 & 0.085 & 0.091 & 0.096 & 0.090 & 0.076 &       \\
& size (simulated cv)           & 0.041 & 0.053 & 0.054 & 0.057 & 0.052 & 0.062 & 0.069 & 0.057 & 0.067 & 0.059 & 0.053 & 0.058 \\
& $\widehat{k}$ frequency                & 0.000 & 0.000 & 0.000 & 0.000 & 0.000 & 0.002 & 0.013 & 0.013 & 0.005 & 0.195 & 0.772 &       \\
& p(sim\_size\textgreater{}.05) & 0.492 & 0.565 & 0.727 & 0.675 & 0.828 & 0.966 & 0.966 & 0.998 & 1.000 & 1.000 & 0.999 & 0.997 \\
\cline{2-14}
& $\widehat{\alpha}$ quantile &       &       &       &       &       &       &       &       &       &       &       \\
& q10                           & 0.042 & 0.041 & 0.039 & 0.039 & 0.037 & 0.032 & 0.031 & 0.025 & 0.021 & 0.024 & 0.024 & 0.026 \\
& q25                           & 0.046 & 0.045 & 0.042 & 0.043 & 0.039 & 0.036 & 0.035 & 0.028 & 0.024 & 0.027 & 0.027 & 0.029 \\
& q50                           & 0.050 & 0.049 & 0.046 & 0.047 & 0.044 & 0.039 & 0.038 & 0.032 & 0.027 & 0.030 & 0.031 & 0.032 \\
& q75                           & 0.050 & 0.050 & 0.050 & 0.050 & 0.048 & 0.043 & 0.042 & 0.035 & 0.031 & 0.034 & 0.035 & 0.036 \\
& q90                           & 0.050 & 0.050 & 0.050 & 0.050 & 0.050 & 0.047 & 0.046 & 0.039 & 0.035 & 0.038 & 0.038 & 0.040 \\
\hline
CRS & size (usual cv)               & 0.000 & 0.000 & 0.000 & 0.000 & 0.034 & 0.074 & 0.082 & 0.088 & 0.100 & 0.093 & 0.077 &       \\
& size (simulated cv)           & 0.000 & 0.000 & 0.000 & 0.000 & 0.033 & 0.047 & 0.062 & 0.053 & 0.058 & 0.059 & 0.053 & 0.058 \\
& $\widehat{k}$ frequency                & 0.000 & 0.000 & 0.000 & 0.000 & 0.000 & 0.001 & 0.004 & 0.009 & 0.004 & 0.170 & 0.812 &       \\
& p(sim\_size\textgreater{}.05) & 0.000 & 0.000 & 0.000 & 0.000 & 0.027 & 0.882 & 0.923 & 0.993 & 1.000 & 0.998 & 0.999 & 0.992 \\
		\cline{2-14}
& $\widehat{\alpha}$ quantile &       &       &       &       &       &       &       &       &       &       &       \\
& q10                           & 0.050 & 0.050 & 0.050 & 0.050 & 0.050 & 0.047 & 0.031 & 0.027 & 0.021 & 0.023 & 0.023 & 0.024 \\
& q25                           & 0.050 & 0.050 & 0.050 & 0.050 & 0.050 & 0.047 & 0.039 & 0.029 & 0.023 & 0.026 & 0.026 & 0.027 \\
& q50                           & 0.050 & 0.050 & 0.050 & 0.050 & 0.050 & 0.047 & 0.039 & 0.033 & 0.027 & 0.030 & 0.030 & 0.031 \\
& q75                           & 0.050 & 0.050 & 0.050 & 0.050 & 0.050 & 0.047 & 0.047 & 0.035 & 0.031 & 0.034 & 0.034 & 0.035 \\
& q90                           & 0.050 & 0.050 & 0.050 & 0.050 & 0.050 & 0.050 & 0.047 & 0.039 & 0.035 & 0.038 & 0.038 & 0.039 \\
			\hline
		\end{tabular}
		\begin{tablenotes}
			\item Notes: Simulation results for the design described in Section \ref{simulation}. 
Inferential properties are presented for the estimators, CCE, IM, CRS, described in the text. 
Columns under $k$ report results with the number of groups fixed at a certain $k$.
$\hat k$ is the number of clusters chosen by the criterion on size-power tradeoff described in the text. 
Sizes under both usual critical value (``size (usual cv)'') and adjusted critical value (``size (simulated cv)'') are reported. 
The rows ``$\hat k$ frequency'' is the frequency of a particular $k$ achieving the highest simulated power among candidate $k$'s in the setting. 
The rows ``p(sim\_size\textgreater.05)'' report the empirical frequency of $\widehat {\mathrm{Err}}_{\mathrm{Type-I}}(\bullet(.05), k)>.05$.
The row ``$\widehat{\alpha}$ quantile'' report the quantiles of the selected number of groups. 
		\end{tablenotes}
	\end{threeparttable}
\end{table}
\end{landscape}

\begin{landscape}

\begin{table}[H]
	\scriptsize
	\caption{Clustering: OLS - SAR ($N_{\mathrm{pan}}=820$)}
	\label{table_clustering_ols_sar_large}
	\centering
	\begin{threeparttable}	\begin{tabular}{cccccccccccccc}
			\hline
			\hline
			& & \multicolumn{11}{c}{$k$} & $\hat k$\\
\cline{3-13}
& & 2 & 3 & 4 & 5 & 6 & 7 & 8 & 9 & 10 & 11 & 12 \\
			\hline 
CCE & size (usual cv)               & 0.040 & 0.043 & 0.037 & 0.033 & 0.033 & 0.035 & 0.035 & 0.045 & 0.047 & 0.048 & 0.052 &       \\
& size (simulated cv)           & 0.036 & 0.037 & 0.026 & 0.020 & 0.021 & 0.017 & 0.012 & 0.017 & 0.019 & 0.020 & 0.020 & 0.020 \\
& $\widehat{k}$ frequency                & 0.000 & 0.000 & 0.004 & 0.008 & 0.028 & 0.052 & 0.105 & 0.113 & 0.167 & 0.220 & 0.303 &       \\
& p(sim\_size\textgreater{}.05) & 0.767 & 0.987 & 0.992 & 1.000 & 1.000 & 1.000 & 1.000 & 1.000 & 1.000 & 1.000 & 1.000 & 1.000 \\
\cline{2-14}
& $\widehat{\alpha}$ quantile &       &       &       &       &       &       &       &       &       &       &       \\
& q10                           & 0.038 & 0.030 & 0.028 & 0.024 & 0.021 & 0.023 & 0.022 & 0.020 & 0.019 & 0.018 & 0.016 & 0.020 \\
& q25                           & 0.041 & 0.033 & 0.031 & 0.027 & 0.024 & 0.025 & 0.024 & 0.022 & 0.021 & 0.020 & 0.019 & 0.022 \\
& q50                           & 0.045 & 0.036 & 0.035 & 0.030 & 0.027 & 0.028 & 0.027 & 0.025 & 0.024 & 0.023 & 0.021 & 0.026 \\
& q75                           & 0.050 & 0.040 & 0.038 & 0.033 & 0.030 & 0.031 & 0.030 & 0.028 & 0.027 & 0.026 & 0.024 & 0.029 \\
& q90                           & 0.050 & 0.044 & 0.041 & 0.036 & 0.033 & 0.034 & 0.033 & 0.031 & 0.029 & 0.028 & 0.026 & 0.033 \\
\hline
IM  & size (usual cv)               & 0.051 & 0.058 & 0.027 & 0.024 & 0.027 & 0.037 & 0.040 & 0.039 & 0.034 & 0.034 & 0.042 &       \\
& size (simulated cv)           & 0.047 & 0.055 & 0.026 & 0.020 & 0.024 & 0.034 & 0.037 & 0.036 & 0.033 & 0.032 & 0.040 & 0.035 \\
& $\widehat{k}$ frequency                & 0.000 & 0.000 & 0.002 & 0.006 & 0.073 & 0.044 & 0.409 & 0.223 & 0.061 & 0.067 & 0.115 &       \\
& p(sim\_size\textgreater{}.05) & 0.495 & 0.366 & 0.515 & 0.383 & 0.471 & 0.522 & 0.535 & 0.491 & 0.453 & 0.482 & 0.483 & 0.354 \\
\cline{2-14}
& $\widehat{\alpha}$ quantile &       &       &       &       &       &       &       &       &       &       &       \\
& q10                           & 0.042 & 0.043 & 0.041 & 0.044 & 0.042 & 0.041 & 0.041 & 0.042 & 0.043 & 0.042 & 0.042 & 0.047 \\
& q25                           & 0.045 & 0.047 & 0.045 & 0.048 & 0.046 & 0.045 & 0.045 & 0.046 & 0.046 & 0.046 & 0.046 & 0.049 \\
& q50                           & 0.050 & 0.050 & 0.050 & 0.050 & 0.050 & 0.049 & 0.049 & 0.050 & 0.050 & 0.050 & 0.050 & 0.050 \\
& q75                           & 0.050 & 0.050 & 0.050 & 0.050 & 0.050 & 0.050 & 0.050 & 0.050 & 0.050 & 0.050 & 0.050 & 0.050 \\
& q90                           & 0.050 & 0.050 & 0.050 & 0.050 & 0.050 & 0.050 & 0.050 & 0.050 & 0.050 & 0.050 & 0.050 & 0.050 \\
\hline
CRS & size (usual cv)               & 0.000 & 0.000 & 0.000 & 0.000 & 0.022 & 0.053 & 0.050 & 0.054 & 0.048 & 0.047 & 0.058 &       \\
& size (simulated cv)           & 0.000 & 0.000 & 0.000 & 0.000 & 0.022 & 0.049 & 0.048 & 0.054 & 0.046 & 0.040 & 0.053 & 0.062 \\
& $\widehat{k}$ frequency                & 0.000 & 0.000 & 0.000 & 0.000 & 0.000 & 0.019 & 0.199 & 0.143 & 0.127 & 0.188 & 0.324 &       \\
& p(sim\_size\textgreater{}.05) & 0.000 & 0.000 & 0.000 & 0.000 & 0.001 & 0.355 & 0.342 & 0.383 & 0.487 & 0.515 & 0.522 & 0.277 \\
\cline{2-14}
& $\widehat{\alpha}$ quantile &       &       &       &       &       &       &       &       &       &       &       \\
& q10                           & 0.050 & 0.050 & 0.050 & 0.050 & 0.050 & 0.047 & 0.047 & 0.043 & 0.042 & 0.042 & 0.041 & 0.047 \\
& q25                           & 0.050 & 0.050 & 0.050 & 0.050 & 0.050 & 0.047 & 0.047 & 0.047 & 0.046 & 0.045 & 0.045 & 0.050 \\
& q50                           & 0.050 & 0.050 & 0.050 & 0.050 & 0.050 & 0.050 & 0.050 & 0.050 & 0.050 & 0.050 & 0.050 & 0.050 \\
& q75                           & 0.050 & 0.050 & 0.050 & 0.050 & 0.050 & 0.050 & 0.050 & 0.050 & 0.050 & 0.050 & 0.050 & 0.050 \\
& q90                           & 0.050 & 0.050 & 0.050 & 0.050 & 0.050 & 0.050 & 0.050 & 0.050 & 0.050 & 0.050 & 0.050 & 0.050\\
			\hline
		\end{tabular}
		\begin{tablenotes}
			\item Notes: Simulation results for the design described in Section \ref{simulation}. 
Inferential properties are presented for the estimators, CCE, IM, CRS, described in the text. 
Columns under $k$ report results with the number of groups fixed at a certain $k$.
$\hat k$ is the number of clusters chosen by the criterion on size-power tradeoff described in the text. 
Sizes under both usual critical value (``size (usual cv)'') and adjusted critical value (``size (simulated cv)'') are reported. 
The rows ``$\hat k$ frequency'' is the frequency of a particular $k$ achieving the highest simulated power among candidate $k$'s in the setting. 
The rows ``p(sim\_size\textgreater.05)'' report the empirical frequency of $\widehat {\mathrm{Err}}_{\mathrm{Type-I}}(\bullet(.05), k)>.05$.
The row ``$\widehat{\alpha}$ quantile'' report the quantiles of the selected number of groups. 
		\end{tablenotes}
	\end{threeparttable}
\end{table}
\end{landscape}

\begin{landscape}

\begin{table}[H]
	\scriptsize
	\caption{Clustering: IV - BASELINE ($N_{\mathrm{pan}}=820$)}
	\label{table_clustering_iv_baseline_large}
	\centering
	\begin{threeparttable}	\begin{tabular}{cccccccccccccc}
			\hline
			\hline
			& & \multicolumn{11}{c}{$k$} & $\hat k$\\
\cline{3-13}
& & 2 & 3 & 4 & 5 & 6 & 7 & 8 & 9 & 10 & 11 & 12 \\
			\hline
CCE & size (usual cv)               & 0.059 & 0.091 & 0.113 & 0.150 & 0.171 & 0.156 & 0.169 & 0.186 & 0.201 & 0.197 & 0.210 &       \\
& size (simulated cv)           & 0.044 & 0.047 & 0.048 & 0.043 & 0.046 & 0.047 & 0.054 & 0.052 & 0.047 & 0.040 & 0.041 & 0.051 \\
& $\widehat{k}$ frequency                & 0.000 & 0.002 & 0.008 & 0.009 & 0.053 & 0.038 & 0.069 & 0.096 & 0.185 & 0.207 & 0.333 &       \\
& p(sim\_size\textgreater{}.05) & 0.994 & 1.000 & 1.000 & 1.000 & 1.000 & 1.000 & 1.000 & 1.000 & 1.000 & 1.000 & 1.000 & 1.000 \\
\cline{2-14}
& $\widehat{\alpha}$ quantile &       &       &       &       &       &       &       &       &       &       &       \\
& q10                           & 0.030 & 0.021 & 0.015 & 0.010 & 0.007 & 0.008 & 0.007 & 0.005 & 0.004 & 0.004 & 0.003 & 0.004 \\
& q25                           & 0.033 & 0.023 & 0.017 & 0.011 & 0.008 & 0.009 & 0.008 & 0.006 & 0.005 & 0.005 & 0.004 & 0.005 \\
& q50                           & 0.036 & 0.025 & 0.019 & 0.012 & 0.009 & 0.010 & 0.009 & 0.008 & 0.006 & 0.006 & 0.005 & 0.007 \\
& q75                           & 0.039 & 0.028 & 0.021 & 0.014 & 0.011 & 0.012 & 0.011 & 0.009 & 0.007 & 0.007 & 0.006 & 0.009 \\
& q90                           & 0.042 & 0.030 & 0.023 & 0.016 & 0.012 & 0.014 & 0.012 & 0.010 & 0.009 & 0.009 & 0.007 & 0.011 \\
\hline
IM  & size (usual cv)               & 0.050 & 0.054 & 0.058 & 0.064 & 0.062 & 0.069 & 0.081 & 0.087 & 0.093 & 0.079 & 0.074 &       \\
& size (simulated cv)           & 0.048 & 0.051 & 0.057 & 0.059 & 0.054 & 0.056 & 0.073 & 0.063 & 0.065 & 0.057 & 0.054 & 0.055 \\
& $\widehat{k}$ frequency                & 0.000 & 0.000 & 0.000 & 0.000 & 0.000 & 0.000 & 0.001 & 0.009 & 0.000 & 0.128 & 0.862 &       \\
& p(sim\_size\textgreater{}.05) & 0.483 & 0.544 & 0.636 & 0.528 & 0.694 & 0.890 & 0.903 & 0.988 & 0.998 & 0.993 & 0.998 & 0.993 \\
\cline{2-14}
& $\widehat{\alpha}$ quantile &       &       &       &       &       &       &       &       &       &       &       \\
& q10                           & 0.042 & 0.042 & 0.039 & 0.041 & 0.039 & 0.035 & 0.035 & 0.028 & 0.025 & 0.028 & 0.028 & 0.028 \\
& q25                           & 0.046 & 0.045 & 0.043 & 0.045 & 0.043 & 0.039 & 0.038 & 0.032 & 0.028 & 0.031 & 0.031 & 0.031 \\
& q50                           & 0.050 & 0.049 & 0.048 & 0.049 & 0.047 & 0.043 & 0.042 & 0.035 & 0.032 & 0.034 & 0.034 & 0.035 \\
& q75                           & 0.050 & 0.050 & 0.050 & 0.050 & 0.050 & 0.046 & 0.046 & 0.039 & 0.035 & 0.038 & 0.038 & 0.039 \\
& q90                           & 0.050 & 0.050 & 0.050 & 0.050 & 0.050 & 0.050 & 0.050 & 0.043 & 0.039 & 0.042 & 0.041 & 0.042 \\
\hline
CRS & size (usual cv)               & 0.000 & 0.000 & 0.000 & 0.000 & 0.034 & 0.072 & 0.084 & 0.087 & 0.107 & 0.098 & 0.089 &       \\
& size (simulated cv)           & 0.000 & 0.000 & 0.000 & 0.000 & 0.033 & 0.043 & 0.060 & 0.058 & 0.060 & 0.057 & 0.053 & 0.056 \\
& $\widehat{k}$ frequency                & 0.000 & 0.000 & 0.000 & 0.000 & 0.000 & 0.001 & 0.001 & 0.008 & 0.001 & 0.125 & 0.864 &       \\
& p(sim\_size\textgreater{}.05) & 0.000 & 0.000 & 0.000 & 0.000 & 0.030 & 0.889 & 0.941 & 0.996 & 1.000 & 1.000 & 1.000 & 0.996 \\
\cline{2-14}
& $\widehat{\alpha}$ quantile &       &       &       &       &       &       &       &       &       &       &       \\
& q10                           & 0.050 & 0.050 & 0.050 & 0.050 & 0.050 & 0.047 & 0.031 & 0.027 & 0.021 & 0.022 & 0.022 & 0.022 \\
& q25                           & 0.050 & 0.050 & 0.050 & 0.050 & 0.050 & 0.047 & 0.039 & 0.029 & 0.023 & 0.025 & 0.024 & 0.025 \\
& q50                           & 0.050 & 0.050 & 0.050 & 0.050 & 0.050 & 0.047 & 0.039 & 0.031 & 0.026 & 0.028 & 0.028 & 0.028 \\
& q75                           & 0.050 & 0.050 & 0.050 & 0.050 & 0.050 & 0.047 & 0.047 & 0.035 & 0.029 & 0.031 & 0.031 & 0.032 \\
& q90                           & 0.050 & 0.050 & 0.050 & 0.050 & 0.050 & 0.050 & 0.047 & 0.039 & 0.033 & 0.035 & 0.034 & 0.036 \\
			\hline
		\end{tabular}
		\begin{tablenotes}
			\item Notes: Simulation results for the design described in Section \ref{simulation}. 
Inferential properties are presented for the estimators, CCE, IM, CRS, described in the text. 
Columns under $k$ report results with the number of groups fixed at a certain $k$.
$\hat k$ is the number of clusters chosen by the criterion on size-power tradeoff described in the text. 
Sizes under both usual critical value (``size (usual cv)'') and adjusted critical value (``size (simulated cv)'') are reported. 
The rows ``$\hat k$ frequency'' is the frequency of a particular $k$ achieving the highest simulated power among candidate $k$'s in the setting. 
The rows ``p(sim\_size\textgreater.05)'' report the empirical frequency of $\widehat {\mathrm{Err}}_{\mathrm{Type-I}}(\bullet(.05), k)>.05$.
The row ``$\widehat{\alpha}$ quantile'' report the quantiles of the selected number of groups. 
		\end{tablenotes}
	\end{threeparttable}
\end{table}
\end{landscape}

\begin{landscape}

\begin{table}[H]
	\scriptsize
	\caption{Clustering: IV - SAR ($N_{\mathrm{pan}}=820$)}
	\label{table_clustering_iv_sar_large}
	\centering
	\begin{threeparttable}	\begin{tabular}{cccccccccccccc}
			\hline
			\hline
			& & \multicolumn{11}{c}{$k$} & $\hat k$\\
\cline{3-13}
& & 2 & 3 & 4 & 5 & 6 & 7 & 8 & 9 & 10 & 11 & 12 \\
			\hline
CCE & size (usual cv)               & 0.034 & 0.042 & 0.033 & 0.035 & 0.041 & 0.030 & 0.040 & 0.041 & 0.041 & 0.043 & 0.047 &       \\
& size (simulated cv)           & 0.029 & 0.035 & 0.021 & 0.021 & 0.019 & 0.020 & 0.022 & 0.025 & 0.027 & 0.029 & 0.030 & 0.030 \\
& $\widehat{k}$ frequency                & 0.000 & 0.000 & 0.000 & 0.000 & 0.007 & 0.016 & 0.046 & 0.076 & 0.151 & 0.268 & 0.436 &       \\
& p(sim\_size\textgreater{}.05) & 0.800 & 0.992 & 0.994 & 1.000 & 1.000 & 0.999 & 1.000 & 1.000 & 1.000 & 1.000 & 1.000 & 1.000 \\
\cline{2-14}
& $\widehat{\alpha}$ quantile &       &       &       &       &       &       &       &       &       &       &       \\
& q10                           & 0.037 & 0.030 & 0.028 & 0.024 & 0.022 & 0.022 & 0.022 & 0.020 & 0.019 & 0.018 & 0.016 & 0.019 \\
& q25                           & 0.041 & 0.033 & 0.031 & 0.026 & 0.024 & 0.025 & 0.024 & 0.022 & 0.021 & 0.020 & 0.019 & 0.021 \\
& q50                           & 0.044 & 0.036 & 0.035 & 0.029 & 0.027 & 0.028 & 0.027 & 0.025 & 0.024 & 0.023 & 0.021 & 0.024 \\
& q75                           & 0.049 & 0.040 & 0.039 & 0.033 & 0.030 & 0.031 & 0.031 & 0.028 & 0.027 & 0.026 & 0.024 & 0.027 \\
& q90                           & 0.050 & 0.043 & 0.042 & 0.036 & 0.033 & 0.034 & 0.034 & 0.031 & 0.030 & 0.029 & 0.027 & 0.031 \\
\hline
IM  & size (usual cv)               & 0.041 & 0.056 & 0.029 & 0.021 & 0.020 & 0.037 & 0.030 & 0.028 & 0.026 & 0.019 & 0.026 &       \\
& size (simulated cv)           & 0.039 & 0.055 & 0.025 & 0.019 & 0.019 & 0.035 & 0.029 & 0.028 & 0.026 & 0.019 & 0.026 & 0.022 \\
& $\widehat{k}$ frequency                & 0.000 & 0.004 & 0.008 & 0.008 & 0.153 & 0.003 & 0.777 & 0.046 & 0.001 & 0.000 & 0.000 &       \\
& p(sim\_size\textgreater{}.05) & 0.481 & 0.321 & 0.399 & 0.234 & 0.256 & 0.224 & 0.240 & 0.176 & 0.115 & 0.104 & 0.139 & 0.160 \\
\cline{2-14}
& $\widehat{\alpha}$ quantile &       &       &       &       &       &       &       &       &       &       &       \\
& q10                           & 0.042 & 0.045 & 0.043 & 0.046 & 0.046 & 0.047 & 0.046 & 0.048 & 0.049 & 0.050 & 0.049 & 0.049 \\
& q25                           & 0.046 & 0.048 & 0.047 & 0.050 & 0.050 & 0.050 & 0.050 & 0.050 & 0.050 & 0.050 & 0.050 & 0.050 \\
& q50                           & 0.050 & 0.050 & 0.050 & 0.050 & 0.050 & 0.050 & 0.050 & 0.050 & 0.050 & 0.050 & 0.050 & 0.050 \\
& q75                           & 0.050 & 0.050 & 0.050 & 0.050 & 0.050 & 0.050 & 0.050 & 0.050 & 0.050 & 0.050 & 0.050 & 0.050 \\
& q90                           & 0.050 & 0.050 & 0.050 & 0.050 & 0.050 & 0.050 & 0.050 & 0.050 & 0.050 & 0.050 & 0.050 & 0.050 \\
\hline
CRS & size (usual cv)               & 0.000 & 0.000 & 0.000 & 0.000 & 0.025 & 0.054 & 0.051 & 0.045 & 0.042 & 0.039 & 0.052 &       \\
& size (simulated cv)           & 0.000 & 0.000 & 0.000 & 0.000 & 0.025 & 0.044 & 0.041 & 0.042 & 0.040 & 0.032 & 0.038 & 0.037 \\
& $\widehat{k}$ frequency                & 0.000 & 0.000 & 0.000 & 0.000 & 0.000 & 0.031 & 0.469 & 0.320 & 0.090 & 0.048 & 0.042 &       \\
& p(sim\_size\textgreater{}.05) & 0.000 & 0.000 & 0.000 & 0.000 & 0.002 & 0.389 & 0.446 & 0.607 & 0.780 & 0.899 & 0.957 & 0.352 \\
\cline{2-14}
& $\widehat{\alpha}$ quantile &       &       &       &       &       &       &       &       &       &       &       \\
& q10                           & 0.050 & 0.050 & 0.050 & 0.050 & 0.050 & 0.047 & 0.047 & 0.039 & 0.038 & 0.036 & 0.034 & 0.045 \\
& q25                           & 0.050 & 0.050 & 0.050 & 0.050 & 0.050 & 0.047 & 0.047 & 0.043 & 0.041 & 0.039 & 0.036 & 0.047 \\
& q50                           & 0.050 & 0.050 & 0.050 & 0.050 & 0.050 & 0.050 & 0.050 & 0.047 & 0.045 & 0.042 & 0.040 & 0.050 \\
& q75                           & 0.050 & 0.050 & 0.050 & 0.050 & 0.050 & 0.050 & 0.050 & 0.050 & 0.050 & 0.046 & 0.043 & 0.050 \\
& q90                           & 0.050 & 0.050 & 0.050 & 0.050 & 0.050 & 0.050 & 0.050 & 0.050 & 0.050 & 0.050 & 0.047 & 0.050\\
			\hline
		\end{tabular}
		\begin{tablenotes}
			\item Notes: Simulation results for the design described in Section \ref{simulation}. 
Inferential properties are presented for the estimators, CCE, IM, CRS, described in the text. 
Columns under $k$ report results with the number of groups fixed at a certain $k$.
$\hat k$ is the number of clusters chosen by the criterion on size-power tradeoff described in the text. 
Sizes under both usual critical value (``size (usual cv)'') and adjusted critical value (``size (simulated cv)'') are reported. 
The rows ``$\hat k$ frequency'' is the frequency of a particular $k$ achieving the highest simulated power among candidate $k$'s in the setting. 
The rows ``p(sim\_size\textgreater.05)'' report the empirical frequency of $\widehat {\mathrm{Err}}_{\mathrm{Type-I}}(\bullet(.05), k)>.05$.
The row ``$\widehat{\alpha}$ quantile'' report the quantiles of the selected number of groups. 
		\end{tablenotes}
	\end{threeparttable}
\end{table}
\end{landscape}


\end{document}